\documentclass[reqno, 12pt]{amsart}
\usepackage{mathtools}
\usepackage[T1]{fontenc}
\usepackage{amsmath, amsthm, amssymb}
\usepackage[ansinew]{inputenc}
\usepackage{xcolor}
\usepackage{bbm}
\usepackage{mathrsfs}
\usepackage{cite}
\usepackage{xparse}
\usepackage{fullpage}
\usepackage{pdflscape}
\usepackage{stmaryrd}
\usepackage{slashed}
\usepackage{multirow}
\usepackage{lscape}
\usepackage{dsfont}
\usepackage{amsfonts}
\usepackage{scrextend} 
\usepackage{longtable}
\usepackage{verbatim}
\usepackage{subfig}
\usepackage{yfonts}
\usepackage{enumitem} 
\usepackage{tabu}
\usepackage[colorlinks=true,allcolors=red]{hyperref}
\usepackage{graphicx} 
\usepackage{accents}

\usepackage{array}
\newcolumntype{P}[1]{>{\centering\arraybackslash}p{#1}}
\newcolumntype{M}[1]{>{\centering\arraybackslash}m{#1}}

\usepackage[left=2.1cm, right=2.1cm]{geometry}

\DeclareSymbolFont{sfletters}{OML}{cmbrm}{m}{it}
\DeclareMathSymbol{\somega}{\mathord}{sfletters}{"21}
\DeclareSymbolFont{Letters} {U}{zeur}{m}{n} 
\DeclareMathSymbol{\ssomega}{\mathalpha}{Letters}{"21} 
\DeclareSymbolFont{Letters} {U}{zeur}{m}{n} 
\DeclareMathSymbol{\pPhi}{\mathalpha}{Letters}{"1E}

\makeatletter
\newcommand*\bigcdot{\mathpalette\bigcdot@{.5}}
\newcommand*\bigcdot@[2]{\mathbin{\vcenter{\hbox{\scalebox{#2}{$\m@th#1\bullet$}}}}}
\makeatother

\newtheorem{theorem}{Theorem}[section]
\newtheorem{lemma}[theorem]{Lemma}
\newtheorem{prop}[theorem]{Proposition}

\theoremstyle{definition}

\theoremstyle{remark}
\newtheorem{remark}[theorem]{Remark}

\allowdisplaybreaks[4]

\numberwithin{equation}{section}

\usepackage{scalerel,stackengine}
\stackMath
\newcommand\reallywidehat[1]{
\savestack{\tmpbox}{\stretchto{
  \scaleto{
    \scalerel*[\widthof{\ensuremath{#1}}]{\kern.1pt\mathchar"0362\kern.1pt}
    {\rule{0ex}{\textheight}}
  }{\textheight}
}{2.4ex}}
\stackon[-6.9pt]{#1}{\tmpbox}
}
\makeatletter
\newcommand*{\rom}[1]{\expandafter\@slowromancap\romannumeral #1@}
\makeatother

\begin{document}

\title{On the Fourier analysis of the Einstein-Klein-Gordon system: Limiting behaviour of the Fourier coefficients}

%  Author 1 information
\author{Athanasios Chatzikaleas}
\address{Department of Mathematics, National and Kapodistrian University of Athens, Greece}
\email{achatzik@uoa.gr}
%  Author 2 information
\author{Jacques Smulevici}
\address{Laboratory Jacques-Louis Lions (LJLL), Sorbonne Universit\'e, 4 place Jussieu, 75252 Paris, France}
\email{jacques.smulevici@ljll.math.upmc.fr}

\dedicatory{} 

\begin{abstract}
We revisit the problem of computing the Fourier coefficients associated with the spherically symmetric   massless-scalar field system in $(1+3)$ dimensions initiated in \cite{2020arXiv200411049C}. We rigorously derive computationally efficient formulas for these coefficients and then use them to establish their asymptotic behavior as the mode numbers tend to infinity, including the explicit determination of the leading-order coefficients. Our approach builds on the method developed in \cite{chatzikaleas2023timeperiodicsolutionsnekhoroshev} and relies on two powerful algorithms for hypergeometric summation, due to Zeilberger and Gosper. These algorithms enable us to systematically derive first-order linear recurrence relations and obtain closed-form discrete antiderivatives respectively  for  intricate hypergeometric sequences that arise in our analysis. 
\end{abstract}

\maketitle
\tableofcontents 
\addtocontents{toc}{\protect\setcounter{tocdepth}{1}}

\section{Introduction}

\subsection{Instability of Anti de-Sitter spacetime}
 \noindent    
The Anti-de Sitter (AdS) spacetime is the unique (up to quotients by isometry subgroups) maximally symmetric solution to the Einstein equations with a negative cosmological constant, 
 \begin{align*}  
 		\text{Ric}(\text{g}) = \Lambda  \text{g} , \quad  \Lambda<0.
 \end{align*} 
Due to the fact that AdS is not globally hyperbolic, in order to determine of evolution of fields on AdS one has to prescribe, in addition to initial data on the $\{t=0\}$ slice, suitable boundary conditions at the conformal infinity. For reflective boundary conditions, there is no flux of energy across the conformal boundary as it behaves like a mirror at which perturbations propagating outwards bounce off and return to the interior of the spacetime. In this setting, Dafermos-Holzegel \cite{DafermosHolzegel,DafermosTalk} conjectured the instability of the AdS solution to the Einstein equations. The first who pioneered in this direction were Bizo\'{n}-Rostworowski \cite{PhysRevLett.107.031102} who numerically studied small perturbations of AdS for the Einstein-massless-scalar field system  
\begin{align}\label{EKGintro} 
 		\begin{dcases}
 		\text{Ric}(\text{g})  -\frac{1}{2}\text{R}(\text{g}) \text{g}    +  \Lambda  \text{g} = 8 \pi \left ( d  \phi \otimes d \phi - \frac{1}{2} |\partial \phi|_{\text{g}}^2    \text{g} \right ), \quad \Lambda<0,  \\
 			 \Box_{\text{g}}   \phi = 0,
 	\end{dcases}
 	\end{align} 
in spherical symmetry. On top of their numerical analysis, they also proposed a mechanism responsible for the instability. Using a perturbative series expansion, they suggested that energy is transferred from low to high frequencies through resonant interactions, which, in the setting of the Einstein equations, eventually become responsible for the instability leading to the formation of a small black hole. This important work then triggered a substantial amount of numerical and heuristic studies, such as the analysis of the resonant interactions outside from spherical symmetry in \cite{MR2978943}, or further numerical analysis of the instability such as  \cite{PhysRevLett.115.081103, PhysRevLett.114.071102, PhysRevLett.119.191103, PhysRevLett.111.041102, MR4016456,  MR3708611, MR3205859, MR3491620, MR3329398}. The resonant mode mixing mechanism of Bizo\'{n}-Rostworowski \cite{PhysRevLett.107.031102} is relevant only at the first stage of the instability where perturbation theory is still valid. A non-perturbative approach that lead to a rigorous proof of the AdS instability conjecture was developed by Moschidis for the Einstein--massless Vlasov system and for the Einstein--null dust system \cite{MR4150259, MR3828498, MoschidisNew, Moschidis2020NullDust, 2017arXiv170408685M}. Instead of looking at the Fourier space and the mode mixing, Moschidis prescribed a large number of pulses in the physical space as initial data where each one initially moves inwards towards the center of symmetry, reemerges as an outer beam and bounces back into the interior of the spacetime once it reaches infinity (due to the reflective boundary conditions) leading to complicated non-linear interactions and overlappings. In fact, Moschidis showed that energy flows from one beam to the other in such a way so that a trapped surface is eventually formed forcing the arbitrary small perturbations of the AdS initial data to necessarily collapse into a black hole.  \\ \\
Although AdS is unstable for generic perturbations, there should also exist ``special'' initial data that evolve into time-periodic solutions that stay close to AdS for all times. This is a fascinating improvement of the instability conjecture put forward by Rostworowski-Maliborski \cite{PhysRevLett.111.051102} who  numerically constructed time-periodic solutions to the Einstein-massless-scalar field equations in spherical symmetry bifurcating from initial data dominated by 1-modes. Their analysis also suggests that the time-periodic solutions should in fact enjoy better stability properties than the AdS spacetime itself. Following their work, a large class of time-periodic solutions that extend the single-mode ones constructed in \cite{PhysRevLett.107.031102, PhysRevLett.111.051102} was determined in \cite{PhysRevLett.113.071601} and their stability was analyzed in \cite{PhysRevD.92.084001}. Following \cite{PhysRevLett.107.031102}, Rostworowski-Maliborski \cite{PhysRevLett.111.051102} constructed their time-periodic solutions using a perturbative series expansion by expanding the dynamical variables in \eqref{EKGintro} in terms of the eigenfunctions $\{e_n: n \geq 0\}$ to the associated linearized operator. Then, the non-linearities of the system define mode couplings via oscillatory integrals, the so-called Fourier coefficients, which are essentially integrals of products of the eigenfunctions, see for example \cite{2020arXiv200411049C,PhysRevLett.107.031102,PhysRevLett.111.051102}. The initial data $\phi(0,\cdot)$ for the wave function, 
\begin{align}\label{PerturbativeSeriesInitialData}
	\phi(0,\cdot)=\phi_1(0,\cdot) \epsilon +\sum_{\lambda=1}^{\infty}\phi_{\lambda}(0,\cdot) \epsilon^{2\lambda+1},
\end{align}
were also expanded into modes. After projecting onto the set of all eigenmodes, one obtains an infinite system of harmonic oscillators in which secular\footnote{These are terms of the form $\tau^{\alpha} \sin( \beta \tau)$ and $\tau^{\alpha} \cos( \beta \tau)$ for some $\alpha>0$ and $\beta \in \mathbb{R}$.} terms appear naturally and spoil the periodicity of the solutions. According to \cite{PhysRevLett.107.031102}, if these terms are not removed, they should eventually become responsible for exponential instability and a possible breakdown of the solutions. Rostworowski-Maliborski \cite{PhysRevLett.111.051102} were able to fine tune the initial data $\phi(0,\cdot)$ so that the resulting infinite system of harmonic oscillators described above is free of secular terms. To this end, they discovered that such a ``special'' choice of the initial data follows from a recurrence relation that yields order-by-order the coefficients $\{\phi _{\lambda}(0,\cdot): \lambda \geq 1 \}$ in the perturbative expansion \eqref{PerturbativeSeriesInitialData}. In order to initiate their recursive procedure, they chose a fixed 1-mode dominant leading order term in the expansion \eqref{PerturbativeSeriesInitialData}, namely 
\begin{align*}
	\phi_1(0,\cdot)=e_{\gamma},
\end{align*}
for some integer $\gamma \geq 0$.  However, given $\phi _{\lambda}(0,\cdot)$, the recurrence relation above yields the next order coefficient $\phi _{\lambda+1}(0,\cdot)$, only if the so called  non-degeneracy condition holds true, that is a system of infinitely many (linear and non-linear) conditions for the Fourier coefficients. Despite the considerable computational effort involved, Rostworowski-Maliborski \cite{PhysRevLett.111.051102} numerically computed the Fourier coefficients so that they can verify the non-degeneracy condition up to some high order. This way, they were able to extract step-by-step the coefficients of the ``special'' initial data only up to a high order $2\lambda+1$ that was satisfactory enough so that they can safely assume that the recurrence relation keeps on going. \\ \\
Despite the fact that the method of Rostworowski-Maliborski \cite{PhysRevLett.111.051102} relies on the numerical computation of the Fourier coefficients, it surely provides strong numerical evidence indicating the existence of time-periodic solutions for the spherically symmetric Einstein-Klein-Gordon system \eqref{EKGintro}. However, it does not provide insight into the  asymptotic behaviour of the Fourier coefficients that is crucial to understand the asymptotic behavior of solutions of the resonant system. Moreover, their method does not allow a priori to understand in which sense does the perturbative series expansion \eqref{PerturbativeSeriesInitialData} converge. For this reason, several research efforts have been focused on either explicitly computing the Fourier coefficients or deriving ultraviolet asymptotics \cite{2020arXiv200411049C, MR3435559,MR3435560,MR3512262,MR3270350,PhysRevLett.113.071601,PhysRevD.92.084001,zbMATH07227949}. Having a robust method that indicates how one can treat the Fourier coefficients would not only be beneficial for numerical studies but also for rigorous analysis.  \\ \\
We have previously proposed a novel systematic way to derive closed formulas the Fourier coefficients associated to quadratic and cubic nonlinearities for simpler wave equations. We refer the reader to \cite{SakisJacques} where we used KAM tools and constructed time-periodic solutions to semi-linear toy models on fixed AdS as well as to \cite{chatzikaleas2023timeperiodicsolutionsnekhoroshev} where we rigorously proved their non-linear stability over exponentially long times. In these works, both the existence and the non-linear stability of the time-periodic solutions depend on non-degeneracy conditions similar to the ones described above in the numerical work of Rostworowski-Maliborski \cite{PhysRevLett.111.051102}. There, in order to rigorously verify the validity of these conditions, we used a powerful algorithm for hypergeometric summation due to Zeilberger \cite{MR1048463} that allowed us to systematically derive 1-step linear recurrence relations for several complicated hypergeometric sequences that occur in the process. \\ \\
The aim of this paper is to build upon the novel systematic way proposed in \cite{chatzikaleas2023timeperiodicsolutionsnekhoroshev} and derive computationally efficient formulas for the Fourier coefficients to the Einstein equations \eqref{EKGintro} in spherical symmetry (Theorem \ref{TheoremMain2}) and then use them to establish their asymptotic behavior as the mode numbers tend to infinity, including the explicit determination of the leading-order coefficients. To do so, besides Zeilberger's algorithm \cite{MR1048463}, we also use another powerful algorithm for hypergeometric summation due to Gosper \cite{MR485674} that allows us to systematically find discrete antiderivatives in closed formulas for several complicated hypergeometric sequences.

\subsection{The Einstein-massless-scalar field system in the Fourier space}\label{SubSectionEKGinFourierSpace}
With respect to the compactified coordinates
\begin{align*}
	(t,x,\omega) \in \mathbb{R} \times \left[ 0 , \frac{\pi }{2} \right) \times \mathbb{S}^2,
\end{align*}
the AdS metric is given by
\begin{align*}
	\text{g}_{\text{AdS}} (t,x,\omega) = \frac{1}{\cos^2(x)} \left( -dt^2 + dx^2 + \sin^2(x) d\omega^2 \right),
\end{align*}
where $d\omega^2 $ denotes the standard round metric on $\mathbb{S}^2$. Following \cite{PhysRevLett.111.051102,PhysRevLett.107.031102}, we parametrize the spacetime metric in \eqref{EKGintro} by
\begin{align*}
	\text{g}  (t,x,\omega) = \frac{1}{\cos^2(x)} \left( -A(t,x)e^{-2\delta(t,x)}dt^2 + A^{-1}(t,x)dx^2 + \sin^2(x) d\omega^2 \right),
\end{align*}
where, for convenience, we fixed the cosmological constant $\Lambda=-3$. With this ansatz, the wave equation for the self-gravitating massless scalar field $\phi (t,\cdot)$ in \eqref{EKGintro} is given by the system
\begin{align}
	\partial_{t} \Phi (t,x) & = \partial_{x} \big( A(t,x)e^{-\delta(t,x) } \Pi(t,x)  \big), \label{EquationEKG1} \\
	\partial_{t} \Pi (t,x) & = e^{-\delta (t,x)} 
	\left(A(t,x) \partial_{x} \Phi(t,x) +2\Phi (t,x)\frac{2-(1-A(t,x))\cos(2x)}{ \sin(2x)} \right), \label{EquationEKG2}
\end{align}
coupled to the Einstein equation in \eqref{EKGintro} given by
\begin{align}
	\partial_{x} A(t,x) &= \frac{1+2\sin^2(x)}{\sin(x) \cos(x)} (1-A(t,x)) - \sin(x) \cos(x) A(t,x) \big(\Phi^2 (t,x)+\Pi^2(t,x) \big),\label{EquationEKG3} \\
	\partial_{x} \delta (t,x)&=  -\sin(x) \cos(x)\big(\Phi^2 (t,x)+\Pi^2(t,x) \big).\label{EquationEKG4}
\end{align}
From now on, we omit the dependence on $(t,x)$ from the dynamical variables. The operator that governs the behaviour of the linearized solutions in \eqref{EquationEKG1}-\eqref{EquationEKG2} is given by
\begin{align}\label{LinearizedOperatorIntro}
	-L\phi=  \frac{1}{\tan^2(x)} \partial_{x}\left( \tan^2(x) \partial_{x}\phi  \right) =\partial_{x}^2 \phi + \frac{4}{ \sin(2x)}\partial_{x} \phi=\partial_{x} \Phi + \frac{4}{ \sin(2x)} \Phi, 
\end{align} 
for all $\phi \in \mathcal{D}(L)\simeq H^2 $.

\begin{lemma}[The Einstein-massless-scalar field system in the Fourier space]\label{LemmaEinsteinMassLessFourier}
 We define the auxiliary dynamical variables
\begin{align*}
	\xi :=e^{-\delta } (1-A  )   , \quad
	\zeta  :=e^{-\delta }   -1
\end{align*}
as well as the Fourier coefficients	
\begin{align*}
	\phi^{(m)}(t) &:=(\phi(t,\cdot)|e_m), \quad 
	\sigma^{(m)}(t):=(\Pi(t,\cdot)|e_m), \\ 
	\xi^{(m)}(t) &:=(\xi(t,\cdot)|e_m),  \quad 
	\zeta^{(m)}(t):=(\zeta(t,\cdot)|e_m),
\end{align*}
for all integers $m \geq 0$. Then, the Einstein-massless-scalar field system \eqref{EquationEKG1}-\eqref{EquationEKG4} can be written in the Fourier space as
\begin{align*} 
	\partial_{t} \begin{pmatrix}
		\phi^{(m)} (t)   \\
		\sigma^{(m)} (t)   
	\end{pmatrix}+
	\begin{pmatrix}
		0&-1 \\
		\omega_{m}^2 & 0 
	\end{pmatrix} 
	\begin{pmatrix}
		\phi^{(m)} (t)   \\
		\sigma^{(m)} (t)   
	\end{pmatrix}=
	\begin{pmatrix}
		F_1^{(m)} (t)   \\
		F_2^{(m)} (t)   
	\end{pmatrix},
\end{align*}
for all integers $m \geq 0$, where 
\begin{align*}
	F_1^{(m)} (t)   &= \sum_{i,j=0}^{\infty}\mathbf{C}_{ijm} \big( \zeta^{(i)} (t)  -\xi^{(i)} (t)   \big)\sigma^{(j)}   (t)     ,   \\ 
	F_2^{(m)}  (t)  &= \sum_{i,j=0}^{\infty}\left[  
	 \left( 4\mathbf{D}_{ijm }
	+\omega_{j}^2  \mathbf{C}_{ijm}-2  \mathbf{E}_{ijm } \right) \xi^{(i)} (t)  
	-  
	\omega_{j}^2  \mathbf{C}_{ijm}  \zeta^{(i)}    (t) 
	    \right]\phi^{(j)} (t)  . 
\end{align*}
and the Fourier coefficients are defined by
\begin{align*} 
	\mathbf{C}_{ijk} := \left(e_{i}e_{j} \big|e_{k}\right) ,  \quad  
\mathbf{D}_{ijk } := \left(e_{i}e_{j}^{\prime} \Big|\frac{e_k}{\sin (2 \cdot)}\right)  , \quad  	
\mathbf{E}_{ijk } := \left(e_{i}e_{j}^{\prime} \Big|\frac{e_{k}}{\tan (2 \cdot)}\right).  
\end{align*}
\end{lemma}
\begin{proof}
	Using the auxiliary dynamical variables $	\xi$ and $\zeta$, \eqref{EquationEKG1}-\eqref{EquationEKG2} can be written as
\begin{align}
	\partial_{t} \Phi  & = \Pi ^{\prime}  +\big( \zeta  \Pi  -\xi  \Pi   \big)^{\prime}, \label{EquationEKG1New}  \\
	\partial_{t} \Pi   
	& =-L\phi +  \zeta  \Phi^{\prime}  -\xi   \Phi ^{\prime}
	+ \frac{4}{ \sin(2x)}\zeta\Phi
	- \frac{2}{ \tan(2x)}\xi\Phi ,\label{EquationEKG2New} 
\end{align}
%and
%\begin{align}
%	\left(\frac{\sin(x)}{\cos^2(x)} \xi  \right)^{\prime} & =  (\zeta+1) (\Phi^2+\Pi^2) \tan^2(x)  , \label{EquationEKG3New}  \\
%	\zeta^{\prime} &= (\zeta+1)(\Phi^2+\Pi^2)\sin(x)\cos(x) , \label{EquationEKG4New} 
%\end{align}
% respectively. 
Here, we use a prime to denote the derivative operator $ \partial_{x}$. Let $ \{e_{ n }:  n  \geq 0\} $ be the set of all eigenfunctions to the linearized operator $L$ given by in \eqref{LinearizedOperatorIntro}. In Section \ref{Sectionlineareigenvalueproblem}, we show that $ \{e_{ n }:  n  \geq 0\} $  forms a orthogonal and complete basis for the weighted function space $  L^2 ( (0,\pi/2);\tan^2(x)dx ) $. Hence, for $\phi (t,\cdot) \in L^2 ( (0,\pi/2);\tan^2(x)dx ) $, there exist coefficients $\{\phi^{(n)}  (t): n \geq 0 \} \in l^2(\mathbb{N})$ so that
\begin{align*}
	\phi   =\sum_{n=0}^{\infty} \phi^{(n)} e_{n}. 
\end{align*}
This leads us to further expand 
\begin{align*}
	\Pi    = \sum_{n=0}^{\infty} \sigma^{(n)}    e_{n} , \quad 
	\Phi   =\sum_{n=0}^{\infty} \phi^{(n)}   e_{n}^{\prime} ,  \quad 
	\xi =\sum_{n=0}^{\infty} \xi^{(n)}   e_{n} , \quad
	\zeta  =\sum_{n=0}^{\infty} \zeta^{(n)}   e_{n}, 
\end{align*}
for some coefficients $\{\sigma^{(n)}  (t): n \geq 0 \}$, $\{\xi^{(n)}  (t): n \geq 0 \}$ and $\{ \zeta^{(n)}  (t): n \geq 0 \}$. In Section \ref{Sectionlineareigenvalueproblem}, we also show that $ \{\omega_{ n }^2:  n  \geq 0\} $, with $\omega_{n}=2n+3$, for all integers $n \geq 0$, is the set of all eigenvalues to the linearized operator $L$. Moreover, we introduce an inner product $(\cdot |\cdot)$ for functions in $L^2 ( (0,\pi/2);\tan^2(x)dx )$ so that the orthogonality conditions $\left( e_{n} | e_{m} \right) =\mathds{1}(n=m) $ and $\left( e_{n}^{\prime} | e_{m}^{\prime} \right) =\omega_{m}^2\mathds{1}(n=m) $ hold true, for all integers $n,m \geq 0$. In the following, we fix an integer $m \geq 0$ and use these orthogonality conditions to obtain  \eqref{EquationEKG1New}-\eqref{EquationEKG2New} in the Fourier space. Decomposing equations \eqref{EquationEKG1New} and \eqref{EquationEKG2New} into modes, we infer
\begin{align*}  
	%  \sum_{n=0}^{\infty} \partial_{t}  \phi^{(n)}   e_{n}^{\prime}  & = \sum_{n=0}^{\infty} \sigma^{(n)}    e_{n}^{\prime}+  \sum_{i,j=0}^{\infty} \zeta^{(i)}   \sigma^{(j)}     	( e_{i} e_{j} )^{\prime}	-\sum_{i,j=0}^{\infty} \xi^{(i)}   \sigma^{(j)}    	( e_{i} e_{j} )^{\prime}	 \Longrightarrow \\
	    \omega_{m}^2    \partial_{t} \phi^{(m)}     & =   \omega_{m}^2   \sigma^{(m)}    +
	\sum_{i,j=0}^{\infty} \zeta^{(i)}   \sigma^{(j)}   
	\left(  
	( e_{i} e_{j} )^{\prime} | e_{m}^{\prime} \right)
	-\sum_{i,j=0}^{\infty} \xi^{(i)}   \sigma^{(j)}   \left(  
	( e_{i} e_{j} )^{\prime} | e_{m}^{\prime} \right). 
\end{align*}
and
\begin{align*} 
	% \sum_{n=0}^{\infty} \partial_{t}  \sigma^{(n)}    e_{n}     & =  - \sum_{n=0}^{\infty} \phi^{(n)}   \omega_{n}^2 e_{n}   +	 \sum_{i,j=0}^{\infty} \zeta^{(i)}    \phi^{(j)}   e_{i}  e_{j}^{\prime\prime}    -	 \sum_{i,j=0}^{\infty} \xi^{(i)}     \phi^{(j)}   e_{i}  e_{j}^{\prime\prime}   \\	&+\sum_{i,j=0}^{\infty}  \phi^{(i)}    \zeta^{(j)}    \frac{4}{ \sin(2x)}e_{i}  	 e_{j}^{\prime}    	-2\sum_{i,j=0}^{\infty} \xi^{(i)}    \phi^{(j)}     \frac{1}{ \tan(2x)}e_{i}  e_{j}^{\prime}   	  \Longrightarrow \\
	      \partial_{t} \sigma^{(m)}       & = -\omega_{m}^2 \phi^{(m)}   +
	 \sum_{i,j=0}^{\infty} \zeta^{(i)}    \phi^{(j)}    
	 \left(e_{i}e_{j}^{\prime\prime} | e_{m} \right) -
	 \sum_{i,j=0}^{\infty} \xi^{(i)}     \phi^{(j)}   
	 \left( e_{i} e_{j}^{\prime\prime} | e_{m} \right)  \\
	&+\sum_{i,j=0}^{\infty}   \zeta^{(i)}   \phi^{(j)}   
	\left( \frac{4}{ \sin(2 \cdot)}e_{i} 	 e_{j}^{\prime} \Bigg| e_{m} \right)  
	-2\sum_{i,j=0}^{\infty} \xi^{(i)}    \phi^{(j)}    
	\left(\frac{1}{ \tan(2 \cdot )}e_{i} e_{j}^{\prime} \Bigg| e_{m} \right) .
\end{align*}
respectively. Since $ \{\omega_{ n }^2:  n  \geq 0\} $ is the set of all eigenvalues to the linearized operator $L$, we have that
\begin{align*}
	L e_n  = \omega^2 _{n} e_n  \Longleftrightarrow -  \left( e_{n}^{\prime}  \tan^2 (x)\right)^{\prime} = \omega_{n}^2  e_n  \tan^2 (x),
\end{align*}
for all integers $n \geq 0$, and integration by parts yields
\begin{align*}
	\left(  ( e_{i} e_{j} )^{\prime} | e_{m}^{\prime} \right) &=
	\int_{0}^{\pi /2}  ( e_{i}(x) e_{j}(x) )^{\prime}   e_{m}^{\prime}(x) \tan^2 (x) dx 
	=
	- \int_{0}^{\pi /2}   e_{i}(x) e_{j}(x)  \left(   e_{m}^{\prime}(x) \tan^2 (x) \right)^{\prime} dx \\
	&=\omega_{m}^2 \int_{0}^{\pi /2}   e_{i}(x) e_{j}(x)   e_m(x) \tan^2 (x) dx 
	=\omega_{m}^2 \left( e_{i} e_{j}  | e_m  \right) , \\
	\left(e_{i}e_{j}^{\prime\prime} | e_{m} \right) &= 
	\int_{0}^{\pi /2}  e_{i}(x)e_{j}^{\prime\prime}(x)  e_{m}(x)\tan^2 (x) dx 
	= -
	\int_{0}^{\pi /2}  e_{i}(x)  e_{m}(x)\left( \frac{4 }{\sin (2 x)}e_{j}^{\prime}(x) + \omega_{j}^2 e_{j}(x)  \right)\tan^2(x)  dx \\
	&=-
	\int_{0}^{\pi /2} \frac{4 }{\sin (2 x)} e_{i}(x)   e_{j}^{\prime}(x)e_{m}(x) \tan^2(x)  dx
	-\omega_{j}^2
	\int_{0}^{\pi /2}  e_{i}(x)   e_{j}(x) e_{m}(x)   \tan^2(x)  dx \\
	&=-\left( \frac{4}{ \sin(2 \cdot)}e_{i} 	 e_{j}^{\prime} \Bigg| e_{m} \right)
	-\omega_{j}^2   \left( e_{i} e_{j}  | e_m  \right). 
\end{align*} 
Putting all together, yields
\begin{align*}
%%%%%%%%%%
	       \partial_{t}  \phi^{(m)}     & =       \sigma^{(m)}    +
	\sum_{i,j=0}^{\infty} \left( \zeta^{(i)}    
	-\xi^{(i)}   \right)\sigma^{(j)}   
	  \left( e_{i} e_{j}  | e_m  \right), \\
%%%%%%%%%%	
	     \partial_{t} \sigma^{(m)}       & = - \omega_{m}^2 \phi^{(m)}   +
	 \sum_{i,j=0}^{\infty} 
	 \left(    \xi^{(i)}   -\zeta^{(i)}   \right)\phi^{(j)}    
	 \left( 4\left( \frac{1}{ \sin(2 \cdot)}e_{i} 	 e_{j}^{\prime} \Bigg| e_{m} \right)
	+\omega_{j}^2   \left( e_{i} e_{j}  | e_m  \right) \right)   \\
	&+4\sum_{i,j=0}^{\infty}   \zeta^{(i)}   \phi^{(j)}   
	\left( \frac{1}{ \sin(2 \cdot)}e_{i} 	 e_{j}^{\prime} \Bigg| e_{m} \right)  
	-2\sum_{i,j=0}^{\infty} \xi^{(i)}    \phi^{(j)}    
	\left(\frac{1}{ \tan(2 \cdot )}e_{i} e_{j}^{\prime} \Bigg| e_{m} \right). 
\end{align*}
Finally, we define the Fourier coefficients $\mathbf{C}_{ijk}$, $\mathbf{D}_{ijk }$ and $\mathbf{E}_{ijk } $ as above and arrive at the first-order infinite dimensional differential system as stated above. 
\end{proof}
\noindent
According to \cite{PhysRevLett.111.051102}, the constrains \eqref{EquationEKG3}-\eqref{EquationEKG4} can also be written in the Fourier space, where in this case the Fourier coefficients are given by
\begin{align}\label{DefinitionAandB} 
\mathbf{A}_{ij} & := \left(e_{i}e_{j} \big|\cos(2\cdot)\right) , \quad \mathbf{B}_{ij} :=\left(e_{i}^{\prime}e_{j} \big|\sin(2\cdot)\right).
\end{align}
We refer the reader to \cite{PhysRevLett.111.051102,PhysRevLett.107.031102} for more details.

%
%\subsection{Organization of the paper}
%\noindent
%In Section \ref{SubSectionEKGinFourierSpace}, we derive the Einstein-massless-scalar field system in the Fourier space. In Section \ref{SubSectionEKGinFourierSpace}, we derive the Einstein-Klein-Gordon system in the Fourier space.
% 
% 

\section{Evaluating sums of hypergeometric functions}  \label{SectionAlgorithms}
\noindent
By definition, each Fourier coefficient is an oscillating integral of products of the eigenfunctions $  e_{ n } $ to the associated linearized operator. In Section \ref{Sectionlineareigenvalueproblem}, we show that the sequence $ \{e_{ n }:  n  \geq 0\} $ forms an orthogonal and complete basis for a weighted function space $  L^2 ( (0,\pi/2);\tan^2(x)dx ) $. Moreover, we also show that each eigenmode is a weighted Jacobi polynomial and hence a weighted hypergeometric function. In general, the parameters of the Jacobi polynomials depend on the choice of the mass and, in the case of zero mass, these are given by\footnote{Here, the symbol $\sim $ means ``up to a normalization constant''.} 
 \begin{align*}
	e_{n}(x) \sim \cos^3 (x)~ P_{n}^{(1/2,3/2)}(\cos(2x)) \sim \cos^3 (x)~ {}_2 F_1 \left( -n,3+n;\frac{3}{2},\sin^2(x) \right) ,
\end{align*}
for all integers $n \geq 0$ and real $x\in [0,\pi/2]$. In other words, one can interpret the Fourier coefficients as oscillating integrals of products of hypergeometric functions. On the one hand, once these are multiplied inside the integrals (possibly with extra weights), one needs to express their resulting product in terms of the original base $ \{e_{ n }:  n  \geq 0\} $. These are called linearization formulas and are of the form 
\begin{align}\label{LinearizationFormulaIntro}
P_{n}^{(1/2,3/2)}(y)P_{m}^{(1/2,3/2)}(y)
= \sum_{k=|n-m|}^{n+m} L_{k}(n,m) P_{k}^{(1/2,3/2)}(y),
\end{align}
for all integers $n,m\geq 0$ and real $y\in[-1,1]$. Although Rahman \cite{MR634149, MR3527696} expressed the linearization coefficients $L_k(n,m)$ in terms of a generalized hypergeometric function ${}_9F_8(1)$, this representation does not, in general, yield explicit closed-form expressions for these coefficients. Nonetheless, one can use identities of Jacobi polynomials, such as  \eqref{AppendixRelationJacobitoGegenbauer}, to rewrite these special functions in terms of ``reduced'' orthogonal polynomials, called Gegenbauer polynomials, which is a subset of the Jacobi family with  better structure. In this case, the linearization coefficients are given explicitly via formulas like \eqref{AppendixLinearizatonFormula}. On the other hand, due to the presence of derivatives and weights, one also needs to express a given Jacobi polynomial with a fixed pair of parameters, say $(\alpha,\beta)$, in terms of the base $ \{e_{ n }:  n  \geq 0\} $ consisting of Jacobi polynomials with the original pair of parameters, namely $(1/2,3/2)$. These are called connection formulas and are of the form
\begin{align}\label{ConnectionFormulaIntro}
	P_{n}^{(\alpha,\beta)}(y)  = \sum_{k=0}^{n} c_{k}(n,\alpha,\beta)  P_{k}^{(1/2,3/2)}(y),
\end{align}
for all integers $n \geq 0$ and real $y\in [-1,1]$. In the case of Gegenbauer polynomials, the connection coefficients are given explicitly, see for example \eqref{AppendixConnectionFormulaExponent5half}-\eqref{AppendixConnectionFormulaExponent7half}.
\subsection{Methodology}
To begin with, we use properties of the Gegenbauer polynomials, such as \eqref{LinearizationFormulaIntro} and \eqref{ConnectionFormulaIntro}, together with orthogonality conditions of orthogonal polynomials \eqref{AppendixOrthogonalityGegenbauerExponent1}-\eqref{AppendixOrthogonalityGegenbauerExponent7half} and several recurrence relations \eqref{AppendixRecrrenceRelation1}-\eqref{AppendixRecrrenceRelation4}, to replace the integral in the definition of each of the Fourier coefficients with a sum obtaining the so called ``computationally efficient formulas''. We find that  the Fourier coefficients $ \mathbf{C}_{ijk}$, $ \mathbf{D}_{ijk}$ and $ \mathbf{E}_{ijk}$ defined in Lemma \ref{LemmaEinsteinMassLessFourier} are of the form
\begin{align}\label{DefinitionSumfGegenralStrategy}
	\mathbf{f}_{ijk}= \sum_{\lambda=0}^{\Lambda}   \mathbf{F}_{ijk}(\lambda),
\end{align} 
for some function $  \mathbf{F}_{ijk}(\lambda)$ given in closed formula depending on the Fourier coefficient we consider. Here, $\Lambda$ is an integer that may depend on $i,j,k$. Crucially,  all functions $ \mathbf{F}_{ijk}(\lambda)$ that appear above turn out to be hypergeometric functions with respect to all their arguments, meaning that all 
\begin{align*}
	\frac{ \mathbf{F}_{ijk}(\lambda+1)}{ \mathbf{F}_{ijk}(\lambda)}, \quad 
	\frac{ \mathbf{F}_{i+1,j,k}(\lambda)}{ \mathbf{F}_{ijk}(\lambda)}, \quad 
	\frac{ \mathbf{F}_{i,j+1,k}(\lambda)}{ \mathbf{F}_{ijk}(\lambda)}, \quad 
	\frac{ \mathbf{F}_{i,j,k+1}(\lambda)}{ \mathbf{F}_{ijk}(\lambda)}
\end{align*}
are rational functions of $i,j,k$ and $\lambda$. Then, we implement  two powerful algorithms that use methods of creative telescoping and yield 
\begin{itemize}[leftmargin=*]
	\item small-step linear recurrence relations (Zeilberger's  \cite{MR1048463} algorithm),
	\item  discrete anti-derivatives (Gosper's \cite{MR485674} algorithm), 
\end{itemize}
for sums like \eqref{DefinitionSumfGegenralStrategy}. Thanks to \cite{MR1048463, MR485674}, as well as to improvements by Paule, Riese, Wegschaider and Schorn \cite{MR1953199, MR2001619, MR2768529, MR1880090, MR1448687, MR1395420,MR1962799, MR1394550}, proving most definite hypergeometric summation and transformation formulas  becomes routine work that can also be performed by Mathematica using  Fast Zeilberger Package \cite{MathematicaFastZeilbergerPackage} and MultiSum Package \cite{MathematicaMultiSumPackage} both in Mathematica. Most importantly, we note that all results obtained by these algorithms can be rigorously proved. 

\subsection{Zeilberger's algorithm}
Assume that we are interested in finding a linear recurrence relations for the sum $\mathbf{f}_{ijk}$ in \eqref{DefinitionSumfGegenralStrategy} with respect to $i$. For simplicity, we ignore the fixed indices $j,k$ and, with a slight abuse of notation, we write $\mathbf{f}_{i}$ for the sum and  $ \mathbf{F}_{i}(\lambda)$ for the summand. We claim that it suffices to find a recurrence relation for the summand $ \mathbf{F}_{i}(\lambda)$ with respect to $i$. In other words, given a hypergeometric function $ \mathbf{F}_{i}(\lambda)$ and an integer $\Lambda$, we would like to find a function $ \mathbf{G}_i(\lambda)$, coefficients $\{ \mathbf{c}_{i}(\alpha): \alpha =0,1,\dots,N\}$, as well as a (possibly large) integer $N$, all in closed forms, so that  
\begin{align}\label{RecurrenceRelationFGegenralStrategy}
	\sum_{\alpha=0}^{N} \mathbf{c}_{i}(\alpha) \mathbf{F}_{i+\alpha }(\lambda)=\mathbf{G}_{i}(\lambda+1)-\mathbf{G}_{i}(\lambda).
\end{align}
Notice that the coefficients $\{ \mathbf{c}_{i}(\alpha): \alpha =0,1,\dots,N\}$ are independent of  $\lambda$ and the right-hand-side of \eqref{RecurrenceRelationFGegenralStrategy} is in a  ``summable'' form. A linear recurrence relation for the sum $\mathbf{f}_i$ then follows immediately by summing both sides of \eqref{RecurrenceRelationFGegenralStrategy} with respect to $\lambda$, 
\begin{align}\label{RecurrenceRelationfGegenralStrategyZeilberger}
	\sum_{\alpha=0}^{N}  \mathbf{c}_{i}(\alpha)\mathbf{f}_{i+\alpha } &=
	\sum_{\alpha=0}^{N}   \mathbf{c}_{i}(\alpha) \sum_{\lambda=0}^{\Lambda} \mathbf{F}_{i+\alpha }(\lambda) 
	%=\sum_{\lambda=0}^{\Lambda}  \left(\mathbf{G}_{i}(\lambda+1)-\mathbf{G}_{i}(\lambda) \right)
	= \mathbf{G}_{i}(\Lambda+1)-\mathbf{G}_{i}(0) .
\end{align}   
  The coefficients $\{\mathbf{c}_{i}(\alpha): \alpha =0,1,\dots,N\}$ and the function $\mathbf{G}_{i}(\lambda)$ can be found explicitly using Zeilberger's algorithm, see \cite{MR1048463} and Chapter 6 in \cite{MR1379802}.  The significance of the recurrence relation \eqref{RecurrenceRelationFGegenralStrategy} for the summand is twofold. First, according to Theorem 6.2.1, page 105, in \cite{MR1379802}, for all proper hypergeometric functions $\mathbf{F}_{i}(\lambda)$, it always exists for a (possibly large) integer $N$. Second, once $\{\mathbf{c}_{i}(\alpha): \alpha =0,1,\dots,N\}$ and $\mathbf{G}_{i}(\lambda)$ are found in closed formulas, then one can rigorously prove the validity of \eqref{RecurrenceRelationFGegenralStrategy}  by a (possibly long but) straightforward computation. Finally, we also note that, in the case where the resulting recurrence relation \eqref{RecurrenceRelationfGegenralStrategyZeilberger} cannot be solved in a reasonable sense, one can use further algorithms. For example, Petkovsek's algorithm, see Chapter 8 in \cite{MR1379802}, will find the solution to a (possibly complicated) recurrence \eqref{RecurrenceRelationfGegenralStrategyZeilberger} provided that it has a solution $\mathbf{f}_i$ that is a linear combination of a fixed number of hypergeometric terms in $i$. We refer the reader to \cite{MR1644447,MR1379802,MR1048463} for more details.

\subsection{Gosper's algorithm}
Assume that we are interested in finding a closed formula for the  sum $\mathbf{f}_{ijk}$ in \eqref{DefinitionSumfGegenralStrategy}. Then, given a hypergeometric function $ \mathbf{F}_{ijk}(\lambda)$, we would like to find a function $ \mathbf{G}_{ijk}(\lambda)$ in closed form that puts $ \mathbf{F}_{ijk}(\lambda)$ in a  ``summable'' form meaning that  
\begin{align}\label{RecurrenceRelationFGegenralStrategy2}
	  \mathbf{F}_{ijk }(\lambda)=\mathbf{G}_{ijk}(\lambda+1)-\mathbf{G}_{ijk}(\lambda).
\end{align}
A closed formula for the sum $\mathbf{f}_{ijk}$ then follows immediately by summing both sides of \eqref{RecurrenceRelationFGegenralStrategy2} with respect to $\lambda$, since
\begin{align} \label{RecurrenceRelationFGegenralStrategyExample2}
\mathbf{f}_{ijk}=\sum_{\lambda=0}^{\Lambda}     \mathbf{F}_{ijk }(\lambda)  &
% =\sum_{\lambda=0}^{\Lambda}  \left( \mathbf{G}_{ijk}(\lambda+1)-\mathbf{G}_{ijk}(\lambda) \right) 
=    \mathbf{G}_{ijk}(\Lambda+1)-\mathbf{G}_{ijk}(0)  .
\end{align}  
The function $\mathbf{G}_{ijk}(\lambda)$ can be found explicitly using Gosper's algorithm, see \cite{MR485674} and Chapter 5 in \cite{MR1379802}. The significance of the formulation \eqref{RecurrenceRelationFGegenralStrategy2} for the summand is again twofold. First, according to Theorem 5.2.1, page 76, in \cite{MR1379802}, finding the hypergeometric functions $\mathbf{G}_{ijk}(\lambda)$ in \eqref{RecurrenceRelationFGegenralStrategy2} is equivalent to finding polynomial solutions of a first order recurrence relation. In addition, according to Theorem 5.6.3, page 94 in \cite{MR1379802}, if Gosper's algorithm does not succeed in finding the discrete anti-derivative $\mathbf{G}_{ijk}(\lambda)$, then this means that the sum was not expressible in closed form in the first place. Second, once $\mathbf{G}_{ijk}(\lambda)$ is found in closed formula, then one can rigorously prove the validity of \eqref{RecurrenceRelationFGegenralStrategy2}  by a (possibly long but) straightforward computation. We refer the reader to \cite{MR1644447,MR1379802,MR1048463} for more details.

\section{Main result}
\noindent
\subsection{Statement of the main result}
In this note, we develop a novel abstract method that derives computationally efficient formulas with the Fourier coefficients. Our main result (Theorem \ref{TheoremMain2}) uses the strategy discussed in Section \ref{SectionAlgorithms} and rigorously establishes computationally efficient formulas for the Fourier coefficients. 

\begin{theorem}[Main result: Computationally efficient formulas for the Fourier coefficients]\label{TheoremMain2}
	We derive the computationally efficient formulas for the Fourier coefficients $\{\mathbf{A}_{ij}:i,j \geq 0\}$, $\{\mathbf{B}_{ij}:i,j\geq 0\}$, $\{\mathbf{C}_{ijk}:i,j,k \geq 0\}$, $\{\mathbf{D}_{ijk}:i,j,k \geq 0\}$ and $\{\mathbf{E}_{ijk}:i,j,k \geq 0\}$ defined in \eqref{DefinitionAandB} and Lemma \ref{LemmaEinsteinMassLessFourier} respectively. 
\end{theorem}
\begin{proof}
	The complete statement is given in Propositions \ref{ComputationallyEfficientFormulaAij}-\ref{ComputationallyEfficientFormulaEijk}. 
\end{proof}
\noindent
As an application, we rigorously establish the limiting behaviour of the Fourier coefficients as the mode numbers go to infinity. We consider interactions involving 1-mode, 2-modes and 3-modes and, for simplicity, whenever an index is fixed, we assume that it is the first zero-mode. Moreover, since the Fourier coefficients  $\mathbf{A}_{ij}$ and $\mathbf{B}_{ij}$ concentrate only on finitely many modes (see Propositions \ref{ComputationallyEfficientFormulaAij} and \ref{ComputationallyEfficientFormulaBij}), we only focus on $\mathbf{C}_{ijk}$, $\mathbf{D}_{ijk}$ and $\mathbf{E}_{ijk}$.  
 
\begin{theorem}[Application: Limiting behaviour for the Fourier coefficients]\label{MainTheorem}
	The following limits for the Fourier coefficients listed in the Table below hold true. 
\begin{center} 
{\tabulinesep=2mm 
\begin{longtabu}{ |p{3cm}|p{4cm}|p{6cm}|p{3.5cm}|  }
\hline
\multicolumn{4}{|c|}{Limiting behaviour for the Fourier coefficients} \\
\hline
Varying modes & \hfil Fourier coefficient & \hfil Leading order terms as $\omega_{  \nu  } \rightarrow \infty$ & \hfil   Behaviour  \\
\hline
  & \hfil  $\mathbf{C}_{ijk}$   & \hfil  {$ \!\begin{aligned} 
               \mathbf{C}_{00  \nu  } &=\frac{9216}{\pi ^{3/2}} \frac{(-1)^  \nu  }{ \omega_{  \nu  }^6} +\mathcal{O} \left (\frac{1}{\omega_{\nu}^{7}}\right )   \end{aligned}$} & \hfil Decay  \\ \cline{2-4}
\center{1-mode}  & \hfil $\mathbf{D}_{ijk}$   & \hfil {$   \begin{aligned} 
               \mathbf{D}_{00  \nu  } &=\frac{768}{\pi ^{3/2}}\frac{ (-1)^  \nu  }{ \omega_{  \nu  }^4}+\mathcal{O}\left (\frac{1}{\omega_{\nu}^{5}}\right )  \\
                \mathbf{D}_{0  \nu  0} &=-\frac{16}{\sqrt{\pi } }  +\mathcal{O}\left (\frac{1}{\omega_{\nu}}\right ) 
                 \end{aligned} $}   & \hfil Decay/Constant  \\ \cline{2-4} 
& \hfil $\mathbf{E}_{ijk}$   & \hfil  {$   \begin{aligned} 
               \mathbf{E}_{00  \nu  } &= -\frac{768}{\pi ^{3/2}}\frac{(-1)^  \nu  }{\omega_{  \nu  }^4}  +\mathcal{O}\left (\frac{1}{\omega_{\nu}^{5}}\right )  \\
                \mathbf{E}_{0  \nu  0} &=-\frac{16}{\sqrt{\pi }}+\mathcal{O}\left (\frac{1}{\omega_{\nu}}\right )
                 \end{aligned} $}  & \hfil Decay/Constant  \\  
\hline 
  & \hfil $\mathbf{C}_{ijk}$   & \hfil   {$ \!\begin{aligned} 
               \mathbf{C}_{0  \nu    \nu  } &=\frac{16 }{3\pi }\sqrt{\frac{2}{\pi }}  +\mathcal{O}\left (\frac{1}{\omega_{\nu}}\right )\end{aligned}$} & \hfil Constant  \\ \cline{2-4}
\center{2-modes}  & \hfil $\mathbf{D}_{ij k }$   & \hfil {$   \begin{aligned} 
               \mathbf{D}_{0  \nu    \nu  } &=-2\sqrt{\frac{2}{\pi }} \omega_  \nu   +\mathcal{O}\left (1\right ) \\   
                \mathbf{D}_{  \nu  0  \nu  } &=-\frac{12}{\pi }\sqrt{\frac{2}{\pi }}   +\mathcal{O}\left (\frac{1}{\omega_{\nu}}\right )  
                 \end{aligned} $}  & \hfil Growth/Constant  \\  \cline{2-4}
& \hfil $\mathbf{E}_{ijk}$   & \hfil  {$   \begin{aligned} 
               \mathbf{E}_{0  \nu    \nu  } &= -2\sqrt{\frac{2}{\pi }}\omega_{  \nu  }    +\mathcal{O}\left ( 1\right )    \\
                \mathbf{E}_{  \nu  0  \nu  } &=-\frac{4  }{\pi } \sqrt{\frac{2}{\pi }}    +\mathcal{O}\left (\frac{1}{\omega_{\nu}}\right )  
                 \end{aligned} $}   & \hfil Growth/Constant  \\
\hline 
  & \hfil $\mathbf{C}_{ijk}$   & \hfil    {$ \!\begin{aligned} 
               \mathbf{C}_{  \nu    \nu    \nu  } &=\frac{2}{\sqrt{\pi }}   +\mathcal{O}\left (\frac{1}{\omega_{\nu}}\right )    \end{aligned}$} 
               & \hfil Constant  \\ \cline{2-4}
\center{3-modes}  & \hfil $\mathbf{D}_{ijk}$   & \hfil  {$   \begin{aligned} 
               \mathbf{D}_{  \nu    \nu    \nu  } &=-\frac{1}{2 \sqrt{\pi }} \omega_  \nu  ^2    +\mathcal{O}\left (\omega_{\nu}\right )    
                 \end{aligned} $}  & \hfil Growth  \\  \cline{2-4}
& \hfil $\mathbf{E}_{ijk}$   & \hfil  {$   \begin{aligned} 
               \mathbf{E}_{  \nu    \nu    \nu  } &= -\frac{1}{2 \sqrt{\pi }} \omega_{  \nu  }^2   +\mathcal{O}\left ( \omega_{\nu}\right )   
                 \end{aligned} $}  & \hfil Growth  \\   
\hline 
\end{longtabu} 
} 
\end{center} 
\end{theorem}
\begin{proof}
	The proof is given in Section \ref{ProofSectionTable}.
\end{proof}

\begin{remark}[Loss of derivatives]
 The growth of the Fourier coefficients is natural due to the quasi-nonlinearity of the system and it yields the following estimate with loss of derivatives, 
\begin{align*} 
		\| F_1 (t,\cdot) \|_{ L^2} +\| F_2 (t,\cdot) \|_{ L^2} &\lesssim 
		\left( \| \xi (t,\cdot) \|_{ L^2  }+\| \zeta (t,\cdot) \|_{ L^2  } \right)\left( \| \sigma (t,\cdot) \|_{ L^2  }+	\| \phi (t,\cdot) \|_{ H^2  } \right),
	\end{align*}
	where $F_1(t,\cdot) $ and $F_2(t,\cdot)$ are the non-linear terms in the physical space whose Fourier representations $ \{F_{1}^{(m)}(t): m \geq 0 \}$ and $ \{F_{2}^{(m)}(t): m \geq 0 \}$ are defined in Lemma \ref{LemmaEinsteinMassLessFourier}.  
\end{remark}

  \subsection{Acknowledgments}

 The authors would like to thank Professor Peter Paule at the Research Institute for Symbolic Computation in Linz  for giving us access information to the Zeilberger package in Mathematica. All computations presented here were done by hand and verified using Fast Zeilberger Package \cite{MathematicaFastZeilbergerPackage} and MultiSum Package \cite{MathematicaMultiSumPackage} in Mathematica \cite{Mathematica}.
 
 \subsection{Declaration of artificial intelligence use}
 No artificial intelligence (AI) was used in the production of these results or in any aspect of the work presented in this publication.

\section{The linear eigenvalue problem} \label{Sectionlineareigenvalueproblem}
\noindent
The operator that governs the linearized perturbations around AdS in $(1+3)$-dimensions with zero mass is given by
\begin{align*} 
	 L \phi = - \frac{1}{\tan^2(x)} \partial_{x} \left( \tan^2(x) \partial_{x} \phi  \right), 
\end{align*}
for a spherically symmetric scalar field $\phi: \mathbb{R} \times (0,\pi/2) \rightarrow \mathbb{R}$.  We consider the Hilbert space $L^2:=L^2 ( (0,\pi/2);\tan^2(x)dx )$ associated with the inner product
 \begin{align}\label{DefinitionInnerProduct}
 	(f|g)=  \int_{0}^{\pi /2} f(x) g(x) \tan^2(x)dx,
 \end{align}
for all $f,g\in L^2$, and define the domain of $L$ by
\begin{align*}
	\mathcal{D}(L)= \left\{\phi \in L^2: L\phi \in L^2\right\}.
\end{align*}
The linearized operator $L$ is generated by a closed sesquilinear form 
\begin{align*}
	\alpha :H^1 \times H^1 \  \rightarrow \mathbb{R}
\end{align*}
defined by
\begin{align*}
	\alpha (\phi,\psi) &= \int_{0}^{ \pi /2} \left( \partial_{x} \phi (x)\partial_{x}\psi(x)  + \phi (x) \psi(x)  \right) \tan^2(x)dx ,
\end{align*}
for all $\phi,\psi\in H^1$, where $H^1:=H^1 ( (0,\pi/2 );\tan^2(x)dx )$. Clearly, one has
	\begin{align*}
	  \alpha (\phi,\phi) = \|\phi \|_{H^1 ((0,\pi/2);\tan^2(x)dx)}^2,  
	\end{align*}
	for all $\phi  \in H^1$, hence the sesquilinear form is elliptic and strictly positive. In particular, $L$ is self-adjoint on its domain $\mathcal{D}(L)$. Therefore, all the eigenvalues of $L$ are real, any two eigenfunctions corresponding to different eigenvalues are orthogonal to each other and the set of eigenfunctions forms a complete (dense) basis for $L^2 $. The linear eigenvalue problem 
\begin{align*}
	L \phi =\omega^2 \phi,
\end{align*}
with $\phi \in \mathcal{D}(L)$, can be transformed into
\begin{align*}
	(1-y^2)   \psi^{\prime \prime}(y) +(1-4y)   \psi^{\prime}(y) +\frac{1}{4}\left(\omega^2 -  9\right  )  \psi(y) = 0, 
\end{align*}
with $\phi(x)=\cos^3(x)\psi(y)$ and $y=\cos(2x)$, that has a non-trivial solution if and only if $\psi(y)$ is $P_{n}^{(1/2,3/2)}(y)$, that is the Jacobi polynomials of degree $n$ with parameters $(1/2,3/2)$. Recall that the set $\{ P_{n}^{(\alpha,\beta)}(y)   : n \geq 0\}$ forms an orthonormal and complete basis for the weighted $L^2 ((-1,1);(1-y)^{\alpha} (1+y)^{\beta}dy)$. We refer the reader to \cite{MR0372517, MR2723248} for a review on Jacobi polynomials. Consequently, the eigenfunctions and eigenvalues to the linearized operators $L  $ are given respectively by
\begin{align*}
	e_{n}(x)=N_{n} \cos^3 (x) P_{n}^{ \left(1/2,3/2 \right)}(\cos(2x)), \quad \omega_{n}^2= (2n+3)^2,
\end{align*} 
for all integers $n \geq 0$. Here, $N_{n}  $ stands for a suitable normalization constant so that
\begin{align*}
	(e_n |e_{m})=\mathds{1}(n=m), \quad 
	\left( e_n^{\prime}   | e_{m}^{\prime} \right)=\omega_{n}^2\mathds{1}(n=m),
\end{align*}
for all integers $n,m \geq 0$, and it is given by
\begin{align*}
 	N_{n}=2\frac{\sqrt{ \Gamma (n+1) \Gamma (n+3)}}{ \Gamma  \left(n+\frac{3}{2} \right)},
 \end{align*}
for all integers $n \geq 0$, where $(\cdot|\cdot)$ denotes the inner product defined in \eqref{DefinitionInnerProduct}. This orthogonality condition follows immediately from the orthogonality of the Jacobi polynomials. We conclude this section, by expressing the eigenfunctions in terms of Gegenbauer polynomials $\{C_{n}^{(2)}(y): n \geq 0\}$. This is a special class of Jacobi polynomials that forms an orthonormal and complete basis for the weighted $L^2 ((-1,1); (1-y^2)^{3/2}dy )$.

\begin{lemma}[Eigenfunctions in terms of Gegenbauer polynomials]\label{LemmaEigenfunctionsWithRespectToGegenbauer}
	Let $n \geq 0$ be an integer and $x\in[0,\pi/2]$ a real number. Then, we have that
	\begin{align*}
		e_{n} (x) & =K_n \frac{(1-y^2)^{3/2}}{y}  C_{2 n+1}^{(2)}(y) ,  \quad
		e_{n}^{\prime} (x)  =K_n(1-y^2)^{1/2}  \frac{d}{dy} \left[ \frac{  (1-y^2 )^{3/2} }{y}C_{2 n+1}^{(2)}(y) \right] ,
	\end{align*}
	where $y=\sin(x)$ and the normalization constant reads
	\begin{align*}
		 K_n = \frac{(-1)^n}{\sqrt{\pi }} \frac{2}{\sqrt{(n+1) (n+2)}} .
	\end{align*}
\end{lemma}
\begin{proof}
	 We set $y=\sin(x)$ and using the symmetries of the Jacobi polynomials \eqref{AppendixSymmetriesJacobi} and the relation \eqref{AppendixRelationJacobitoGegenbauer} to express the Jacobi polynomials in terms of the Gegenbauer ones, we infer  
\begin{align*}
%%%%%%%%%%
	e_{n}(x) &=N_{n} \cos^3 (x) P_{n}^{ \left(1/2,3/2 \right)}(\cos(2x)) 
	=(-1)^n N_{n} \cos^3 (x) P_{n}^{ \left(3/2,1/2 \right)}(-\cos(2x)) \\
	&=(-1)^n N_{n} \cos^3 (x) P_{n}^{ \left(3/2,1/2 \right)}(-\cos(2x))
	= (-1)^n N_{n}   \frac{ \Gamma  \left(n+\frac{3}{2}\right) }{ \sqrt{\pi }  \Gamma  (n+3)   } \frac{\cos^3 (x)}{\sin(x)}  C_{2 n+1}^{(2)}(\sin (x)) \\
	& = K_n \frac{\cos^3 (x)}{\sin(x)}  C_{2 n+1}^{(2)}(\sin (x)) 
	= K_n \frac{(1-y^2)^{3/2}}{y}  C_{2 n+1}^{(2)}(y) , \\
%%%%%%%%%%	
e_{n}^{\prime}(x)  
&=K_n \frac{d}{dx} \left[ \frac{  (1-y^2 )^{3/2} }{y}C_{2 n+1}^{(2)}(y) \right]   =K_n \frac{d}{dy} \left[ \frac{  (1-y^2 )^{3/2} }{y}C_{2 n+1}^{(2)}(y) \right]  \cos(x),
\end{align*}
that completes the proof.
\end{proof}

\section{Computationally efficient formulas for the Fourier coefficients}  \label{SectionFourierCoefficients}
\noindent
In this section, we focus on establishing computationally efficient formulas for the Fourier coefficients (Propositions \ref{ComputationallyEfficientFormulaAij}-\ref{ComputationallyEfficientFormulaEijk}).

\subsection{Integrals of Gegenbauer polynomials}

To begin with, we use Lemma \ref{LemmaEigenfunctionsWithRespectToGegenbauer} to express the Fourier coefficients in terms of oscillating integrals of Gegenbauer polynomials.

\begin{lemma}[Fourier coefficients in terms of Gegenbauer polynomials]\label{LemmaFourierCoefficientsWithRespectToGegenbauer}
	Let $i,j,k \geq 0$ be three integers. Then, the Fourier coefficients are given by
	\begin{align*}
		\mathbf{A}_{ij} &=  K_{i}K_{j}\int_{0}^{1} (1-2y^2)  C_{2 i+1}^{(2)}(y)C_{2 j+1}^{(2)}(y)(1-y^2)^{3/2}   dy, \\
		\mathbf{B}_{ij} &= 2K_i K_j  \int_{0}^{1} y^2  \frac{d}{dy} \left[ \frac{  (1-y^2 )^{3/2} }{y}C_{2 i+1}^{(2)}(y) \right]   C_{2 j+1}^{(2)}(y) (1-y^2)  dy , \\
	\mathbf{C}_{ijk} &= K_{i}K_{j}K_{k}\int_{0}^{1} \frac{1}{y}
	C_{2 i+1}^{(2)}(y)
	C_{2 j+1}^{(2)}(y)
	C_{2 k+1}^{(2)}(y)  (1-y^2)^{3}  dy  , \\
	\mathbf{D}_{ijk } &=  \frac{1}{2} K_i K_j K_k \int_{0}^{1} \frac{1}{y}  C_{2 i+1}^{(2)}(y)
   \frac{d}{dy} \left[ \frac{  (1-y^2 )^{3/2} }{y}C_{2 j+1}^{(2)}(y) \right] C_{2 k+1}^{(2)}(y) (1-y^2)^{3/2} dy   , \\
   \mathbf{E}_{ijk } &= 
		\frac{1}{2}K_i K_j K_k \int_{0}^{1}  \frac{(1-2 y ^2) }{y}  C_{2 i+1}^{(2)}(y)  \frac{d}{dy} \left[ \frac{  (1-y^2 )^{3/2} }{y}C_{2 j+1}^{(2)}(y) \right]C_{2 k+1}^{(2)}(y)  (1-y^2)^{3/2}  dy ,
	\end{align*}
where  normalization constant $K_n$ is defined in Lemma \ref{LemmaEigenfunctionsWithRespectToGegenbauer}.
\end{lemma}
\begin{proof}
We set $y=\sin(x)$ and use Lemma \ref{LemmaEigenfunctionsWithRespectToGegenbauer}  to rewrite the Fourier coefficients as integrals with respect to $y \in [0,1]$. We obtain
\begin{align*}
%%%%%%%%%%%	
\mathbf{A}_{ij} &= \left(e_{i}e_{j} \big|\cos(2\cdot)\right) 
=\int_{0}^{\pi/2}e_{i}(x)e_{j}(x)\cos(2x) \tan^2(x)dx \\
  &=K_{i}K_{j}\int_{0}^{1} (1-y^2)^2 C_{2 i+1}^{(2)}(y)C_{2 j+1}^{(2)}(y) (1-2y^2)  \frac{dy}{(1-y^2)^{1/2}} \\
  &=K_{i}K_{j}\int_{0}^{1} (1-y^2)^{3/2} C_{2 i+1}^{(2)}(y)C_{2 j+1}^{(2)}(y) (1-2y^2)  dy, \\
%%%%%%%%%%%	
\mathbf{B}_{ij} &=\left(e_{i}^{\prime}e_{j} \big|\sin(2\cdot)\right) 
=\int_{0}^{\pi/2}e_{i}^{\prime}(x)e_{j}(x)\sin(2x) \tan^2(x)dx \\ 
&=K_i K_j  \int_{0}^{1}   \frac{d}{dy} \left[ \frac{  (1-y^2 )^{3/2} }{y}C_{2 i+1}^{(2)}(y) \right]\frac{(1-y^2)^{3/2}}{y}  C_{2 j+1}^{(2)}(y)  2y^3\frac{dy}{(1-y^2)^{1/2} } \\
&=2K_i K_j  \int_{0}^{1}   \frac{d}{dy} \left[ \frac{  (1-y^2 )^{3/2} }{y}C_{2 i+1}^{(2)}(y) \right] (1-y^2)   C_{2 j+1}^{(2)}(y)  y^2dy ,  \\  
%%%%%%%%%%%
	\mathbf{C}_{ijk} &= \left(e_{i}e_{j} \big|e_{k}\right) 
	= \int_{0}^{\pi/2}e_{i}(x)e_{j}(x)e_{k}(x)\tan^2(x)dx \\ 
	&    =  K_{i}K_{j}K_{k}\int_{0}^{1} \frac{(1-y^2)^{7/2}}{y}
	C_{2 i+1}^{(2)}(y)
	C_{2 j+1}^{(2)}(y)
	C_{2 k+1}^{(2)}(y)    \frac{dy}{(1-y^2)^{1/2}}  \\
	&    =  K_{i}K_{j}K_{k}\int_{0}^{1} \frac{(1-y^2)^{3}}{y}
	C_{2 i+1}^{(2)}(y)
	C_{2 j+1}^{(2)}(y)
	C_{2 k+1}^{(2)}(y)    dy,  \\
	%%%%%%%%%%%	
\mathbf{D}_{ijk } &= \left(e_{i}e_{j}^{\prime} \Big|\frac{e_k}{\sin (2 \cdot)}\right) = \int_{0}^{\pi/2}e_{i}(x)e_{j}^{\prime} (x)\frac{e_{k}(x)}{\sin  (2 x)} \tan^2(x)dx \\
&= K_i K_j K_k \int_{0}^{1} \frac{(1-y^2)^{3}}{y^2}  C_{2 i+1}^{(2)}(y)
(1-y^2)^{1/2}  \frac{d}{dy} \left[ \frac{  (1-y^2 )^{3/2} }{y}C_{2 j+1}^{(2)}(y) \right]
 C_{2 k+1}^{(2)}(y)\frac{y}{2 (1-y^2)^2} dy   \\
 &=\frac{1}{2} K_i K_j K_k \int_{0}^{1} \frac{(1-y^2)^{3/2}}{y}  C_{2 i+1}^{(2)}(y)
   \frac{d}{dy} \left[ \frac{  (1-y^2 )^{3/2} }{y}C_{2 j+1}^{(2)}(y) \right]
 C_{2 k+1}^{(2)}(y)dy  , \\
%%%%%%%%%%%	
\mathbf{E}_{ijk } &= \left(e_{i}e_{j}^{\prime} \Big|\frac{e_{k}}{\tan (2 \cdot)}\right)
= \int_{0}^{\pi/2}e_{i} (x)e_{j}^{\prime}(x)\frac{e_{k}(x)}{\tan (2 x)} \tan^2(x)dx  \\
&=K_i K_j K_k \int_{0}^{1}  \frac{(1-y^2)^{5/2}}{y}  C_{2 i+1}^{(2)}(y)  
(1-y^2)^{1/2}  \frac{d}{dy} \left[ \frac{  (1-y^2 )^{3/2} }{y}C_{2 j+1}^{(2)}(y) \right]
C_{2 k+1}^{(2)}(y)   \frac{ 1-2 y ^2}{2 (1-y^2)^{3/2}} dy  \\
&=\frac{1}{2}K_i K_j K_k \int_{0}^{1} (1-2 y ^2)  \frac{(1-y^2)^{3/2}}{y}  C_{2 i+1}^{(2)}(y)  
  \frac{d}{dy} \left[ \frac{  (1-y^2 )^{3/2} }{y}C_{2 j+1}^{(2)}(y) \right]
C_{2 k+1}^{(2)}(y)    dy  ,
\end{align*}
that completes the proof.
\end{proof}

\subsection{Linearization formulas}

We begin our analysis by deriving novel linearization formulas for weighted Gegenbauer polynomials (Lemmata \ref{Lemma1auxiliaryForFourierCijk},  \ref{Lemma1auxiliaryForFourierDijk} and \ref{Lemma1auxiliaryForFourierEijk}).

\begin{lemma}[Linearization formula for the Gegenbauer polynomials with inverse polynomial weight]\label{Lemma1auxiliaryForFourierCijk}
	Let $n \geq 0$ be an integer. We have that
	\begin{align}\label{Asxeto1FourierCijk}
		\frac{(1-y^2)^{3/2}}{y}
	C_{2 n+1}^{(2)}(y)=\sum_{l=0}^{\infty} \mathfrak{L}_{n}(l) C_{2l}^{(2)}(y) ,
	\end{align}
	where the coefficient is given by
	\begin{align*}
	\mathfrak{L}_{n}(l) &=   (-1)^{n} (2n+3)\frac{3072}{\pi} \Bigg( 
		\frac{  (l+1)}{   (2 l-3) (2 l-1) (2 l+1)^2 (2 l+3)^2 (2 l+5) (2 l+7)}
		 \\
		 &+(l+1) \sum_{p=1}^{n}  \frac{   (-1)^{p}(p+1)^2}{  
	 \prod_{j=0}^3 (2 l-2 p+2j-3) 
	 \prod_{j=0}^3 (2 l+2 p+2j+1) 
	 }  \Bigg).
	\end{align*}
\end{lemma}
\begin{proof}
	Let $n \geq 0$ be any integer. Since $C_{2 n+1}^{(2)}(0)=0$, the function at the left-hand-side of \eqref{Asxeto1FourierCijk} is clearly in $L^2 ((-1,1);(1-y^2)^{3/2}dy )$, hence there exists a sequence $\{\mathfrak{L}_{n}(l): l \geq 0\}$ so that \eqref{Asxeto1FourierCijk} holds. To find it in closed form, we proceed as follows. For any fix integer $m \geq 0$, we multiply both sides of \eqref{Asxeto1FourierCijk} with $ C_{2 m}^{(2)}(y) (1-y^2)^{3/2} $, integrate with respect to $y \in [-1,1]$ and use the orthogonality condition \eqref{AppendixOrthogonalityGegenbauerExponent2} to obtain
\begin{align*}  
	\mathfrak{L}_{n}(m) &=\frac{1}{\xi_{2m} }\int_{-1}^{1} \frac{(1-y^2)^{3}}{y}
	C_{2 n+1}^{(2)}(y)C_{2 m}^{(2)}(y) dy.
\end{align*}
According to the recurrence relation \eqref{AppendixRecrrenceRelation5}, we get that
\begin{align*}
	\mathfrak{L}_{n}(m) &=\frac{1}{\xi_{2m} }\int_{-1}^{1} \frac{(1-y^2)^{3}}{y}
	C_{2 n+1}^{(2)}(y)C_{2 m}^{(2)}(y) dy \\
	&=\frac{1}{\xi_{2m} }\int_{-1}^{1} \frac{(1-y^2)^{3}}{y}
	\left[\frac{4 n+4 }{2 n+1}y C_{2 n}^{(2)}(y)
	-\frac{2 n+3 }{2 n+1} C_{2 n-1}^{(2)}(y) \right]C_{2 m}^{(2)}(y) dy \\
	&=\frac{4 n+4 }{2 n+1}\frac{1}{\xi_{2m} }\int_{-1}^{1} (1-y^2)^{3}  
	 C_{2 n}^{(2)}(y) C_{2 m}^{(2)}(y) dy
	 -\frac{2 n+3 }{2 n+1}
	\frac{1}{\xi_{2m} }\int_{-1}^{1} \frac{(1-y^2)^{3}}{y}
	 C_{2 n-1}^{(2)}(y) C_{2 m}^{(2)}(y) dy \\
	 &=\frac{4 n+4 }{2 n+1}\frac{1}{\xi_{2m} }\int_{-1}^{1} (1-y^2)^{3}  
	 C_{2 n}^{(2)}(y) C_{2 m}^{(2)}(y) dy
	 -\frac{2 n+3 }{2 n+1} \mathfrak{L}_{n-1} (m).
\end{align*}
For the first term, we use the connection formula \eqref{AppendixConnectionFormulaExponent7half} together with the orthogonality of the Gegenbauer polynomials \eqref{AppendixOrthogonalityGegenbauerExponent7half} to infer
\begin{align*}
	& \int_{-1}^{1} (1-y^2)^{3}  
	 C_{2 n}^{(2)}(y) C_{2 m}^{(2)}(y) dy =
	 \sum _{l_1 =0}^n \eta^{(3)}_{n}(l_1)
	 \sum _{l_2 =0}^m \eta^{(3)}_{m}(l_2) \int_{-1}^{1} (1-y^2)^{3}  
	C_{2 l_1}^{\left( 7/2  \right)}(y)C_{2 l_2}^{\left( 7/2  \right)}(y)  dy \\
	 & = \sum _{l_1 =0}^n \eta^{(3)}_{n}(l_1)
	 \sum _{l_2 =0}^m \eta^{(3)}_{m}(l_2)\theta_{ 2l_1}  \mathds{1} (l_1=l_2) 
	  = \sum _{l =0}^{\min (n,m)} \eta^{(3)}_{n}(l)
	 \eta^{(3)}_{m}(l )\theta_{ 2l}  .
\end{align*}
Now, we run Gosper's algorithm and find a (very long) function $R_{nm}(l)$ in closed form so that $ \eta^{(3)}_{n}(l)
	 \eta^{(3)}_{m}(l )\theta_{ 2l}  = R_{nm}(l+1)-R_{nm}(l)$ and hence 
\begin{align*}
	\sum _{l =0}^{\min (n,m)} \eta^{(3)}_{n}(l)
	 \eta^{(3)}_{m}(l )\theta_{ 2l} 
	 % &=\sum _{l =0}^{\min (n,m)} \left(R_{nm}(l+1)-R_{nm}(l) \right) 
	 = R_{nm}(\min(n,m)+1)-R_{nm}(0).
\end{align*}
Since $R_{nm}(l)$ is given by a closed formula, we compute $R_{nm}(\min(n,m)+1)-R_{nm}(0)$ and obtain 
%\begin{align*}
%	R_{nm}(\min(n,m)+1)-R_{nm}(0)=\frac{96 (m+1) (2 m+1) (2 m+3) (n+1) (2 n+1) (2 n+3)}{
%	\prod_{j=0}^3  (2 m-2 n+2j-3) 
%	% (2 m-2 n-3) (2 m-2 n-1) (2 m-2 n+1) (2 m-2 n+3) 
%	\prod_{j=0}^3  (2 m+2 n+2j+1) 
%	% (2 m+2 n+1) (2 m+2 n+3) (2 m+2 n+5) (2 m+2 n+7)
%	}
%\end{align*}
%and hence
\begin{align*}
	\int_{-1}^{1} (1-y^2)^{3}  
	 C_{2 n}^{(2)}(y) C_{2 m}^{(2)}(y) dy= 
	 \frac{96 (m+1) (2 m+1) (2 m+3) (n+1) (2 n+1) (2 n+3)}{
	\prod_{j=0}^3  (2 m-2 n+2j-3)   
	\prod_{j=0}^3  (2 m+2 n+2j+1)  
	}.
\end{align*}
This leads us to the recurrence relation
\begin{align*}
	\mathfrak{L}_{n}(m)  +\frac{2 n+3 }{2 n+1} 
	\mathfrak{L}_{n-1} (m)
	 &=\frac{3072}{\pi} \frac{ (m+1) (n+1)^2 (2 n+3)}{  
	 \prod_{j=0}^3 (2 m-2 n+2j-3)  
	 \prod_{j=0}^3 (2 m+2 n+2j+1) 
	 }.
\end{align*}
Next, we find the initial condition to the recurrence relation above. Using the fact that $C_1^{(2)}(y)=4 y$ and the computation of the integral in \eqref{AppendixComputationIntegralHyperGeometricFunction1}, we get  
\begin{align*}
	\mathfrak{L}_{0}(m) &=\frac{1}{\xi_{2m} }\int_{-1}^{1} \frac{(1-y^2)^{3}}{y}
	C_{1}^{(2)}(y)C_{2 m}^{(2)}(y) dy 
	 = \frac{4}{\xi_{2m} }\int_{-1}^{1}(1-y^2)^{3}
	C_{2 m}^{(2)}(y) dy \\
	& =\frac{4}{\xi_{2m} } \sum _{l=0}^{2 m} \frac{96 \sqrt{\pi } (-1)^l (l+3) (l+2) (l+1)  \Gamma  (l+2 m+4)}{ \Gamma  \left(l+\frac{5}{2}\right)  \Gamma  (l+8)  \Gamma  (-l+2 m+1)}.
\end{align*}
For our purposes, it is necessary to compute the sum above in closed form. To do so, we run Zeilberger's algorithm and find the 1-step recurrence relation
\begin{align*}
	\mathfrak{L}_{0}(m+1) = \frac{(m+2) (2 m-3) (2 m+1)}{(m+1) (2 m+5) (2 m+9)} \mathfrak{L}_{0}(m)
\end{align*}
and the initial condition can be easily found,
\begin{align*}
	\mathfrak{L}_{0}(0) &= \frac{1}{\xi_{0} }\int_{-1}^{1} \frac{(1-y^2)^{3}}{y}
	C_{1}^{(2)}(y)C_{0}^{(2)}(y) dy 
	= \frac{4}{\xi_{0} }\int_{-1}^{1} (1-y^2)^{3}   dy=\frac{1024}{105 \pi },
\end{align*}
since $C_{0}^{(2)}(y)=1$ and $C_1^{(2)}(y)=4 y$. Now, the solution to the recurrence relation
\begin{align*}
	\begin{dcases}
	\mathfrak{L}_{0}(m+1) = \frac{(m+2) (2 m-3) (2 m+1)}{(m+1) (2 m+5) (2 m+9)} \mathfrak{L}_{0}(m),  \\
		\mathfrak{L}_{0}(0) = \frac{1024}{105 \pi }
	\end{dcases}
\end{align*}
is given by
\begin{align*}
	\mathfrak{L}_{0}(m) &= \frac{1024 }{105 \pi } \prod _{j=1}^{m} \frac{(j+1) (2 j-5) (2 j-1)}{j (2 j+3) (2 j+7)} 
	\\
	&=\frac{9216 (m+1)}{\pi  (2 m-3) (2 m-1) (2 m+1)^2 (2 m+3)^2 (2 m+5) (2 m+7)}.
\end{align*}
Consequently, we obtain the recurrence relation
\begin{align*}
\begin{dcases}
	\mathfrak{L}_{n}(m)  +\frac{2 n+3 }{2 n+1} 
	\mathfrak{L}_{n-1} (m)
	 =\frac{3072}{\pi} \frac{ (m+1) (n+1)^2 (2 n+3)}{  
	 \prod_{j=0}^3 (2 m-2 n+2j-3) 
	 \prod_{j=0}^3 (2 m+2 n+2j+1) 
	 },   \\
	 \mathfrak{L}_{0}(m)= \frac{9216 (m+1)}{\pi  (2 m-3) (2 m-1) (2 m+1)^2 (2 m+3)^2 (2 m+5) (2 m+7)}
\end{dcases}
\end{align*}
and its unique solution is given by (Lemma \ref{AppendixLemma1stOrderLinearRecurrence})
\begin{align*}
	\mathfrak{L}_{n}(m) = (-1)^{n+1} (2n+3) \left( -\frac{\mathfrak{L}_{0}(m) }{3}+\frac{3072}{\pi} \sum_{i=1}^{n}  \frac{ (m+1)  (-1)^{i+1}(i+1)^2}{  
	 \prod_{j=0}^3 (2 m-2 i+2j-3) 
	 \prod_{j=0}^3 (2 m+2 i+2j+1) 
	 }  \right),
\end{align*}
that completes the proof.
\end{proof}

\begin{lemma}[Linearization formulas for the Gegenbauer polynomials with inverse polynomial weights]\label{Lemma1auxiliaryForFourierDijk}
	Let $n \geq 0$ be an integer. We have that
	\begin{align}
	\frac{d}{dy} \left(	\frac{(1-y^2)^{3/2}}{y}
	C_{2 n+1}^{(2)}(y) \right) &=\sum_{r=0}^{\infty} \mathfrak{M}_{n}(r) C_{2r+1}^{(2)}(y) , \label{Asxeto1FourierDijk} \\
	\frac{1}{y}
	C_{2 n+1}^{(2)}(y) &=\sum_{s=0}^{n} \mathfrak{N}_{n}(s) C_{2s}^{(2)}(y) \label{Asxeto2FourierDijk}
	\end{align}
	where the coefficients are given by
	\begin{align*}
	\mathfrak{M}_{n}(r) &=(-1)^{n} (2 n+3)\frac{256}{\pi} \Bigg(   \frac{1}{  (2 r-1) (2 r+1) (2 r+3) (2 r+5) (2 r+7)}  \\
	&+   (2 r+3) \sum_{p=1}^{n} \frac{(-1)^{p} (p+1)^2 }{ \prod _{j=0}^2 (2 r-2 p+2 j-1)  \prod _{j=0}^2 (2 r+2p+2 j+3)} \Bigg)  
\end{align*}
	and
	\begin{align*}
	\mathfrak{N}_{n}(s) 
	&=  (2 n+3)   \frac{ 4(-1)^{s+n} (s+1)}{(2 s+1) (2 s+3)}    . 
\end{align*} 
\end{lemma}
\begin{proof}
	Let $n \geq 0$ be an integer. Since $C_{2 n+1}^{(2)}(0)=0$, both functions at the left-hand-sides of \eqref{Asxeto1FourierDijk} and \eqref{Asxeto2FourierDijk} are clearly in $L^2 ((-1,1);(1-y^2)^{3/2}dy )$, hence there exist sequences $\{\mathfrak{M}_{n}(r): r \geq 0\}$ and $\{\mathfrak{N}_{n}(l): l \geq 0\}$ so that \eqref{Asxeto1FourierDijk} and \eqref{Asxeto2FourierDijk} hold respectively and the series in \eqref{Asxeto1FourierDijk} converges. To find them in closed forms, we proceed as follows. For any fix integer $l \geq 0$, we multiply both sides of \eqref{Asxeto1FourierDijk} with $ C_{2 l+1}^{(2)}(y) (1-y^2)^{3/2} $ and \eqref{Asxeto2FourierDijk} with $ C_{2 l}^{(2)}(y) (1-y^2)^{3/2} $, integrate with respect to $y \in [-1,1]$ and use the orthogonality condition \eqref{AppendixOrthogonalityGegenbauerExponent2} to obtain
\begin{align*}  
	\mathfrak{M}_{n}(l) &=\frac{1}{\xi_{2l+1} }\int_{-1}^{1} \frac{d}{dy} \left(	\frac{(1-y^2)^{3/2}}{y}
	C_{2 n+1}^{(2)}(y) \right) C_{2 l+1}^{(2)}(y) (1-y^2)^{3/2} dy, \\
	\mathfrak{N}_{n}(l) &=\frac{1}{\xi_{2l} }\int_{-1}^{1}\frac{1}{y}
	C_{2 n+1}^{(2)}(y) C_{2 l}^{(2)}(y) (1-y^2)^{3/2} dy.
\end{align*}
Firstly, we focus on $\mathfrak{M}_{n}(l)$. Integration by parts and the recurrence relation \eqref{AppendixRecrrenceRelation4} yield
\begin{align*}
	\mathfrak{M}_{n}(l) &=\frac{1}{\xi_{2l+1} }\int_{-1}^{1} \frac{d}{dy} \left(	\frac{(1-y^2)^{3/2}}{y}
	C_{2 n+1}^{(2)}(y) \right) C_{2 l+1}^{(2)}(y) (1-y^2)^{3/2} dy \\
	&=-\frac{1}{\xi_{2l+1} }\int_{-1}^{1} 	\frac{(1-y^2)^{3/2}}{y}
	C_{2 n+1}^{(2)}(y)  \frac{d}{dy} \left( C_{2 l+1}^{(2)}(y) (1-y^2)^{3/2} \right) dy \\
	&=\frac{  2 (l+1) (l+2)}{\xi_{2l+1} }\int_{-1}^{1} 	\frac{(1-y^2)^{3/2}}{y}
	C_{2 n+1}^{(2)}(y)  \left(1-y^2\right)^{1/2}  C_{2( l+1)}^{(1)}(y) dy \\
	&=\frac{  2 (l+1) (l+2)}{\xi_{2l+1} }\int_{-1}^{1} 	\frac{(1-y^2)^{2}}{y}
	C_{2 n+1}^{(2)}(y)  C_{2( l+1)}^{(1)}(y) dy.
\end{align*}
Now, the recurrence relation \eqref{AppendixRecrrenceRelation5} yields
\begin{align*}
	\mathfrak{M}_{n}(l) &=\frac{  2 (l+1) (l+2)}{\xi_{2l+1} }\int_{-1}^{1} 	\frac{(1-y^2)^{2}}{y}
	C_{2 n+1}^{(2)}(y)  C_{2( l+1)}^{(1)}(y) dy \\
	&=\frac{  2 (l+1) (l+2)}{\xi_{2l+1} }\int_{-1}^{1} 	\frac{(1-y^2)^{2}}{y}
	\left[ \frac{4 n+4 }{2 n+1}y C_{2 n}^{(2)}(y)
	-\frac{2 n+3 }{2 n+1} C_{2 n-1}^{(2)}(y) \right]  C_{2( l+1)}^{(1)}(y) dy \\
	&= \frac{4 n+4 }{2 n+1} \frac{  2 (l+1) (l+2)}{\xi_{2l+1} }\int_{-1}^{1} 	 (1-y^2)^{2} 
	  C_{2 n}^{(2)}(y)  C_{2( l+1)}^{(1)}(y) dy \\
	& -\frac{2 n+3 }{2 n+1}
	  \frac{  2 (l+1) (l+2)}{\xi_{2l+1} }\int_{-1}^{1} 	\frac{(1-y^2)^{2}}{y} 
	 C_{2 n-1}^{(2)}(y)  C_{2( l+1)}^{(1)}(y) dy \\
	 &= \frac{4 n+4 }{2 n+1} \frac{  2 (l+1) (l+2)}{\xi_{2l+1} }\int_{-1}^{1} 	 (1-y^2)^{2} 
	  C_{2 n}^{(2)}(y)  C_{2( l+1)}^{(1)}(y) dy 
	   -\frac{2 n+3 }{2 n+1}\mathfrak{M}_{n-1}(l).
\end{align*}
 For the first term, we use the connection formulas \eqref{AppendixConnectionFormulaExponent5half} and \eqref{AppendixConnectionFormulaExponent5halfAgain} together with the orthogonality of the Gegenbauer polynomials \eqref{AppendixOrthogonalityGegenbauerExponent5half} to infer
\begin{align*}
	&\int_{-1}^{1} 	 (1-y^2)^{2} 
	  C_{2 n}^{(2)}(y)  C_{2( l+1)}^{(1)}(y) dy = 
	  \sum _{\mu_1=0}^n \eta_{n}^{(1)} (\mu_1 ) \sum _{\mu_2 =0}^{l+1}\eta_{l}^{(2)} (\mu_2 ) \int_{-1}^{1} 	 (1-y^2)^{2} 
	 C_{2 \mu_1}^{\left( 5/2  \right)}(y)  C_{2 \mu_2 }^{\left( 5/2  \right)}(y) dy  \\
	 & =  \sum _{\mu_1=0}^n \eta_{n}^{(1)} (\mu_1 ) \sum _{\mu_2 =0}^{l+1}\eta_{l}^{(2)} (\mu_2 ) \rho_{2\mu_1} \mathds{1}(\mu_1=\mu_2) 
	  =  \sum _{\mu =0}^{\min (n,l+1)} \eta_{n}^{(1)} (\mu )  \eta_{l}^{(2)} (\mu ) \rho_{2\mu} .
\end{align*}
As in Lemma \ref{Lemma1auxiliaryForFourierCijk}, we run Gosper's algorithm and find a (very long) function $R_{nl}(\mu)$ in closed form so that $\eta_{n}^{(1)} (\mu )  \eta_{l}^{(2)} (\mu ) \rho_{2\mu}  =R_{nl}(\mu+1)-R_{nl}(\mu)$, hence 
\begin{align*}
	\sum _{\mu =0}^{\min (n,l+1)} \eta_{n}^{(1)} (\mu )  \eta_{l}^{(2)} (\mu ) \rho_{2\mu}=  R_{nl}(\min(n,l+1)+1)-R_{nl}(0).
\end{align*}
Since $R_{nl}(\mu)$ is given by a closed formula, we compute $R_{nl}(\min(n,l+1)+1)-R_{nl}(0)$ and obtain  
\begin{align*}
	\int_{-1}^{1} 	 (1-y^2)^{2} 
	  C_{2 n}^{(2)}(y)  C_{2( l+1)}^{(1)}(y) dy=\frac{16 (2 l+3) (n+1) (2 n+1) (2 n+3)}{ \prod _{j=0}^2 (2 j+2 l-2 n-1)  \prod _{j=0}^2 (2 j+2 l+2 n+3)}.
\end{align*}
This leads us to the recurrence relation
\begin{align*}
	\mathfrak{M}_{n}(l)+\frac{2 n+3 }{2 n+1}\mathfrak{M}_{n-1}(l) &=\frac{256}{\pi}  \frac{ (2 l+3) (n+1)^2 (2 n+3)}{   \prod _{j=0}^2 (2 j+2 l-2 n-1)  \prod _{j=0}^2 (2 j+2 l+2 n+3)}.
\end{align*}
 Next, we find the initial condition to the recurrence relation above. Using the fact that $C_1^{(2)}(y)=4 y$ and the computation of the integral in \eqref{AppendixComputationIntegralHyperGeometricFunction2}, we get 
\begin{align*}
	\mathfrak{M}_{0}(l) 
	&=\frac{  2 (l+1) (l+2)}{\xi_{2l+1} }\int_{-1}^{1} 	\frac{(1-y^2)^{2}}{y}
	C_{1}^{(2)}(y)  C_{2( l+1)}^{(1)}(y) dy \\
	&=\frac{  8 (l+1) (l+2)}{\xi_{2l+1} }\int_{-1}^{1}  (1-y^2)^{2} 
	  C_{2( l+1)}^{(1)}(y) dy \\
	&=\frac{  2 (l+1) (l+2)}{\xi_{2l+1} } \sum_{\mu=0}^{2(l+1)} \frac{(-1)^{\mu } 4^{\mu +3}  \Gamma  (2 l+\mu +4)}{(\mu +3) (\mu +4) (\mu +5)  \Gamma  (2 \mu +2)  \Gamma  (2 l-\mu +3)}.    
\end{align*}
For our purposes, it is necessary to compute the sum above in closed form. To do so, we run Zeilberger's algorithm and find the 1-step recurrence relation
\begin{align*}
	\mathfrak{M}_{0}(l+1) =  \frac{2 l-1}{2 l+9} \mathfrak{M}_{0}(l)
\end{align*}
and the initial condition can be easily found,
\begin{align*}
	\mathfrak{M}_{0}(0) = \frac{  4}{\xi_{1} }\int_{-1}^{1} 	\frac{(1-y^2)^{2}}{y}
	C_{1}^{(2)}(y)  C_{2}^{(1)}(y) dy =
	\frac{16}{\pi} \int_{-1}^{1} \left(y^2-1\right)^2 \left(4 y^2-1\right) dy=-\frac{256}{35 \pi },
\end{align*}
since $C_1^{(2)}(y)=4 y$ and $C_{2}^{(1)}(y)=4 y^2-1$. Now, the solution to the recurrence relation
\begin{align*}
	\begin{dcases}
	\mathfrak{M}_{0}(l+1) =  \frac{2 l-1}{2 l+9} \mathfrak{M}_{0}(l),  \\
		\mathfrak{M}_{0}(0) =-\frac{256}{35 \pi }
	\end{dcases}
\end{align*}
is given by
\begin{align*} 
	\mathfrak{M}_{0}(l) &=  -\frac{256}{35 \pi } \prod _{j=0}^{l-1} \frac{2 j-1}{2 j+9}=\frac{768}{\pi  (2 l-1) (2 l+1) (2 l+3) (2 l+5) (2 l+7)}.
\end{align*}
Consequently, we obtain the recurrence relation 
\begin{align*}
\begin{dcases}
	\mathfrak{M}_{n}(l)+\frac{2 n+3 }{2 n+1}\mathfrak{M}_{n-1}(l) =\frac{256}{\pi}  \frac{ (2 l+3) (n+1)^2 (2 n+3)}{   \prod _{j=0}^2 (2 j+2 l-2 n-1)  \prod _{j=0}^2 (2 j+2 l+2 n+3)},   \\
	 \mathfrak{M}_{0}(l)=\frac{768}{\pi  (2 l-1) (2 l+1) (2 l+3) (2 l+5) (2 l+7)}
\end{dcases}
\end{align*}
and its unique solution is given by (Lemma \ref{AppendixLemma1stOrderLinearRecurrence})
 \begin{align*}
	\mathfrak{M}_{n}(l) =(-1)^{n+1} (2 n+3) \left(  -\frac{\mathfrak{M}_{0}(l)}{ 3 }  +  \frac{256}{\pi} \sum_{i=1}^{n} \frac{(-1)^{i+1} (i+1)^2 (2 l+3)}{ \prod _{j=0}^2 (-2 i+2 j+2 l-1)  \prod _{j=0}^2 (2 i+2 j+2 l+3)} \right)  .
\end{align*} 
that completes the proof for $\mathfrak{M}_{n}(l)$. For $\mathfrak{N}_{n}(l)$, the recurrence relation \eqref{AppendixRecrrenceRelation5} and the orthogonality condition \eqref{AppendixOrthogonalityGegenbauerExponent2} yield
\begin{align*}
	\mathfrak{N}_{n}(l) &=\frac{1}{\xi_{2l} }\int_{-1}^{1}\frac{1}{y}
	C_{2 n+1}^{(2)}(y) C_{2 l}^{(2)}(y) (1-y^2)^{3/2} dy \\
	&=\frac{1}{\xi_{2l} }\int_{-1}^{1}\frac{1}{y}
	\left[ \frac{4 n+4 }{2 n+1}y C_{2 n}^{(2)}(y)
	-\frac{2 n+3 }{2 n+1} C_{2 n-1}^{(2)}(y) \right] C_{2 l}^{(2)}(y) (1-y^2)^{3/2} dy \\
	&=\frac{4 n+4 }{2 n+1} \frac{1}{\xi_{2l} }\int_{-1}^{1} 
	    C_{2 n}^{(2)}(y) C_{2 l}^{(2)}(y) (1-y^2)^{3/2} dy -\frac{2 n+3 }{2 n+1} \frac{1}{\xi_{2l} }\int_{-1}^{1}\frac{1}{y} C_{2 n-1}^{(2)}(y)   C_{2 l}^{(2)}(y) (1-y^2)^{3/2} dy \\
	    &=\frac{4 n+4 }{2 n+1} \mathds{1}(n=l) -\frac{2 n+3 }{2 n+1} \mathfrak{N}_{n-1}(l)
\end{align*}
coupled to the initial condition
\begin{align*}
	\mathfrak{N}_{0}(l) &=\frac{1}{\xi_{2l} }\int_{-1}^{1}\frac{1}{y}
	C_{1}^{(2)}(y) C_{2 l}^{(2)}(y) (1-y^2)^{3/2} dy 
	=\frac{4 }{\xi_{2l} }\int_{-1}^{1} C_{2 l}^{(2)}(y) (1-y^2)^{3/2} dy \\
	&=\frac{4 }{\xi_{2l} }\mathds{1}(l=0)\int_{-1}^{1}  (1-y^2)^{3/2} dy =4\mathds{1}(l=0),
\end{align*}
since $C_1^{(2)}(y)=4 y$ and $C_0^{(2)}(y)=1$. Here, we also used the fact that the set $\{C_{l}^{(2)}(y): l\geq 0\}$ forms an orthogonal and complete basis for the $L^2 ( (-1,1);(1-y^2)^{3/2}dy  )$ (see also relation 1, section 7.311 in \cite{MR3307944}). Consequently, we obtain the recurrence relation 
\begin{align*}
\begin{dcases}
	\mathfrak{N}_{n}(l)+\frac{2 n+3 }{2 n+1} \mathfrak{N}_{n-1}(l)=\frac{4 n+4 }{2 n+1} \mathds{1}(n=l),   \\
	 \mathfrak{N}_{0}(l)=4\mathds{1}(l=0)
\end{dcases}
\end{align*}
and its unique solution is given by (Lemma \ref{AppendixLemma1stOrderLinearRecurrence})
\begin{align*}
	\mathfrak{N}_{n}(l) &=(-1)^{n} (2 n+3) \left(  \frac{\mathfrak{N}_{0}(l)}{ 3 }  +  \sum_{i=1}^{n} \frac{(-1)^{i}}{ (2 i+3)}\frac{4 i+4 }{2 i+1} \mathds{1}(i=l) \right)  \\
	&=4(-1)^{n} (2 n+3)  \left(  \frac{1}{ 3 }\mathds{1}(l=0)  +\frac{ (-1)^l (l+1)}{(2 l+1) (2 l+3)}\mathds{1}(1\leq l \leq n)    \right)  \\
	& =4(-1)^{n} (2 n+3)   \frac{ (-1)^l (l+1)}{(2 l+1) (2 l+3)} \mathds{1}(0\leq l \leq n)    ,
\end{align*} 
that completes the proof. 
\end{proof}

\begin{lemma}[Linearization formula for the Gegenbauer polynomials with weight]\label{Lemma1auxiliaryForFourierEijk}
	Let $n \geq 0$ be an integer. We have that
	\begin{align}\label{Asxeto1FourierEijk} 
		(1-y^2) ^{1/2} C_{2(n+1)}^{(1)}(y) =\sum_{l=0}^{\infty} \mathfrak{P}_{n}(l) C_{2l}^{(2)}(y) ,
	\end{align}
	where the coefficient is given by
	\begin{align*}
	\mathfrak{P}_{n}(l) &= \frac{128}{\pi } \frac{(l+1) (2 n+3)}{ 
	\prod _{r=0}^2 (2 n-2 l+2 r-1)  
	\prod _{r=0}^2 (2 n+2 l+2 r+3)} .
	\end{align*}
\end{lemma}
\begin{proof}
	Let $n \geq 0$ be an integer. Clearly, the left-hand-side of \eqref{Asxeto1FourierEijk} is in $L^2 ( (-1,1);(1-y^2)^{3/2}dy  )$, hence there exists a sequence $\{  \mathfrak{P}_{n}(l)  : l \geq 0\}$ such that \eqref{Asxeto1FourierEijk} holds and the series converges. To find it in closed form, we proceed as follows. For any fix integer $m \geq 0$, we multiply both sides of \eqref{Asxeto1FourierEijk} with $ C_{2 m}^{(2)}(y) (1-y^2)^{3/2} $, integrate with respect to $y \in [-1,1]$ and use the orthogonality condition \eqref{AppendixOrthogonalityGegenbauerExponent2} to obtain
	\begin{align*}
		 \mathfrak{P}_{n}(m)  = \frac{1}{\xi_{2m}}\int_{-1}^{1}  C_{2(n+1)}^{(1)}(y)C_{2 m}^{(2)}(y) (1-y^2)^{2} dy .
	\end{align*}
	To compute this integral, we use the connection formulas \eqref{AppendixConnectionFormulaExponent5half} and \eqref{AppendixConnectionFormulaExponent5halfAgain} together with the orthogonality of the Gegenbauer polynomials \eqref{AppendixOrthogonalityGegenbauerExponent5half} to infer
\begin{align*}
	& \int_{-1}^{1} 	 (1-y^2)^{2} 
	  C_{2 m}^{(2)}(y)  C_{2( n+1)}^{(1)}(y) dy = 
	  \sum _{\mu_1=0}^m \eta_{m}^{(1)} (\mu_1 ) \sum _{\mu_2 =0}^{n+1}\eta_{n}^{(2)} (\mu_2 ) \int_{-1}^{1} 	 (1-y^2)^{2} 
	 C_{2 \mu_1}^{\left( 5/2  \right)}(y)  C_{2 \mu_2 }^{\left( 5/2  \right)}(y) dy  \\
	 & =  \sum _{\mu_1=0}^m \eta_{m}^{(1)} (\mu_1 ) \sum _{\mu_2 =0}^{n+1}\eta_{n}^{(2)} (\mu_2 ) \rho_{2\mu_1} \mathds{1}(\mu_1=\mu_2)
	  =  \sum _{\mu =0}^{\min (m,n+1)} \eta_{m}^{(1)} (\mu )  \eta_{n}^{(2)} (\mu ) \rho_{2\mu} .
\end{align*}
As in Lemma \ref{Lemma1auxiliaryForFourierCijk},  we run Gosper's algorithm and find a (very long) function $R_{mn}(\mu)$ in closed form so that $\eta_{m}^{(1)} (\mu )  \eta_{n}^{(2)} (\mu ) \rho_{2\mu}  =R_{mn}(\mu+1)-R_{mn}(\mu)$, hence 
\begin{align*}
	\sum _{\mu =0}^{\min (m,n+1)} \eta_{m}^{(1)} (\mu )  \eta_{n}^{(2)} (\mu ) \rho_{2\mu}=
	% \sum _{\mu =0}^{\min (m,n+1)} \left(R_{mn}(\mu+1)-R_{mn}(\mu) \right)= 
	R_{mn}(\min(m,n+1)+1)-R_{mn}(0).
\end{align*}
Since $R_{mn}(\mu)$ is given by a closed formula, we compute $R_{mn}(\min(m,n+1)+1)-R_{mn}(0)$ and obtain  
\begin{align*}
	\int_{-1}^{1} 	 (1-y^2)^{2} 
	  C_{2 m}^{(2)}(y)  C_{2( n+1)}^{(1)}(y) dy=\frac{16 (2 n+3) (m+1) (2 m+1) (2 m+3)}{ \prod _{j=0}^2 (2 j+2 n-2 m-1)  \prod _{j=0}^2 (2 j+2 n+2 m+3)}
\end{align*}
that completes the proof. 
\end{proof}

\subsection{Computationally efficient formulas} 

Next, we use the linearization formulas for the Gegenbauer polynomials (Lemmata \ref{Lemma1auxiliaryForFourierCijk},  \ref{Lemma1auxiliaryForFourierDijk} and \ref{Lemma1auxiliaryForFourierEijk}) to establish computationally efficient formulas for the Fourier coefficients (Propositions \ref{ComputationallyEfficientFormulaAij}-\ref{ComputationallyEfficientFormulaEijk}). To begin with, we focus on the Fourier coefficient $\mathbf{A}_{ij}$.

\begin{prop}[Computationally efficient formula for $\mathbf{A}_{ij}$]\label{ComputationallyEfficientFormulaAij}
	Let $i,j \geq 0$ be integers. Then, the Fourier coefficient $\mathbf{A}_{ij}$ is given by
	\begin{align*}
		\mathbf{A}_{ij} &=  \frac{\pi}{8} K_{i}K_{j} \Big[  \mathds{1}(i=j) -    j (j+2) \mathds{1}(i=j-1) -    (j+1) (j+3) \mathds{1}(i=j+1)  \Big],
	\end{align*} 
where  normalization constant $K_n$ is defined in Lemma \ref{LemmaEigenfunctionsWithRespectToGegenbauer}. 
	\end{prop}
\begin{proof}
	Let $i,j \geq 0$ be integers. Then, according to Lemma \ref{LemmaFourierCoefficientsWithRespectToGegenbauer}, the symmetries of the Gegenbauer polynomials \eqref{AppendixSymmetriesGegenbauer} and the orthogonality condition \eqref{AppendixOrthogonalityGegenbauerExponent2}, we have that
	\begin{align*}
		\mathbf{A}_{ij} &=  K_{i}K_{j}\int_{0}^{1} (1-y^2)^{3/2} C_{2 i+1}^{(2)}(y)C_{2 j+1}^{(2)}(y) (1-2y^2)  dy \\
		&=  \frac{1}{2}K_{i}K_{j}\int_{-1}^{1} (1-y^2)^{3/2} C_{2 i+1}^{(2)}(y)C_{2 j+1}^{(2)}(y) (1-2y^2)  dy \\
		&= \frac{1}{2} K_{i}K_{j}\int_{-1}^{1} (1-y^2)^{3/2} C_{2 i+1}^{(2)}(y)C_{2 j+1}^{(2)}(y) (2(1-y^2)-1)  dy \\
		&=  K_{i}K_{j}\int_{-1}^{1} (1-y^2)^{5/2} C_{2 i+1}^{(2)}(y)C_{2 j+1}^{(2)}(y)   dy - \frac{1}{2}K_{i}K_{j}\int_{-1}^{1} (1-y^2)^{3/2} C_{2 i+1}^{(2)}(y)C_{2 j+1}^{(2)}(y)  dy \\ 
		&= K_{i}K_{j}\int_{-1}^{1} (1-y^2)^{5/2} C_{2 i+1}^{(2)}(y)C_{2 j+1}^{(2)}(y)   dy - \frac{1}{2}K_{i}K_{j}\xi _{2i+1} \mathds{1}(i=j).
	\end{align*}
	To compute the integral, we use the recurrence relation \eqref{AppendixRecrrenceRelation1} together with the orthogonality condition \eqref{AppendixOrthogonalityGegenbauerExponent3} to obtain
	\begin{align*}
		& \int_{-1}^{1} (1-y^2)^{5/2} C_{2 i+1}^{(2)}(y) C_{2 j+1}^{(2)}(y) dy 
		= \frac{4 }{(2 i+3) (2 j+3)} \Bigg[ 
		\int_{-1}^{1} (1-y^2)^{5/2} C_{2 i-1}^{(3)}(y) C_{2 j-1}^{(3)}(y) dy \\
		& - \int_{-1}^{1} (1-y^2)^{5/2}  C_{2 i+1}^{(3)}(y) C_{2 j-1}^{(3)}(y)dy
		-\int_{-1}^{1} (1-y^2)^{5/2}  C_{2 i-1}^{(3)}(y) C_{2 j+1}^{(3)}(y) dy \\
		& + \int_{-1}^{1} (1-y^2)^{5/2} C_{2 i+1}^{(3)}(y) C_{2 j+1}^{(3)}(y) dy\Bigg] \\
		&= \frac{4 }{(2 i+3) (2 j+3)} \Big[ 
		\zeta_{2 i-1}\mathds{1}(i=j)  - \zeta_{2 i+1}\mathds{1}(i=j-1) 
		-\zeta_{2 i-1}\mathds{1}(i=j+1)   + \zeta_{2 i+1}\mathds{1}(i=j) \Big]
	\end{align*}
and putting all together yields the claim.	
\end{proof} 
\noindent
Secondly, we establish a computationally efficient formula for $\mathbf{B}_{ij}$.

\begin{prop}[Computationally efficient formula for $\mathbf{B}_{ij}$]\label{ComputationallyEfficientFormulaBij}
	Let $i,j \geq 0$ be integers. Then, the Fourier coefficient $\mathbf{B}_{ij}$ is given by
	\begin{align*}
		\mathbf{B}_{ij} &= \frac{\pi}{8} K_{i}K_{j} \Big[ -(3+2j)^2 \mathds{1}(i=j) -    j (j+2)(2j+1) \mathds{1}(i=j-1) +    (j+1) (j+3) (2j+5) \mathds{1}(i=j+1)  \Big],
	\end{align*} 
where  normalization constant $K_n$ is defined in Lemma \ref{LemmaEigenfunctionsWithRespectToGegenbauer}. 
	\end{prop}
\begin{proof}
	Let $i,j \geq 0$ be integers.  Then, according to Lemma \ref{LemmaFourierCoefficientsWithRespectToGegenbauer} and the recurrence relation \eqref{AppendixRecrrenceRelation4} we have that
	\begin{align*}
		\mathbf{B}_{ij} &= 
2K_i K_j  \int_{0}^{1}  y^2 \frac{d}{dy} \left[ \frac{  (1-y^2 )^{3/2} }{y}C_{2 i+1}^{(2)}(y) \right] (1-y^2)   C_{2 j+1}^{(2)}(y)  dy \\
&=2K_i K_j  \int_{0}^{1} y^2 \left( \frac{1}{y}\frac{d}{dy} \left[    (1-y^2 )^{3/2} C_{2 i+1}^{(2)}(y) \right] - \frac{1}{y^2} (1-y^2 )^{3/2} C_{2 i+1}^{(2)}(y)   \right) (1-y^2)   C_{2 j+1}^{(2)}(y)  y^2dy \\
&=2K_i K_j  \int_{0}^{1}  \left( -2y (i+1) (i+2)   (1-y^2 )^{1/2}  C_{2( i+1)}^{(1)}(y) -   (1-y^2 )^{3/2} C_{2 i+1}^{(2)}(y)   \right) (1-y^2)   C_{2 j+1}^{(2)}(y) dy \\
&=-2\Lambda_{i}K_{j}\int_{0}^{1}  
 y  (1-y^2)^{3/2}     
C_{2i+2}^{(1)}(y)    C_{2 j+1}^{(2)}(y)  dy -2K_iK_{j}\int_{0}^{1}  
 (1-y^2)^{5/2}    
C_{2 i+1}^{(2)}(y)  C_{2 j+1}^{(2)}(y)   dy, 
	\end{align*}
where we set
\begin{align*}
	\Lambda_{i} = \frac{ 4 (-1)^i  \sqrt{(i+2) (i+1)}}{\sqrt{\pi }},
\end{align*}	
and the symmetries of the Gegenbauer polynomials \eqref{AppendixSymmetriesGegenbauer} yield
	\begin{align*}
		\mathbf{B}_{ij}  
&= -K_iK_{j}\int_{-1}^{1}  
 (1-y^2)^{5/2}    
C_{2 i+1}^{(2)}(y)  C_{2 j+1}^{(2)}(y)   dy 
-\Lambda_{i}K_{j}\int_{-1}^{1}  
 y  (1-y^2)^{3/2}     
C_{2i+2}^{(1)}(y)    C_{2 j+1}^{(2)}(y)  dy .
	\end{align*}
In Lemma \ref{ComputationallyEfficientFormulaAij}, we already showed that
	\begin{align*}
		 \int_{-1}^{1} (1-y^2)^{5/2} C_{2 i+1}^{(2)}(y) C_{2 j+1}^{(2)}(y) dy 
	&= \frac{4 }{(2 i+3) (2 j+3)} \Big[ 
		\zeta_{2 i-1}\mathds{1}(i=j)  - \zeta_{2 i+1}\mathds{1}(i=j-1) \\
		&
		-\zeta_{2 i-1}\mathds{1}(i=j+1)   + \zeta_{2 i+1}\mathds{1}(i=j) \Big].
	\end{align*}
To compute the second integral, we use the recurrence relation \eqref{AppendixRecrrenceRelation2} to obtain
\begin{align*}
	& \int_{-1}^{1}  
 y  (1-y^2)^{3/2}     
C_{2i+2}^{(1)}(y)    C_{2 j+1}^{(2)}(y)  dy = \\
& \int_{-1}^{1}  
 y  (1-y^2)^{1/2}     
C_{2i+2}^{(1)}(y)   (1-y^2) C_{2 j+1}^{(2)}(y)  dy = \\
& \int_{-1}^{1}  
 y  (1-y^2)^{1/2}     
C_{2i+2}^{(1)}(y)  \left(\frac{j+2}{2}  C_{2 j+1}^{(1)}(y)-
\frac{j+1}{2}  C_{2 j+3}^{(1)}(y)  \right)  dy = \\
& \frac{j+2}{2}\int_{-1}^{1}  
 y (1-y^2)^{1/2}     
C_{2i+2}^{(1)}(y)   C_{2 j+1}^{(1)}(y) dy -
\frac{j+1}{2} 
 \int_{-1}^{1}  
 y (1-y^2)^{1/2}     
C_{2i+2}^{(1)}(y)   C_{2 j+3}^{(1)}(y)   dy  .
\end{align*}	
Now, the recurrence relation	 \eqref{AppendixRecrrenceRelation3} together with the orthogonality condition \eqref{AppendixOrthogonalityGegenbauerExponent1} yield
\begin{align*}
	 \int_{-1}^{1}  
 y (1-y^2)^{1/2}     
C_{2i+2}^{(1)}(y)   C_{2 j+1}^{(1)}(y) dy & = 
\frac{1}{2}  \int_{-1}^{1}  
 (1-y^2)^{1/2}     
 \left(  C_{2 i+1}^{(1)}(y)+C_{2 i+3}^{(1)}(y) \right)  C_{2 j+1}^{(1)}(y) dy \\
 & = \frac{\pi }{4} \left( \mathds{1}(i=j)  + \mathds{1}(i=j-1)  \right) ,\\
  \int_{-1}^{1}  
 y (1-y^2)^{1/2}     
C_{2i+2}^{(1)}(y)   C_{2 j+3}^{(1)}(y) dy & = 
\frac{1}{2}  \int_{-1}^{1}  
 (1-y^2)^{1/2}     
 \left(  C_{2 i+1}^{(1)}(y)+C_{2 i+3}^{(1)}(y) \right)  C_{2 j+3}^{(1)}(y) dy \\
 & = \frac{\pi }{4} \left( \mathds{1}(i=j+1)  + \mathds{1}(i=j)  \right) 
\end{align*}	
and putting all together yields the claim.	 
\end{proof}
\noindent
Next, we use the linearization formula of Lemma \ref{Lemma1auxiliaryForFourierCijk} to establish a computationally efficient formula for $\mathbf{C}_{ijk}$. Notice that, by its definition (Lemma \ref{LemmaFourierCoefficientsWithRespectToGegenbauer}), the Fourier coefficient $\mathbf{C}_{ijk}$ is invariant under permutations of the indices $i$, $j$ and $k$. Hence,  no generality is lost in Proposition \ref{ComputationallyEfficientFormulaCijk} by making the simplifying assumption $i \leq j \leq k$.

\begin{prop}[Computationally efficient formula for $\mathbf{C}_{ijk}$]\label{ComputationallyEfficientFormulaCijk}
	Let $i,j,k \geq 0$ be integers with $i \leq j \leq k$. Then, the Fourier coefficient $\mathbf{C}_{ijk}$ defined in Lemma \ref{LemmaEinsteinMassLessFourier} is given by
	\begin{align*}
		\mathbf{C}_{ijk} &=  \frac{1}{2} K_{i}K_{j}K_{k} \mathbf{w} _{ijk}\left(    \frac{1}{\mathbf{f}_{jk} }+ \sum_{p=1}^{i} \frac{(-1)^{p} (p+1)^2 }{\mathbf{g}_{jk}(p)} \right),
	\end{align*}
	where the normalization constant $K_n$ is defined in Lemma \ref{LemmaEigenfunctionsWithRespectToGegenbauer} and the coefficients are given by
	\begin{align*}
		\mathbf{w} _{ijk} &= 512 (-1)^i (2 i+3) (j+1) (j+2) (2 j+3) (k+1) (k+2) (2 k+3), \\
		 \mathbf{f}_{jk} &=  (2 j-2 k-3) (2 j-2 k-1) (2 j-2 k+1) (2 j-2 k+3)  \cdot \\
		 & \cdot (2 j+2 k+3) (2 j+2 k+5) (2 j+2 k+7) (2 j+2 k+9), \\
		 \mathbf{g}_{jk}(p) &= (2 (j+k-p)+3) (2 (j+k-p)+5) (2 j-2 k-2 p-3) (2 j-2 k-2 p-1) \cdot \\
		 & \cdot (2 (j-k+p)+1) (2 (j-k+p)+3) (2 (j+k+p)+7) (2 (j+k+p)+9).
	\end{align*}  
	\end{prop} 	
\begin{proof}
	Let $i,j,k \geq 0$ be integers with $i \leq j \leq k$. Then, according to Lemmata \ref{LemmaFourierCoefficientsWithRespectToGegenbauer} and \ref{Lemma1auxiliaryForFourierCijk}, the symmetries of the Gegenbauer polynomials \eqref{AppendixSymmetriesGegenbauer} and Dougall's formula \eqref{AppendixDougallFormula}, we have that
	\begin{align*}
		\mathbf{C}_{ijk} &= K_{i}K_{j}K_{k}\int_{0}^{1} \frac{(1-y^2)^{3}}{y}
	C_{2 i+1}^{(2)}(y)
	C_{2 j+1}^{(2)}(y)
	C_{2 k+1}^{(2)}(y)    dy \\
	&= \frac{1}{2}K_{i}K_{j}K_{k}\int_{-1}^{1} \frac{(1-y^2)^{3}}{y}
	C_{2 i+1}^{(2)}(y)
	C_{2 j+1}^{(2)}(y)
	C_{2 k+1}^{(2)}(y)    dy \\
	&=\frac{1}{2}K_{i}K_{j}K_{k} \int_{-1}^{1} \frac{(1-y^2)^{3/2}}{y}
	C_{2 i+1}^{(2)}(y)
	C_{2 j+1}^{(2)}(y)
	C_{2 k+1}^{(2)}(y)  (1-y^2)^{3/2}  dy \\
	&= \frac{1}{2}K_{i}K_{j}K_{k}\sum_{l=0}^{\infty} \mathfrak{L}_{i}(l) \int_{-1}^{1} C_{2l}^{(2)}(y)
	C_{2 j+1}^{(2)}(y)
	C_{2 k+1}^{(2)}(y)  (1-y^2)^{3/2}  dy \\ 
	&=\frac{1}{2} K_{i}K_{j}K_{k}\sum_{l=0}^{\infty} \mathfrak{L}_{i}(l) D(2l,2j+1,2k+1) \mathds{1}( 2|j-k|\leq 2l \leq 2(j+k+1) ) \\ 
	&= \frac{1}{2} K_{i}K_{j}K_{k}\sum_{l=|j-k|}^{j+k+1} \mathfrak{L}_{i}(l) D(2l,2j+1,2k+1).
	\end{align*}
Now, using the closed formulas for $\mathfrak{L}_{i}(l)$ and $D(2l,2j+1,2k+1)$, we compute
\begin{align*}
		\mathfrak{L}_{i}(l) D(2l,2j+1,2k+1) &=  \mathfrak{F}_{l}(i,j,k) +\sum_{p=1}^{i} \mathfrak{G}_{pl}(i,j,k),
	\end{align*}
for some functions $ \mathfrak{F}_{l}(i,j,k) $ and $\mathfrak{G}_{pl}(i,j,k)$ given by explicit formulas. With this notation, we need to compute
\begin{align*}
	\sum_{l=|j-k|}^{j+k+1} \mathfrak{L}_{i}(l) D(2l,2j+1,2k+1)=
	\sum_{l=|j-k|}^{j+k+1}  \mathfrak{F}_{l}(i,j,k) +
	\sum_{p=1}^{i} \left(\sum_{l=|j-k|}^{j+k+1}  \mathfrak{G}_{pl}(i,j,k) \right).
\end{align*} 
To do so, we use Gosper's algorithm and compute both sums
\begin{align*}
	\sum_{l=|j-k|}^{j+k+1}  \mathfrak{F}_{l}(i,j,k), \quad \sum_{l=|j-k|}^{j+k+1}  \mathfrak{G}_{pl}(i,j,k) .
\end{align*}
%Specifically, we find functions $R^{\mathfrak{F}}_{l}(i,j,k)$ and $R^{\mathfrak{G}}_{l}(i,j,k,p)$, both in closed forms, so that
%\begin{align*}
%	 \mathfrak{F}_{l}(i,j,k) &= R^{\mathfrak{F}}_{l+1}(i,j,k)-R^{\mathfrak{F}}_{l}(i,j,k), \\
%	  \mathfrak{G}_{pl}(i,j,k) &= R^{\mathfrak{G}}_{l+1}(i,j,k,p)-R^{\mathfrak{G}}_{l}(i,j,k,p),
%\end{align*}
%hence we obtain
%\begin{align*}
%	  \sum_{l=|j-k|}^{j+k+1}  \mathfrak{F}_{l}(i,j,k) &= R^{\mathfrak{F}}_{j+k+2}(i,j,k)-R^{\mathfrak{F}}_{|j-k|}(i,j,k), \\
%	   \sum_{l=|j-k|}^{j+k+1}  \mathfrak{G}_{pl}(i,j,k) &= R^{\mathfrak{G}}_{j+k+2}(i,j,k,p)- R^{\mathfrak{G}}_{|j-k|}(i,j,k,p),
%\end{align*}
Specifically, the coefficient $1/\mathbf{f}_{jk}$ comes from the first sum whereas the summand $\mathbf{C}_{ijk}$ in comes from the second sum. This completes the proof.
\end{proof}		
\noindent
Next, we use the linearization formula of Lemma \ref{Lemma1auxiliaryForFourierDijk} to establish a computationally efficient formula for $\mathbf{D}_{ijk}$. Notice that, by its definition (Lemma \ref{LemmaFourierCoefficientsWithRespectToGegenbauer}), the Fourier coefficient $\mathbf{D}_{ijk}$ is invariant under permutations of the indices $i$ and $k$. Hence,  no generality is lost in Proposition \ref{ComputationallyEfficientFormulaDijk} by making the simplifying assumption $i \leq k$.

\begin{prop}[Computationally efficient formula for $\mathbf{D}_{ijk}$]\label{ComputationallyEfficientFormulaDijk}
	Let $i,j,k \geq 0$ be integers with $i \leq k$. Then, the Fourier coefficient $\mathbf{D}_{ijk}$ defined in Lemma \ref{LemmaEinsteinMassLessFourier} is given by
	\begin{align*}
		\mathbf{D}_{ijk } &=4 (-1)^{j+1} \omega_{j} K_{i} K_{j}   K_{k} \Bigg[\frac{2}{3} \frac{ \mathbf{a}(i,k)}{ \mathbf{r}(i,k)} +2\sum_{l=0}^{2i} \frac{ l^4}{\mathbf{h}_{l}(i,k)} 
		+\sum_{l=0}^{2i} \sum _{p=1}^j (-1)^p  (p+1)^2 \mathbf{Z}_{pl}(i,j,k)  \Bigg],
	\end{align*}
	where the normalization constant $K_n$ is defined in Lemma \ref{LemmaEigenfunctionsWithRespectToGegenbauer} and the coefficients are given by
	\begin{align*}
		\mathbf{Z}_{pl}(i,j,k)  &=\frac{(2 (k-i+l)+3) \left[ (-1)^l \omega_{i} \omega_{k} 
		(2 (k+l-i)+3)- \mathbf{b}_{l}(i,k)   \right]}{
		\prod _{d=0}^2  (2 (i-k-l+p+d)-3)  \prod _{e=0}^2 (2 (k-i+l+p+e)+3)}
	\end{align*}
	and
	\begin{align*}
	 \mathbf{h}_{l}(i,k) = \prod _{a=0}^4   (2 (k-i+l+a)-1 ) , \quad 
		\mathbf{r}(i,k) =  \prod _{c=0}^3 (2 (i-k+c)-5)  \prod _{d=0}^3 (2 (i+k+d)+1)  
	\end{align*} 
	with
	\begin{align*}
		\mathbf{a}(i,k) &= 512 i^6 k^2+128 i^5 k^3+1536 i^5 k^2+896 i^4 k^3+4512 i^4 k^2-128 i^3 k^5-960 i^3 k^4-1760 i^3 k^3+720 i^3 k^2 \\
		& -576 i^2 k^5-4320 i^2 k^4-12560 i^2 k^3-17672 i^2 k^2-768 i^7 k-576 i^6 k+2624 i^5 k+8432 i^4 k+4968 i^3 k \\
		&-12200 i^2 k+256 i^8-320 i^7-1664 i^6-48 i^5+4784 i^4+3908 i^3-3304 i^2-832 i k^5-6240 i k^4 \\
		&-18216 i k^3-25812 i k^2-17882 i k-4953 i-384 k^5-2880 k^4-8160 k^3-10800 k^2-6576 k-1440 , \\ 
		\mathbf{b}_{l}(i,k) &= 16 i^2 k l+24 i^2 k+8 i^2 l^2+48 i^2 l +52 i^2-16 i k^2 l-24 i k^2-24 i k l^2-72 i k l-52 i k-8 i l^3-48 i l^2 \\
		& -52 i l+8 k^2 l^2-20 k^2+8 k l^3+24 k l^2-20 k l-60 k+2 l^4+12 l^3+8 l^2-30 l-37 .
	\end{align*}
	\end{prop} 		
\begin{proof}
	Let $i,j,k \geq 0$ be integers with $i \leq k$. Then, according to Lemma \ref{LemmaFourierCoefficientsWithRespectToGegenbauer}  and the symmetries of the Gegenbauer polynomials \eqref{AppendixSymmetriesGegenbauer}, we have that
	\begin{align*}
		\mathbf{D}_{ijk } 
 &=  \frac{1}{2 } K_{i}  K_{j} K_{k} \int_{0}^{1}  \frac{1}{y}C_{2 i+1}^{(2)}(y)       \frac{d}{dy} \left(\frac{ (1-y^2)^{3/2}}{y}        C_{2j+1}^{(2)}(y)\right)  
 C_{2 k+1}^{(2)}(y) (1-y^2)^{3/2}  dy \\
 &=  \frac{1}{4 } K_{i}  K_{j} K_{k} \int_{-1}^{1}  \frac{1}{y}C_{2 i+1}^{(2)}(y)       \frac{d}{dy} \left(\frac{ (1-y^2)^{3/2}}{y}        C_{2j+1}^{(2)}(y)\right)  
 C_{2 k+1}^{(2)}(y) (1-y^2)^{3/2}  dy.
	\end{align*}
Moreover, we use linearization formulas \eqref{Asxeto1FourierDijk} and \eqref{Asxeto2FourierDijk} from Lemma \ref{Lemma1auxiliaryForFourierDijk} together with Dougall's formula \eqref{AppendixDougallFormula} to obtain
\begin{align*}
	\mathbf{D}_{ijk } &=\frac{1}{4 } K_{i}  K_{j} K_{k} \int_{-1}^{1}  \frac{1}{y}C_{2 i+1}^{(2)}(y)       \frac{d}{dy} \left(\frac{ (1-y^2)^{3/2}}{y}        C_{2j+1}^{(2)}(y)\right)  
 C_{2 k+1}^{(2)}(y) (1-y^2)^{3/2}  dy \\
 &=\frac{1}{4 } K_{i}  K_{j} K_{k} \sum_{r=0}^{\infty} \mathfrak{M}_{j}(r) \sum_{s=0}^{i} \mathfrak{N}_{i}(s)
 \int_{-1}^{1}  C_{2s}^{(2)}(y)        C_{2r+1}^{(2)}(y)   
 C_{2 k+1}^{(2)}(y) (1-y^2)^{3/2}  dy\\
 &=\frac{1}{4 } K_{i}  K_{j} K_{k} \sum_{r=0}^{\infty} \mathfrak{M}_{j}(r) \sum_{s=0}^{i} \mathfrak{N}_{i}(s)
 D(2r+1,2s,2k+1)\mathds{1} ( |2(k-s)+1|  \leq  2r+1 \leq 2(k+s)+1 ).
\end{align*}
Notice that, since $i \leq k$, we have that $s  \leq k$ and hence $|2(k-s)+1| =2(k-s)+1 $ that leads us to 
\begin{align*}
	\mathbf{D}_{ijk }  
 &=\frac{1}{4 } K_{i}  K_{j} K_{k} \sum_{s=0}^{i} \mathfrak{N}_{i}(s) \sum_{r=k-s}^{k+s} \mathfrak{M}_{j}(r)  D(2r+1,2s,2k+1) \\
 &=\frac{1}{4 } K_{i}  K_{j} K_{k}\sum_{r=k-i}^{k+i}  \mathfrak{M}_{j}(r)  \sum_{s=|k-r|}^{i} \mathfrak{N}_{i}(s)  D(2r+1,2s,2k+1) \\
 &= K_{i}  K_{j} K_{k}\sum_{r=k-i}^{k+i}   \sum_{s=|k-r|}^{i} \mathfrak{J}_{rs}(i,j,k),
\end{align*}
where we set
\begin{align*}
	\mathfrak{J}_{rs}(i,j,k)=\frac{1}{4 } \mathfrak{M}_{j}(r)\mathfrak{N}_{i}(s)  D(2r+1,2s,2k+1).
\end{align*}
Using the closed formulas for $\mathfrak{N}_{i}(s)$, $\mathfrak{M}_{j}(r)$ and $D(2r+1,2s,2k+1)$, we compute 
\begin{align*}
	\mathfrak{J}_{rs}(i,j,k)= \mathfrak{H}_{rs}(i,j,k) + \sum_{p=1}^{j} \mathfrak{K}_{prs}(i,j,k)
\end{align*}
 for some functions $\mathfrak{H}_{rs}(i,j,k) $ and $ \mathfrak{K}_{prs}(i,j,k)$ given by explicit formulas. With this notation, we need to compute
\begin{align*}
	\sum_{s=|k-r|}^{i} \mathfrak{J}_{rs}(i,j,k) = \sum_{s=|k-r|}^{i} \mathfrak{H}_{rs}(i,j,k)
	+\sum_{p=1}^{j}  \left( \sum_{s=|k-r|}^{i}\mathfrak{K}_{prs}(i,j,k)\right).
\end{align*}
To do so, we use Gosper's algorithm to compute both sums
\begin{align*}
	\mathfrak{X}_{r}(i,j,k)=\sum_{s=|k-r|}^{i} \mathfrak{H}_{rs}(i,j,k), \quad 
	\mathfrak{Y}_{pr}(i,j,k)=\sum_{s=|k-r|}^{i}\mathfrak{K}_{prs}(i,j,k).
\end{align*}
We find that
\begin{align*}
	\mathfrak{X}_{r}(i,j,k)&=\mathfrak{X}_{r}^{(1)}(i,j,k)+\mathfrak{X}_{r}^{(2)}(i,j,k)  , \\
	\mathfrak{X}_{r}^{(1)}(i,j,k)&=4(-1)^{i+j+k} (2 i+3) (2 j+3) (2 k+3)\frac{(-1)^r (2 r+3)}{\prod _{a=0}^4 (2 a+2 r-1)} , \\
	\mathfrak{X}_{r}^{(2)}(i,j,k)&=4 (-1)^{j+1} (2 j+3)\frac{\mathfrak{q}_{r}(i,k)}{\prod _{a=0}^4 (2 a+2 r-1)} , \\
	\mathfrak{Y}_{pr}(i,j,k)& = 4(-1)^{j+p+1} (2 j+3) (p+1)^2 (2 r+3)  \frac{(-1)^{i+k+r}(2 i+3) (2 k+3) (2 r+3)-\mathfrak{q}_{r}(i,k) }{ \prod _{a=0}^2 (2 p-2 r+2 a-3)  \prod _{b=0}^2 (2 p+2 r+2 b+3)},
\end{align*} 
where
\begin{align*}
	\mathfrak{q}_{r}(i,k) &= -4 i^2 k^2-12 i^2 k-4 i^2 r^2-12 i^2 r+2 i^4+12 i^3+8 i^2-12 i k^2-36 i k-12 i r^2-36 i r-30 i-4 k^2 r^2 \\
	&-12 k^2 r+2 k^4+12 k^3+8 k^2-12 k r^2-36 k r-30 k+2 r^4+12 r^3+8 r^2-30 r-37.
\end{align*}
Putting them all together, yields
\begin{align*}
	\mathbf{D}_{ijk }  
 % &= K_{i}  K_{j} K_{k}\sum_{r=k-i}^{k+i}   \sum_{s=|k-r|}^{i} \mathfrak{J}_{rs}(i,j,k) \\
 % &= K_{i}  K_{j} K_{k}\sum_{r=k-i}^{k+i} \left(   \sum_{s=|k-r|}^{i} \mathfrak{H}_{rs}(i,j,k)	+\sum_{p=1}^{j}  \left( \sum_{s=|k-r|}^{i}\mathfrak{K}_{prs}(i,j,k)\right) \right) \\
%	&= K_{i}  K_{j} K_{k}\sum_{r=k-i}^{k+i} \left(   \mathfrak{X}_{r}(i,j,k)	+\sum_{p=1}^{j} \mathfrak{Y}_{pr}(i,j,k) \right) \\
%	&= K_{i}  K_{j} K_{k}\sum_{r=k-i}^{k+i} \left(   \mathfrak{X}_{r}^{(1)}(i,j,k)+\mathfrak{X}_{r}^{(2)}(i,j,k)	+\sum_{p=1}^{j} \mathfrak{Y}_{pr}(i,j,k) \right) \\
	&= K_{i}  K_{j} K_{k}\Bigg[  \sum_{r=k-i}^{k+i}  \mathfrak{X}_{r}^{(1)}(i,j,k)+\sum_{r=k-i}^{k+i} \left(\mathfrak{X}_{r}^{(2)}(i,j,k)
	+\sum_{p=1}^{j} \mathfrak{Y}_{pr}(i,j,k) \right) \Bigg].
\end{align*}
An additional application of the Gosper's algorithm shows that
\begin{align*}
	 \sum_{r=k-i}^{k+i}  \mathfrak{X}_{r}^{(1)}(i,j,k)=\frac{4 (2 i+3) (-1)^j (2 j+3) (2 k+3) \mathfrak{p}(i,k) }{ \prod _{c=0}^3 (2 (i-k)+2 c-5) \prod _{d=0}^3 (2 (i+k)+2 d+1)},
\end{align*}
where
\begin{align*}
	 \mathfrak{p}(i,k) &= 
 45 + 184 i + 200 i^2 + 32 i^3 + 16 i^4 + 168 k + 288 i k + 
  288 i^2 k + 200 k^2 + 96 i k^2 + 96 i^2 k^2 + 96 k^3 + 16 k^4,
\end{align*}
and, by changing variables from $r$ to $r=l+k-i$, we deduce
\begin{align*}
		 \mathbf{D}_{ijk}  &=4 (-1)^j\omega_{j} K_i K_j  K_k \Bigg[
		\frac{\omega_{i} \omega_{k} \mathfrak{a}(i,k)}{\mathfrak{r}(i,k)  }
		-\sum _{l=0}^{2 i} \left( 
		\frac{\mathfrak{b}_{l}(i,k)}{\mathfrak{h}_{l}(i,k)}
		+\sum _{p=1}^j   (p+1)^2 (-1)^p 
		\mathfrak{Z}_{pl}(i,j,k) \right)\Bigg],
	\end{align*}
	where the coefficients are given by
	\begin{align*}
\mathfrak{h}_{l}(i,k) = \prod _{a=0}^4 (2 (k-i+l+a)-1) , \quad
\mathfrak{r}(i,k)=  \prod _{b=0}^3 (2 (i-k+b)-5) \prod _{c=0}^3 (2 (i+k+c)+1),    
	\end{align*}  
	and
	\begin{align*}
		\mathfrak{Z}_{pl}(i,j,k)  &=\frac{(2 (k-i+l)+3) \left[ (-1)^l \omega_{i} \omega_{k} 
		(2 (k+l-i)+3)- \mathfrak{b}_{l}(i,k)   \right]}{
		\prod _{d=0}^2  (2 (i-k-l+p+d)-3)  \prod _{e=0}^2 (2 (k-i+l+p+e)+3)}
	\end{align*}
	with
	\begin{align*}
	\mathfrak{a}(i,k) &= 96 i^2 k^2+288 i^2 k+16 i^4+32 i^3+200 i^2+96 i k^2+288 i k+184 i+16 k^4+96 k^3+200 k^2+168 k+45, \\
		\mathfrak{b}_{l}(i,k) &= 16 i^2 k l+24 i^2 k+8 i^2 l^2+48 i^2 l+52 i^2-16 i k^2 l-24 i k^2-24 i k l^2-72 i k l-52 i k-8 i l^3-48 i l^2-52 i l \\
		&+8 k^2 l^2-20 k^2+8 k l^3+24 k l^2-20 k l-60 k+2 l^4+12 l^3+8 l^2-30 l-37.
	\end{align*}
Finally, we note that, since
\begin{align*}
	\mathfrak{b}_{l}(i,k) = \sum_{q=0}^{4} \mathfrak{b}^{(q)}(i,k) l^q,
\end{align*}	
for some coefficients $\mathfrak{b}^{(q)}(i,k)$, Gosper's algorithm computes all sums
\begin{align*}
	\sum _{l=0}^{2 i}  \frac{
		l^q}{\mathfrak{h}_{l}(i,k)}
\end{align*}
with $0 \leq q \leq 3$, that simplifies the expression above and completes the proof.
\end{proof}
\noindent
Finally, we conclude this section by using the linearization formula of Lemma \ref{Lemma1auxiliaryForFourierEijk} to establish a computationally efficient formula for $\mathbf{E}_{ijk}$. Notice that, by its definition (Lemma \ref{LemmaFourierCoefficientsWithRespectToGegenbauer}), the Fourier coefficient $\mathbf{E}_{ijk}$ is invariant under permutations of the indices $i$ and $k$. Hence,  no generality is lost in Proposition \ref{ComputationallyEfficientFormulaDijk} by making the simplifying assumption $i \leq k$.

\begin{prop}[Computationally efficient formula for $\mathbf{E}_{ijk}$]\label{ComputationallyEfficientFormulaEijk}
	Let $i,j,k \geq 0$ be integers with $i \leq k$. Then, the Fourier coefficient $\mathbf{E}_{ijk}$ defined in Lemma \ref{LemmaEinsteinMassLessFourier} is given by
	\begin{align*}
		\mathbf{E}_{ijk} &=  \mathbf{C}_{ijk }  +\mathbf{D}_{ijk }    -16\omega_{j} (j+1) (j+2) K_i K_j K_k \sum _{\lambda =0}^{2 i+1} \frac{\mathbf{u}_{ik}(\lambda)}{\mathbf{v}_{ijk}(\lambda)} ,
	\end{align*} 
	where the normalization constant $K_n$ is defined in Lemma \ref{LemmaEigenfunctionsWithRespectToGegenbauer} and the coefficients are given by
	\begin{align*}
		\mathbf{u}_{ik}(\lambda) &= (\lambda +1) (-2 i+\lambda -2) (2 k+\lambda +4) (-i+k+\lambda +1) (-2 i+2 k+\lambda +1), \\
		\mathbf{v}_{ijk}(\lambda) & = (2 i+2 j-2 k-2 \lambda -1) (2 i+2 j-2 k-2 \lambda +1) (2 i+2 j-2 k-2 \lambda +3) \cdot \\
		&\cdot (-2 i+2 j+2 k+2 \lambda +3) (-2 i+2 j+2 k+2 \lambda +5) (-2 i+2 j+2 k+2 \lambda +7).
	\end{align*}
	\end{prop} 
\begin{proof}
	Let $i,j,k \geq 0$ be integers with $i \leq k$ and recall from Lemma \ref{LemmaFourierCoefficientsWithRespectToGegenbauer} that
	\begin{align*}
		\mathbf{C}_{ijk} &= K_{i}K_{j}K_{k}\int_{0}^{1} \frac{1}{y}
	C_{2 i+1}^{(2)}(y)
	C_{2 j+1}^{(2)}(y)
	C_{2 k+1}^{(2)}(y)  (1-y^2)^{3}  dy  , \\
	\mathbf{D}_{ijk } &=  \frac{1}{2} K_i K_j K_k \int_{0}^{1} \frac{1}{y}  C_{2 i+1}^{(2)}(y)
   \frac{d}{dy} \left[ \frac{  (1-y^2 )^{3/2} }{y}C_{2 j+1}^{(2)}(y) \right] C_{2 k+1}^{(2)}(y) (1-y^2)^{3/2} dy  .
	\end{align*}
	 Then, we expand the derivative using the recurrence relation \eqref{AppendixRecrrenceRelation4} to obtain
\begin{align*}
		 &\mathbf{E}_{ijk }= 
		\frac{1}{2}K_i K_j K_k \int_{0}^{1}  \frac{(1-2 y ^2) }{y}  C_{2 i+1}^{(2)}(y)  \frac{d}{dy} \left[ \frac{  (1-y^2 )^{3/2} }{y}C_{2 j+1}^{(2)}(y) \right]C_{2 k+1}^{(2)}(y)  (1-y^2)^{3/2}  dy \\
		 &= \mathbf{D}_{ijk } - K_i K_j K_k \int_{0}^{1}   y   C_{2 i+1}^{(2)}(y)  \frac{d}{dy} \left[ \frac{  (1-y^2 )^{3/2} }{y}C_{2 j+1}^{(2)}(y) \right]C_{2 k+1}^{(2)}(y)  (1-y^2)^{3/2}  dy \\
		 &= \mathbf{D}_{ijk } - K_i K_j K_k \int_{0}^{1}   y   C_{2 i+1}^{(2)}(y) \left( \frac{1}{y}
		 \frac{d}{dy} \left[    (1-y^2 )^{3/2} C_{2 j+1}^{(2)}(y)
		 \right]     -\frac{1}{y^2}  (1-y^2 )^{3/2} C_{2 j+1}^{(2)}(y)  \right)C_{2 k+1}^{(2)}(y)  (1-y^2)^{3/2}  dy \\
		 % &= \mathbf{D}_{ijk } - K_i K_j K_k \int_{0}^{1}      C_{2 i+1}^{(2)}(y) \left(  -2 (j+1) (j+2)   (1-y^2 )^{1/2}  C_{2( j+1)}^{(1)}(y)   -\frac{1}{y}  (1-y^2 )^{3/2} C_{2 j+1}^{(2)}(y)  \right)C_{2 k+1}^{(2)}(y)  (1-y^2)^{3/2}  dy \\
		 &= \mathbf{D}_{ijk } - K_i K_j K_k \int_{0}^{1}      C_{2 i+1}^{(2)}(y) \left(  -2 (j+1) (j+2)   (1-y^2 )^{2}  C_{2( j+1)}^{(1)}(y)   -\frac{1}{y}  (1-y^2 )^{3} C_{2 j+1}^{(2)}(y)  \right)C_{2 k+1}^{(2)}(y)     dy \\
		 &= \mathbf{D}_{ijk } +\mathbf{C}_{ijk } + 2 K_i K_j   (j+1) (j+2)  K_k \int_{0}^{1}      C_{2 i+1}^{(2)}(y)     C_{2( j+1)}^{(1)}(y)     C_{2 k+1}^{(2)}(y) (1-y^2 )^{2}     dy .
	\end{align*}	
To compute the last integral, we use the symmetries of the Gegenbauer polynomials \eqref{AppendixSymmetriesGegenbauer}, the linearization formula of Lemma \ref{Lemma1auxiliaryForFourierEijk} as well as Dougall's formula \eqref{AppendixDougallFormula} to obtain 
\begin{align*}
	& \int_{0}^{1}      C_{2 i+1}^{(2)}(y)     C_{2( j+1)}^{(1)}(y)     C_{2 k+1}^{(2)}(y) (1-y^2 )^{2}     dy  = 
	 \frac{1}{2}\int_{-1}^{1}      C_{2 i+1}^{(2)}(y)    (1-y^2 )^{1/2}  C_{2( j+1)}^{(1)}(y)     C_{2 k+1}^{(2)}(y) (1-y^2 )^{3/2}     dy= \\
	& \frac{1}{2}\sum_{l=0}^{\infty} \mathfrak{P}_{j}(l)  \int_{-1}^{1}      C_{2 i+1}^{(2)}(y)    C_{2l}^{(2)}(y)    C_{2 k+1}^{(2)}(y) (1-y^2 )^{3/2}     dy = 
	% & \frac{1}{2}\sum_{l=0}^{\infty} \mathfrak{P}_{j}(l)D(2i+1,2l,2k+1) \mathds{1}( 2(k-i) \leq 2l \leq 2(i+k+1)   )   = \\
	 \frac{1}{2}\sum_{l=k-i}^{i+k+1} \mathfrak{P}_{j}(l)D(2i+1,2l,2k+1),
	 % & =\frac{1}{2}\sum_{\lambda=0}^{2i+1} \mathfrak{P}_{j}(\lambda+k-i)D(2i+1,2(\lambda+k-i),2k+1), 
\end{align*}
that completes the proof.
\end{proof}
\noindent

\section{Limiting behaviour for the Fourier coefficients}\label{ProofSectionTable}
\noindent
In this section, we prove Theorem \ref{MainTheorem}. To do so, we shall use the computationally efficient formulas for the Fourier coefficients (Propositions \ref{ComputationallyEfficientFormulaAij}-\ref{ComputationallyEfficientFormulaEijk}) to compute their limits as $i,j,k$ go to infinity.

\subsection{Interactions involving 1-mode}

Firstly, we focus on the Fourier coefficients with only one large index.

\begin{lemma}[Limit of  $\mathbf{C}_{ijk}$ involving 1-mode]\label{LemmaLimitC1mode}
	We have that
	\begin{align*}
	\lim_{\omega_{  \nu  } \rightarrow \infty}	(-1)^  \nu   \omega_{  \nu  }^6 \mathbf{C}_{00  \nu  }=\frac{9216}{\pi ^{3/2}}.
	\end{align*} 
\end{lemma}
\begin{proof}
	 According to Proposition \ref{ComputationallyEfficientFormulaCijk}, we find
	\begin{align*}
		\mathbf{C}_{00  \nu  } &=\frac{18432 (-1)^  \nu   \sqrt{(  \nu  +1) (  \nu  +2)}}{\pi ^{3/2} (2   \nu  -3) (2   \nu  -1) (2   \nu  +1) (2   \nu  +3) (2   \nu  +5) (2   \nu  +7) (2   \nu  +9)},
	\end{align*}
	for all integers $  \nu   \geq 1$, from which the result follows.
%	\begin{align*}
%		\mathbf{C}_{00  \nu  }=\frac{9216}{\pi ^{3/2}} \frac{(-1)^  \nu  }{ \omega_{  \nu  }^6}  + \mathcal{O} \left(\frac{1}{\omega_{  \nu  }^8} \right),
%	\end{align*}	
%	as $\omega_{  \nu  } \rightarrow \infty$.
\end{proof}
\noindent

\begin{lemma}[Limit of  $\mathbf{D}_{ijk}$ involving 1-mode]\label{LemmaLimitD1mode}
	We have that
	\begin{align*}
	\lim_{\omega_{  \nu  } \rightarrow \infty}	(-1)^  \nu   \omega_{  \nu  }^4 \mathbf{D}_{00  \nu  } =\frac{768}{\pi ^{3/2}}, \quad 
	\lim_{\omega_{  \nu  } \rightarrow \infty}\mathbf{D}_{0  \nu  0} =-\frac{16}{\sqrt{\pi }}.
	\end{align*} 
\end{lemma}
\begin{proof}
	 According to Proposition \ref{ComputationallyEfficientFormulaDijk}, we find
	\begin{align*}
		\mathbf{D}_{00  \nu  } &= \frac{1536 (-1)^  \nu   \sqrt{(  \nu  +1) (  \nu  +2)}}{\pi ^{3/2} (2   \nu  -1) (2   \nu  +1) (2   \nu  +3) (2   \nu  +5) (2   \nu  +7)},
	\end{align*}
	for all integers $  \nu   \geq 1$, from which the first of the results above follows. In addition, according to Proposition \ref{ComputationallyEfficientFormulaDijk}, we also find
	\begin{align*}
		 \mathbf{D}_{0  \nu  0}  &= -\frac{16}{\pi ^{3/2}}\frac{2   \nu  +3}{\sqrt{(  \nu  +1) (  \nu  +2)}}\left(\frac{64}{105}+\sum _{p=1}^  \nu   \frac{192 (-1)^p (p+1)^2}{(2 p-3) (2 p-1) (2 p+1) (2 p+3) (2 p+5) (2 p+7)} \right),
	\end{align*}
	for all integers $  \nu   \geq 1$. Now, Gosper's algorithm yields
	\begin{align*}
		\sum _{p=1}^{\infty} \frac{192 (-1)^p (p+1)^2}{(2 p-3) (2 p-1) (2 p+1) (2 p+3) (2 p+5) (2 p+7)}=\frac{\pi }{2}-\frac{64}{105},
	\end{align*}
	from which the second of the results above follows and completes the proof.
\end{proof}
\noindent

\begin{lemma}[Limit of  $\mathbf{E}_{ijk}$ involving 1-mode]\label{LemmaLimitE1mode}
	We have that
	\begin{align*}
	\lim_{\omega_{  \nu  } \rightarrow \infty}(-1)^  \nu  	\omega_{  \nu  }^4 \mathbf{E}_{00  \nu  } = -\frac{768}{\pi ^{3/2}}  , \quad 
	\lim_{\omega_{  \nu  } \rightarrow \infty}\mathbf{E}_{0  \nu  0} = -\frac{16}{\sqrt{\pi }}.
	\end{align*} 
\end{lemma}
\begin{proof}
	According to Proposition \ref{ComputationallyEfficientFormulaEijk} and the closed formulas for $\mathbf{C}_{00  \nu  }$ and $\mathbf{D}_{00  \nu  }$ from Lemmata \ref{LemmaLimitC1mode} and \ref{LemmaLimitD1mode} respectively, we get that
	\begin{align*}
		\mathbf{E}_{00  \nu  }=-\frac{1536 (-1)^  \nu   \sqrt{(  \nu  +1) (  \nu  +2)} (2   \nu  +3)}{\pi ^{3/2} (2   \nu  -3) (2   \nu  -1) (2   \nu  +1) (2   \nu  +5) (2   \nu  +7) (2   \nu  +9)},
	\end{align*}
	for all integers $  \nu   \geq 1$, from which the first of the results above follows. For the second result, we use Proposition \ref{ComputationallyEfficientFormulaEijk} once again to find
	\begin{align*}
		  \mathbf{E}_{0  \nu  0} &=
		  \mathbf{C}_{0  \nu  0}+\mathbf{D}_{0  \nu  0} +\frac{3072 (-1)^  \nu   \sqrt{(  \nu  +1) (  \nu  +2)} (4   \nu   (  \nu  +3)+3)}{\pi ^{3/2} (2   \nu  -3) (2   \nu  -1) (2   \nu  +1) (2   \nu  +3) (2   \nu  +5) (2   \nu  +7) (2   \nu  +9)} ,
	\end{align*} 
		for all integers $  \nu   \geq 1$. Now, Lemmata \ref{LemmaLimitC1mode} and \ref{LemmaLimitD1mode}, together with the fact that $\mathbf{C}_{0  \nu  0}=\mathbf{C}_{00  \nu  }$, yield
	\begin{align*}
		\lim_{\omega_{  \nu  } \rightarrow \infty} \mathbf{E}_{0  \nu  0} &=
		\lim_{\omega_{  \nu  } \rightarrow \infty} \left(\mathbf{C}_{0  \nu  0}+\mathbf{D}_{0  \nu  0}  \right)=\lim_{\omega_{  \nu  } \rightarrow \infty}\mathbf{D}_{0  \nu  0} =-\frac{16}{\sqrt{\pi }},
	\end{align*} 
	that completes the proof.
\end{proof}
\noindent
\subsection{Interactions involving 2-modes}

Secondly, we focus on the Fourier coefficients with only two large indices.  

\begin{lemma}[Limit of  $\mathbf{C}_{ijk}$ involving 2-modes]\label{LemmaLimitC2modes}
	We have that
	\begin{align*}
		\lim_{\omega_{  \nu  } \rightarrow \infty}\mathbf{C}_{0  \nu    \nu  }=\frac{16 \sqrt{2}}{3 \pi ^{3/2}}.
	\end{align*} 
\end{lemma}
\begin{proof}
	 According to Proposition \ref{ComputationallyEfficientFormulaCijk}, we find
	\begin{align*}
		\mathbf{C}_{0  \nu    \nu  } &=\frac{1024 \sqrt{2} (  \nu  +1) (  \nu  +2) (2   \nu  +3)^2}{3 \pi ^{3/2} (4   \nu  +3) (4   \nu  +5) (4   \nu  +7) (4   \nu  +9)},
	\end{align*}
	for all integers $  \nu  \geq 1$,  from which the result follows.
\end{proof}
\noindent

 \begin{lemma}[Limit of  $\mathbf{D}_{ijk}$ involving 2-modes]\label{LemmaLimitD2modes}
	We have that
	\begin{align*}
		\lim_{\omega_{  \nu  } \rightarrow \infty}\frac{\mathbf{D}_{0  \nu    \nu  }}{\omega_{  \nu  }}= -2 \sqrt{\frac{2}{\pi }} , 
		\quad \lim_{\omega_{  \nu  } \rightarrow \infty}\mathbf{D}_{  \nu  0  \nu  }= -\frac{12 \sqrt{2}}{\pi ^{3/2}}. 
	\end{align*} 
\end{lemma}
\begin{proof}
	According to Proposition \ref{ComputationallyEfficientFormulaDijk}, we find
	\begin{align*}
		\frac{\mathbf{D}_{0  \nu    \nu  }}{\omega_{  \nu  }} &=  \frac{512 \sqrt{2}}{\pi ^{3/2}}\left(\frac{(-1)^  \nu   }{(2   \nu  -1) (2   \nu  +1) (2   \nu  +5) (2   \nu  +7)}  +  \mathbf{q}_{  \nu  } \right) ,
	\end{align*}
	for all integers $  \nu  \geq 1$, where we set
	\begin{align*}
		\mathbf{q}_{  \nu  } &= \sum _{p=1}^  \nu   \frac{(2   \nu  +3)(-1)^{p+  \nu  +1} (p+1)^2}{(-2   \nu  +2 p-3) (-2   \nu  +2 p-1) (-2   \nu  +2 p+1) (2   \nu  +2 p+3) (2   \nu  +2 p+5) (2   \nu  +2 p+7)}.
	\end{align*}
	Now, Zeilberger's algorithm yields the 1-step linear recurrence relation
	\begin{align*}
		 \mathbf{q}_{  \nu  +1}-\mathbf{q}_{  \nu  }=\mathbf{r}_{  \nu  },
	\end{align*}
	valid for all integers $  \nu  \geq 1$, where
	\begin{align*}
		\mathbf{r}_{  \nu  } &=-\frac{ (  \nu  +2)^2}{(4   \nu  +5) (4   \nu  +7) (4   \nu  +9) (4   \nu  +11)}+\frac{4 (-1)^  \nu   (  \nu  +2)}{(2   \nu  -1) (2   \nu  +1) (2   \nu  +3) (2   \nu  +5) (2   \nu  +7) (2   \nu  +9)}.
	\end{align*}
	The recurrence relation above is given in a summable form leading to 
	\begin{align*}
		\mathbf{q}_{  \nu  +1}=\mathbf{q}_{1}+\sum_{i=1}^  \nu   \mathbf{r}_{i},
	\end{align*}
	valid for all integers $  \nu   \geq 1$. On the one hand, one can easily compute the initial value
	\begin{align*}
		\mathbf{q}_{1}=-\frac{20}{2079}  .
	\end{align*}
	On the other hand, Gosper's algorithm yields
	\begin{align*}
		\sum_{i=1}^{\infty} \mathbf{r}_{i}=\frac{20}{2079}-\frac{\pi }{256}.
	\end{align*}
	Putting these together, we conclude that
	\begin{align*}
		\lim_{\omega_{  \nu  } \rightarrow \infty}\frac{\mathbf{D}_{0  \nu    \nu  }}{\omega_{  \nu  }}=
		\frac{512 \sqrt{2}}{\pi ^{3/2}}\lim_{\omega_{  \nu  } \rightarrow \infty}\mathbf{q}_{  \nu  }=
		-\frac{512 \sqrt{2}}{\pi ^{3/2}}\frac{\pi }{256}=-2 \sqrt{\frac{2}{\pi }} ,
	\end{align*}
	that completes the first part of the proof. The second part follows similarly.
\end{proof}
\noindent

\begin{lemma}[Limit of  $\mathbf{E}_{ijk}$ involving 2-modes]
	We have that
	\begin{align*}
		\lim_{\omega_{  \nu  } \rightarrow \infty}\frac{\mathbf{E}_{0  \nu    \nu  }}{\omega_{  \nu  }}= -2 \sqrt{\frac{2}{\pi }} , 
		\quad \lim_{\omega_{  \nu  } \rightarrow \infty}\mathbf{E}_{  \nu  0  \nu  }= -\frac{4 \sqrt{2}}{\pi ^{3/2}} . 
	\end{align*} 
\end{lemma}
\begin{proof}
According to Proposition \ref{ComputationallyEfficientFormulaEijk} and Lemmata \ref{LemmaLimitC2modes} and \ref{LemmaLimitD2modes} respectively, we get that
	\begin{align*}
		\frac{\mathbf{E}_{0  \nu    \nu  }}{\omega_{\nu}}&=\frac{\mathbf{C}_{0  \nu    \nu  }}{\omega_{\nu}} +\frac{\mathbf{D}_{0  \nu    \nu  }}{\omega_{\nu}}+ \frac{512 \sqrt{2} (  \nu  +1) (  \nu  +2) (2   \nu  +3)}{3 \pi ^{3/2}   (4   \nu  +3) (4   \nu  +5) (4   \nu  +7) (4   \nu  +9)},
	\end{align*}
	for all integers $  \nu   \geq 1$, and hence we get that
	\begin{align*}
	\lim_{\omega_{  \nu  } \rightarrow \infty}	\frac{\mathbf{E}_{0  \nu    \nu  }}{\omega_{  \nu  }}&=\lim_{\omega_{  \nu  } \rightarrow \infty}  \left(
	 \frac{\mathbf{C}_{0  \nu    \nu  }}{\omega_{  \nu  }}+
	 \frac{\mathbf{D}_{0  \nu    \nu  }}{\omega_{  \nu  }} \right) =
	 \lim_{\omega_{  \nu  } \rightarrow \infty}    
	 \frac{\mathbf{D}_{0  \nu    \nu  }}{\omega_{  \nu  }}  = -2 \sqrt{\frac{2}{\pi }}.
	\end{align*}
	Similarly, Proposition \ref{ComputationallyEfficientFormulaEijk} also yields that 
	\begin{align*}
		\mathbf{E}_{  \nu  0  \nu  }&=\mathbf{C}_{  \nu  0  \nu  } +\mathbf{D}_{  \nu  0  \nu  }-\mathbf{w}_{  \nu  } , \quad \mathbf{w}_{  \nu  }= \sum_{\lambda=0}^{2  \nu  +1}\widetilde{\mathbf{w}}_{  \nu  }(\lambda),
	\end{align*}
	for all integers $  \nu   \geq 1$, where we set
	\begin{align*}
	  \widetilde{\mathbf{w}}_{  \nu  }(\lambda)=\frac{384 \sqrt{2}}{\pi ^{3/2}}\frac{1}{(  \nu  +1) (  \nu  +2)}\frac{(\lambda +1)^3 (2   \nu  -\lambda +2) (2   \nu  +\lambda +4)}{(2 \lambda -3) (2 \lambda -1) (2 \lambda +1) (2 \lambda +3) (2 \lambda +5) (2 \lambda +7)}.
	\end{align*}
	We claim that
	\begin{align*}
		\lim_{\omega_{  \nu  } \rightarrow \infty} \mathbf{w}_  \nu   =-\frac{8 \sqrt{2}}{3 \pi ^{3/2}}.
	\end{align*}
	Indeed, Zeilberger's algorithm yields the 1-step recurrence relation
	\begin{align}\label{RecurrenceRelationForMathbfw}
		(1+  \nu  )\mathbf{w}_{  \nu  }-(3+  \nu  )\mathbf{w}_{  \nu  +1}= \mathbf{m}_{  \nu  },
	\end{align}
	valid for all integers $  \nu   \geq 1$, where
	\begin{align*}
		\mathbf{m}_{  \nu  }=\frac{16 \sqrt{2}}{3 \pi ^{3/2}}+\frac{48 \sqrt{2} \left(512   \nu  ^4+4096   \nu  ^3+11728   \nu  ^2+14144   \nu  +5977\right)}{\pi ^{3/2} (4   \nu  +3) (4   \nu  +5) (4   \nu  +7) (4   \nu  +9) (4   \nu  +11) (4   \nu  +13)}.
	\end{align*}  
	Now, based on \eqref{RecurrenceRelationForMathbfw} and the fact that
	\begin{align*}
		\mathbf{m}_{\infty}:=\lim_{  \nu   \rightarrow \infty}\mathbf{m}_{  \nu  }=\frac{16 \sqrt{2}}{3 \pi ^{3/2}},
	\end{align*}
	 Lemma \ref{AppendixLemmaLimitofa1stOrderLinearRecurrence} yields that
	\begin{align*}
		\mathbf{w}_{\infty}= -\frac{1}{2}  \mathbf{m}_{\infty}=-\frac{8 \sqrt{2}}{3 \pi ^{3/2}}.
	\end{align*}
	Putting these together, we conclude that 
	\begin{align*}
		\lim_{\omega_{  \nu  } \rightarrow \infty} \mathbf{E}_{  \nu  0  \nu  }&=\lim_{\omega_{  \nu  } \rightarrow \infty}\left( \mathbf{C}_{  \nu  0  \nu  } + \mathbf{D}_{  \nu  0  \nu  }-\mathbf{w}_{  \nu  }  \right)= \frac{16 \sqrt{2}}{3 \pi ^{3/2}} -\frac{12 \sqrt{2}}{\pi ^{3/2}}+\frac{8 \sqrt{2}}{3 \pi ^{3/2}}=-\frac{4 \sqrt{2}}{\pi ^{3/2}} ,
	\end{align*}
	where we also used the symmetries $\mathbf{C}_{  \nu  0  \nu  }=\mathbf{C}_{0  \nu    \nu  }$, that completes the proof. 
\end{proof}
\noindent

\subsection{Interactions involving 3-modes}

Finally, we focus on the Fourier coefficients with all their indices being large. 

\begin{lemma}[Limit of  $\mathbf{C}_{ijk}$ involving 3-modes]\label{LemmaLimitC3modes}
	We have that
	\begin{align*}
		\lim_{\omega_{  \nu  } \rightarrow \infty}\mathbf{C}_{  \nu    \nu    \nu  }= \frac{2}{\sqrt{\pi }} .
	\end{align*} 
\end{lemma}
\begin{proof}
	 According to Proposition \ref{ComputationallyEfficientFormulaCijk}, we find
	\begin{align*}
		\mathbf{C}_{  \nu    \nu    \nu  } &=\frac{2048}{\pi ^{3/2}} \sqrt{(  \nu  +1) (  \nu  +2)}(2   \nu  +3)^3 \Bigg(\frac{1}{9 (4   \nu  +3) (4   \nu  +5) (4   \nu  +7) (4   \nu  +9)} \\
		&+\sum_{p=1}^{  \nu  } \frac{(-1)^p (p+1)^2}{(4 p (p+2)+3)^2 (4   \nu  -2 p+3) (4   \nu  -2 p+5) (4   \nu  +2 p+7) (4   \nu  +2 p+9)}\Bigg),
	\end{align*}
for all integers $  \nu   \geq 1$. 	Note that the summand above is a rational expression and a partial fraction decomposition of it yields 
	\begin{align*}
		\mathbf{C}_{  \nu    \nu    \nu  } &=\mathbf{P}_{  \nu  }+\sum _{p=1}^  \nu   \mathbf{Q}_{  \nu  }(p)+\sum _{p=1}^  \nu   \mathbf{R}_{  \nu  }(p),
	\end{align*} 
	where
	\begin{align*}
		\mathbf{P}_{  \nu  } &= \frac{\mathbf{S}_{  \nu  }}{9 (4   \nu  +3) (4   \nu  +5) (4   \nu  +7) (4   \nu  +9)}, \quad
		\mathbf{S}_{  \nu  } = \frac{2048}{\pi ^{3/2}} \sqrt{(  \nu  +1) (  \nu  +2)} (2   \nu  +3)^3 , \\
		\mathbf{Q}_{  \nu  }(p) &=\frac{	\mathbf{S}_{  \nu  }}{1024 (2   \nu  +3)^2 \left(  \nu  ^2+3   \nu  +2\right)}\frac{(-1)^p}{ (2 p+1)^2}  +\frac{	\mathbf{S}_{  \nu  }}{1024 (2   \nu  +3)^2 \left(  \nu  ^2+3   \nu  +2\right) }\frac{(-1)^p}{(2 p+3)^2} \\
		& + \frac{\mathbf{S}_{  \nu  }}{2048 (  \nu  +1)^2 (  \nu  +2)^2} \frac{(-1)^p}{(2 p+1) (2 p+3)}  , \\
		\mathbf{R}_{  \nu  }(p) &= \frac{	\mathbf{S}_{  \nu  }}{2048 (  \nu  +1)^2 (  \nu  +2)^2 (2   \nu  +3)^2} (-1)^p \mathbf{T}_{  \nu  }(p) , \\
		\mathbf{T}_{  \nu  }(p)&=\frac{-4 (8   \nu   (  \nu  +3)+17) p^2-64   \nu   (  \nu  +3) p+64   \nu   (  \nu  +3) (4   \nu   (  \nu  +3)+17)-136 p+1157}{ (4   \nu  -2 p+3) (4   \nu  -2 p+5) (4   \nu  +2 p+7) (4   \nu  +2 p+9)} .
	\end{align*}
	Notice that all fractions above containing $\mathbf{S}_{  \nu  }$ in the numerator have a limit as $  \nu   \rightarrow \infty$. Moreover, 
	\begin{align*}
		\sum _{p=1}^{\infty} \frac{(-1)^p}{(2 p+1)^2} =  -1+G , \quad 
		\sum _{p=1}^{\infty} \frac{(-1)^p}{(2 p+3)^2} = \frac{8}{9}-G  , \quad 
		\sum _{p=1}^{\infty} \frac{(-1)^p}{(2 p+1) (2 p+3)} = \frac{\pi }{4}-\frac{5}{6}  , 
	\end{align*}
	where $G$ stands for the Catalan's constant\footnote{Its numerical value is approximately $G \simeq 0.915966$. } defined by
	\begin{align*}
		G=\sum_{p=0}^{\infty} \frac{(-1)^p}{(2p+1)^2}.
	\end{align*}
	This leads us to
	\begin{align*}
		\lim_{  \nu   \rightarrow \infty} \mathbf{P}_{  \nu  } = \frac{64}{9 \pi ^{3/2}}, \quad 
		\lim_{  \nu   \rightarrow \infty}\sum_{p=1}^  \nu   \mathbf{Q}_{  \nu  }(p)=\frac{2}{\sqrt{\pi }}-\frac{64}{9 \pi ^{3/2}}, \quad 
		\lim_{  \nu   \rightarrow \infty}\left( \mathbf{P}_{  \nu  } +\sum_{p=1}^  \nu   \mathbf{Q}_{  \nu  }(p) \right) = \frac{2}{\sqrt{\pi }}.
	\end{align*}
	Finally, we claim that the remaining term is an error term meaning that
	\begin{align*}
		\lim_{  \nu   \rightarrow \infty}\sum_{p=1}^  \nu   \mathbf{R}_{  \nu  }(p) = 0.
	\end{align*} 
	Indeed, observe that
	\begin{align*}
	 	\frac{	\mathbf{S}_{  \nu  }}{2048 (  \nu  +1)^2 (  \nu  +2)^2 (2   \nu  +3)^2} =\frac{8}{\pi ^{3/2} (3 + 2   \nu  )^2} + \mathcal{O} \left(\frac{1}{  \nu  ^4} \right),
	\end{align*}
	as $  \nu   \rightarrow \infty$,	hence we need to show that
	\begin{align*}
		\lim_{  \nu   \rightarrow \infty}\frac{1}{(3 + 2   \nu  )^2}\sum_{p=1}^  \nu  (-1)^p \mathbf{T}_{  \nu  }(p) =0.
	\end{align*}
	To this end, we define the auxiliary sequence
	\begin{align*}
		\mathbf{t}_{  \nu  } = \frac{1}{(3 + 2   \nu  )^2}\sum_{p=1}^  \nu  (-1)^p \mathbf{T}_{  \nu  }(p),
	\end{align*}
	and Zeilberger's algorithm finds the 1-step linear recurrence relation
	\begin{align*}
		\mathbf{t}_{  \nu  +1}-\mathbf{t}_{  \nu  }=\mathbf{h}_{  \nu  },
	\end{align*}
	valid for all integers $  \nu   \geq 1$, where
	\begin{align*}
		\mathbf{h}_{  \nu  } &= \frac{4 (  \nu  +2) }{(2   \nu  +3)^2 (2   \nu  +5)^2 (4   \nu  +3) (4   \nu  +5) (4   \nu  +7) (4   \nu  +9) (4   \nu  +11) (4   \nu  +13)}\big[4096   \nu  ^6 \\
		&+49152   \nu  ^5+242944   \nu  ^4+632832   \nu  ^3+916144   \nu  ^2+699072   \nu  +219815\big] \\
		&+\frac{(-1)^{  \nu  +1}}{3 (2   \nu  +3)^4 (2   \nu  +5)^4 (6   \nu  +11) (6   \nu  +13)} \big[ 10752   \nu  ^8+172032   \nu  ^7+1198912   \nu  ^6+4753152   \nu  ^5 \\
		&+11722960   \nu  ^4+18413184   \nu  ^3+17979324   \nu  ^2+9972464   \nu  +2403821 \big].
	\end{align*} 
	The recurrence relation above is given in a summable form leading to 
	\begin{align*}
		\mathbf{t}_{  \nu  +1}=\mathbf{t}_{1}+\sum_{i=1}^  \nu   \mathbf{h}_{i},
	\end{align*}
	valid for all integers $  \nu   \geq 1$. On the one hand, one can easily compute the initial value
	\begin{align*}
		\mathbf{t}_{1}= -\frac{1}{25}  \mathbf{T}_{1}(1)=
		-\frac{9017}{170625}.
	\end{align*}
	On the other hand, Gosper's algorithm yields
	\begin{align*}
		\sum_{i=1}^{\infty} \mathbf{h}_{i}=\frac{9017}{170625}.
	\end{align*}
	Putting these together, we conclude that $\lim _{  \nu   \rightarrow \infty}\mathbf{t}_{  \nu  }=0$ that completes the proof.
\end{proof}
\noindent

 \begin{lemma}[Limit of  $\mathbf{D}_{ijk}$ involving 3-modes]\label{LemmaLimitD3modes}
	We have that
	\begin{align*}
		\lim_{\omega_{  \nu  } \rightarrow \infty} \frac{\mathbf{D}_{  \nu    \nu    \nu  }}{\omega_{  \nu  }^2}= -\frac{1}{2 \sqrt{\pi }} .
	\end{align*} 
\end{lemma}
\begin{proof}
	According to Proposition \ref{ComputationallyEfficientFormulaDijk}, we find
	\begin{align*}
		 \frac{\mathbf{D}_{  \nu    \nu    \nu  }}{\omega_{  \nu  }^2} &=-\mathbf{a}_  \nu   - \mathbf{b}_  \nu    -\frac{2   \nu  +3}{\sqrt{(  \nu  +1) (  \nu  +2)}}\mathbf{c}_  \nu    -\frac{2   \nu  +3}{\sqrt{(  \nu  +1) (  \nu  +2)}}\mathbf{d}_  \nu  ,
	\end{align*}
	for all integers $  \nu   \geq 1$, where
	\begin{align*}
		\mathbf{a}_  \nu   &=\frac{32}{\pi ^{3/2} ((  \nu  +1) (  \nu  +2))^{3/2} (2   \nu  +3)} \frac{(2   \nu  +3)^2 (16   \nu   (  \nu  +2) (8   \nu   (  \nu  +2)+11)+45)}{15 (4   \nu  +1) (4   \nu  +3) (4   \nu  +5) (4   \nu  +7)}   , \\
		\mathbf{b}_  \nu   &=\frac{32}{\pi ^{3/2} ((  \nu  +1) (  \nu  +2))^{3/2} (2   \nu  +3)}\sum _{l=0}^{2   \nu  } \widetilde{\mathbf{b}}_\nu   (l)   ,\\
		\mathbf{c}_  \nu   &=\frac{32}{\pi ^{3/2} (  \nu  +1) (  \nu  +2) }\sum _{p=1}^  \nu   \sum _{l=0}^{2   \nu  } \widetilde{ \mathbf{c}} (p,l) , \\
		\mathbf{d}_  \nu   &=\frac{32}{\pi ^{3/2}  (  \nu  +1) (  \nu  +2) (2   \nu  +3)^2 }\sum _{p=1}^  \nu   \sum _{l=0}^{2   \nu  } \widetilde{\mathbf{d}}_{  \nu  }(p,l)  , 
	\end{align*}
	with
	\begin{align*}
		\widetilde{\mathbf{b}}_\nu   (l)&=  \frac{-4   \nu  ^2 (2 l (l+3)+5)-12   \nu   (2 l (l+3)+5)+2 l (l+3) (l (l+3)-5)-37}{(2 l-1) (2 l+1) (2 l+3) (2 l+5) (2 l+7)} , \\
		\widetilde{ \mathbf{c}} (p,l) &=  \frac{(2 l+3)^2 (p+1)^2 (-1)^{l+p}}{(-2 l+2 p-3) (-2 l+2 p-1) (-2 l+2 p+1) (2 l+2 p+3) (2 l+2 p+5) (2 l+2 p+7)}  , \\
		\widetilde{\mathbf{d}}_{  \nu  }(p,l) &=\frac{(2 l+3) (-1)^p (p+1)^2 \left(4   \nu  ^2 (2 l (l+3)+5)+12   \nu   (2 l (l+3)+5)-2 l^2 (l (l+6)+4)+30 l+37\right)}{(-2 l+2 p-3) (-2 l+2 p-1) (-2 l+2 p+1) (2 l+2 p+3) (2 l+2 p+5) (2 l+2 p+7)} .
	\end{align*}
	For $\{\mathbf{a}_{  \nu  }: \nu \geq 0 \}$, one can easily show that
	\begin{align*}
		\lim_{  \nu   \rightarrow \infty} \mathbf{a}_{  \nu  } = 0.
	\end{align*}
	For $\{\mathbf{b}_{  \nu  } : \nu \geq 0 \}$, we note that the numerator of $\widetilde{\mathbf{b}}_{  \nu  }(l)$ is a second order polynomial with respect to $  \nu  $ and, by computing its discriminant, one can obtain  
	\begin{align*}
	  0< -\widetilde{\mathbf{b}}_\nu   (l)\leq \frac{  \nu  ^2}{l^2},
	\end{align*}
	for all integers $1 \leq l \leq 2  \nu  $ and $  \nu   \geq 1$. Hence, we infer
	\begin{align*}
		\lim_{  \nu  \rightarrow \infty} |\mathbf{b}_  \nu  | & \leq \lim_{  \nu  \rightarrow \infty}\frac{32}{\pi ^{3/2} ((  \nu  +1) (  \nu  +2))^{3/2} (2   \nu  +3)}\sum _{l=0}^{2   \nu  } |\widetilde{\mathbf{b}}_\nu   (l) | \\
		& \leq \lim_{  \nu  \rightarrow \infty}\frac{32}{\pi ^{3/2} ((  \nu  +1) (  \nu  +2))^{3/2} (2   \nu  +3)}\left(|\widetilde{\mathbf{b}}_\nu   (0) |+ \sum _{l=1}^{2   \nu  } |\widetilde{\mathbf{b}}_\nu   (l) | \right) \\
		& \leq \lim_{  \nu  \rightarrow \infty}\frac{32}{\pi ^{3/2} ((  \nu  +1) (  \nu  +2))^{3/2} (2   \nu  +3)}\left(\frac{20   \nu   (  \nu  +3)+37}{105} +   \nu  ^2 \sum _{l=1}^{2   \nu  } \frac{1}{l^2} \right) \\ 
		& = \lim_{  \nu  \rightarrow \infty}\frac{32   \nu  ^2}{\pi ^{3/2} ((  \nu  +1) (  \nu  +2))^{3/2} (2   \nu  +3)}\left(\frac{20   \nu   (  \nu  +3)+37}{105   \nu  ^2} +  \sum _{l=1}^{2   \nu  } \frac{1}{l^2} \right) \\ 
		&= 0\cdot  \left(\frac{20}{105}+\frac{\pi ^2}{6} \right)=0,
	\end{align*}
	and therefore
	\begin{align*}
		\lim_{  \nu   \rightarrow \infty} \mathbf{b}_{  \nu  } = 0.
	\end{align*}
	For $\{\mathbf{c}_{  \nu  }: \nu \geq0 \}$, we write
	\begin{align*}
	\mathbf{c}_  \nu   &=\sum _{p=1}^  \nu   \sum _{l=0}^{2   \nu  } \overline{\mathbf{c}}_{  \nu  }(p,l), \quad \overline{\mathbf{c}}_{  \nu  }(p,l)= \frac{32}{\pi ^{3/2} (  \nu  +1) (  \nu  +2) } \widetilde{ \mathbf{c}} (p,l) 
	\end{align*}
	and  Zeilberger's algorithm yields the 1-step recurrence relation for the summand
	\begin{align*}
		(  \nu  +1) \overline{\mathbf{c}}_{  \nu  }(p,l)-(3+  \nu  ) \overline{\mathbf{c}}_{  \nu  +1}(p,l)=0,
	\end{align*}
	valid for all integers $  \nu   \geq 1$, $0 \leq l \leq 2  \nu  $ and $1\leq p \leq   \nu  $. Summing over $p$ and $l$ as well as by adding and subtracting suitable terms, we obtain
	\begin{align*}
		(  \nu  +1) \sum_{p=1}^{  \nu  } \overline{\mathbf{c}}_{  \nu  }(p,l) &=(3+  \nu  )\sum_{p=1}^{  \nu  } \overline{\mathbf{c}}_{  \nu  +1}(p,l)  \Longrightarrow  \\
		(  \nu  +1) \sum_{p=1}^{  \nu  } \overline{\mathbf{c}}_{  \nu  }(p,l) &=
		-(3+  \nu  )  \overline{\mathbf{c}}_{  \nu  +1}(  \nu  +1,l) 
		+(3+  \nu  )\sum_{p=1}^{  \nu  +1} \overline{\mathbf{c}}_{  \nu  +1}(p,l)  \Longrightarrow \\
		(  \nu  +1) \sum_{p=1}^{  \nu  }\sum_{l=0}^{2  \nu  } \overline{\mathbf{c}}_{  \nu  }(p,l) &=
		-(3+  \nu  ) \sum_{l=0}^{2  \nu  } \overline{\mathbf{c}}_{  \nu  +1}(  \nu  +1,l) 
		+(3+  \nu  )\sum_{p=1}^{  \nu  +1}\sum_{l=0}^{2  \nu  } \overline{\mathbf{c}}_{  \nu  +1}(p,l)  \Longrightarrow \\
		(  \nu  +1) \sum_{p=1}^{  \nu  }\sum_{l=0}^{2  \nu  } \overline{\mathbf{c}}_{  \nu  }(p,l) &=
		-(3+  \nu  ) \sum_{l=0}^{2  \nu  } \overline{\mathbf{c}}_{  \nu  +1}(  \nu  +1,l) 
		-(3+  \nu  )\sum_{p=1}^{  \nu  +1} \overline{\mathbf{c}}_{  \nu  +1}(p,2  \nu  +1)\\
		&-(3+  \nu  )\sum_{p=1}^{  \nu  +1} \overline{\mathbf{c}}_{  \nu  +1}(p,2  \nu  +2)  
		+(3+  \nu  )\sum_{p=1}^{  \nu  +1}\sum_{l=0}^{2  \nu  +2} \overline{\mathbf{c}}_{  \nu  +1}(p,l)	 ,
	\end{align*}
	and hence we arrive at the 1-step recurrence relation for the sum, 
	\begin{align}\label{RecurrenceRelationForMathbfc}
		(  \nu  +1) \mathbf{c}_{  \nu  } -(3+  \nu  )\mathbf{c}_{  \nu  +1}&=\mathbf{g}_{  \nu  }	,
	\end{align}
	valid for all integers $  \nu   \geq 1$, where we set
	\begin{align*}
		\mathbf{g}_{  \nu  } &= -(3+  \nu  ) \sum_{l=0}^{2  \nu  } \overline{\mathbf{c}}_{  \nu  +1}(  \nu  +1,l) 
		-(3+  \nu  )\sum_{p=1}^{  \nu  +1} \overline{\mathbf{c}}_{  \nu  +1}(p,2  \nu  +1) -(3+  \nu  )\sum_{p=1}^{  \nu  +1} \overline{\mathbf{c}}_{  \nu  +1}(p,2  \nu  +2).
	\end{align*}
	In order to apply Zeilberger's algorithm, we must first simplify $\mathbf{g}_{  \nu  }$ so that it can be written as a sum of one variable varying over the same range. To do so, we change variables to the first sum in $\mathbf{g}_{  \nu  }$ and obtain
	\begin{align*}
		\sum_{l=0}^{2  \nu  } \overline{\mathbf{c}}_{  \nu  +1}(  \nu  +1,l) &=
		\sum_{ \substack   { l=0 \\ l \text{ odd} }}^{2  \nu  } \overline{\mathbf{c}}_{  \nu  +1}(  \nu  +1,l) +\sum_{ \substack  { l=0 \\ l \text{ even} }}^{2  \nu  } \overline{\mathbf{c}}_{  \nu  +1}(  \nu  +1,l) \\
		&=
		\sum_{ \substack  { q=0   }}^{  \nu  -1} \overline{\mathbf{c}}_{  \nu  +1}(  \nu  +1,2q+1) +\sum_{ \substack  { q=0  }}^{  \nu  } \overline{\mathbf{c}}_{  \nu  +1}(  \nu  +1,2q) \\
		&=
		\sum_{ \substack  { p=1   }}^{  \nu  } \overline{\mathbf{c}}_{  \nu  +1}(  \nu  +1,2p-1) +\sum_{ \substack  { p=1  }}^{  \nu  +1} \overline{\mathbf{c}}_{  \nu  +1}(  \nu  +1,2p-2) \\
		&=-\overline{\mathbf{c}}_{  \nu  +1}(  \nu  +1,2  \nu  +1)+
		\sum_{ \substack  { p=1   }}^{  \nu  +1}\left(
		\overline{\mathbf{c}}_{  \nu  +1}(  \nu  +1,2p-1)+ 
		\overline{\mathbf{c}}_{  \nu  +1}(  \nu  +1,2p-2) \right)
	\end{align*}
	and putting all together yields
	\begin{align*}
		\mathbf{g}_{  \nu  } &=\mathbf{f}_{  \nu  }+\mathbf{r}_{  \nu  } +\mathbf{s}_{  \nu  } ,
	\end{align*}
	where
	\begin{align*}
		\mathbf{f}_{  \nu  } =\frac{32 (-1)^{  \nu  +1} (  \nu  +2) (4   \nu  +5)^2}{3 \pi ^{3/2} (2   \nu  -1) (2   \nu  +1) (2   \nu  +3)^2 (6   \nu  +7) (6   \nu  +11)},  \quad 
		\mathbf{r}_{  \nu  } = \sum_{p=1}^{  \nu  +1} \mathbf{u}_{  \nu  }(p), \quad 
		\mathbf{s}_{  \nu  } = \sum_{p=1}^{  \nu  +1} \mathbf{v}_{  \nu  }(p),
	\end{align*}
	with
	\begin{align*}
		 \mathbf{u}_{  \nu  }(p) &= \frac{512}{\pi ^{3/2}}\frac{2   \nu  +3}{  \nu  +2} \frac{(-1)^{p+1} (p+1)^2 (16   \nu   (  \nu  +3)+2 p (p+2)+35)}{ \prod _{a=0}^3 (2 a+4   \nu  -2 p+1) \prod _{b=0}^3 (2 b+4   \nu  +2 p+5)}, \\
		  \mathbf{v}_{  \nu  }(p) &=\frac{1024}{\pi ^{3/2}} (  \nu  +2)\frac{(-1)^{  \nu  +1} p \left(2   \nu   (  \nu  +4)+16 p^2+5\right)}{ \prod _{a=0}^3 (2 a+2   \nu  -4 p+1)  \prod _{b=0}^3 (2 b+2   \nu  +4 p+1)}.
	\end{align*}
	Now, Zeilberger's algorithm applies. Firstly, we obtain the 1-step recurrence relation for $\mathbf{s}_{  \nu  }$, that is
	\begin{align*}
		\mathbf{s}_{  \nu  +1}-\mathbf{s}_{  \nu  } &=\mathbf{z}_  \nu  , 
	\end{align*}
	valid for all integers $  \nu   \geq 1$, where
	\begin{align*}
		\mathbf{z}_  \nu   &= \frac{16 (-1)^{  \nu  +1} \left(384   \nu  ^5+4160   \nu  ^4+17832   \nu  ^3+37852   \nu  ^2+39830   \nu  +16633\right)}{3 \pi ^{3/2} (2   \nu  +1) (2   \nu  +3)^2 (2   \nu  +5)^2 (6   \nu  +11) (6   \nu  +13) (6   \nu  +17)} \\
		&+\frac{48 (-1)^{  \nu  +1}}{\pi ^{3/2} (2   \nu  +1) (2   \nu  +3) (2   \nu  +5) (2   \nu  +7) (2   \nu  +9)}.
	\end{align*}
	Since the latter is given in a summable form, we infer
	\begin{align*}
	 	\mathbf{s}_{  \nu  +1} =\mathbf{s}_{1}+ \sum_{\lambda=1}^{  \nu  }\mathbf{z}_{\lambda},
	\end{align*}
    for all integers $  \nu   \geq 1$. On the one hand, one computes
    \begin{align*}
    	\mathbf{s}_{1}=-\frac{6095872}{3828825 \pi ^{3/2}}.
    \end{align*}
    On the other hand, Gosper's algorithm yields
	\begin{align*}
		\sum_{\lambda=1}^{\infty}\mathbf{z}_{\lambda}=-\frac{1}{2 \sqrt{\pi }}+\frac{6095872}{3828825 \pi ^{3/2}}.
	\end{align*}
	Hence, we obtain
	\begin{align*}
		\lim_{  \nu   \rightarrow \infty}\mathbf{s}_{  \nu  } = -\frac{1}{2 \sqrt{\pi }}.
	\end{align*}
	By a similar argument, one infers that $\lim_{  \nu   \rightarrow \infty}\mathbf{r}_{  \nu  } =0$ 	and, since $\lim_{  \nu   \rightarrow \infty}\mathbf{f}_{  \nu  } = 0$,  	 we deduce that
	\begin{align}\label{LimitofMathbfgk}
		\mathbf{g}_{\infty}:=\lim_{  \nu   \rightarrow \infty}\mathbf{g}_{  \nu  } =-\frac{1}{2 \sqrt{\pi }}.
	\end{align}
	Now, based on \eqref{RecurrenceRelationForMathbfc} and \eqref{LimitofMathbfgk}, Lemma \ref{AppendixLemmaLimitofa1stOrderLinearRecurrence} yields that
	\begin{align*}
		-2\mathbf{c}_{\infty}=   \mathbf{g}_{\infty}.
	\end{align*}
A similar approach to the one obtained above for $\{\mathbf{c}_{  \nu  } :  \nu \geq0 \}$ can be applied to $\{\mathbf{d}_{  \nu  } :  \nu \geq0 \}$ and yields 
\begin{align*}
	\lim_{  \nu   \rightarrow \infty}\mathbf{d}_{  \nu  }=0.
\end{align*}
 Finally, by putting all together, we obtain
	\begin{align*}
		\lim_{\omega_{  \nu  }\rightarrow \infty  } \frac{\mathbf{D}_{  \nu    \nu    \nu  }}{\omega_{  \nu  }^2} &=  \lim_{\omega_{  \nu  }\rightarrow \infty  } \left(-  \frac{2   \nu  +3}{\sqrt{(  \nu  +1) (  \nu  +2)}}\mathbf{c}_  \nu   \right)=-2\mathbf{c}_{\infty} =  \mathbf{g}_{\infty}=-\frac{1}{2 \sqrt{\pi }},
	\end{align*}
	that completes the proof.
\end{proof}
\noindent

\begin{lemma}[Limit of  $\mathbf{E}_{ijk}$ involving 3-modes]
	We have that
	\begin{align*}
		\lim_{\omega_{  \nu  } \rightarrow \infty} \frac{\mathbf{E}_{  \nu    \nu    \nu  }}{\omega_{  \nu  }^2}= -\frac{1}{2 \sqrt{\pi }} .
	\end{align*}
\end{lemma}
\begin{proof}
	According to Proposition \ref{ComputationallyEfficientFormulaEijk}, we get that
	\begin{align*}
		\frac{\mathbf{E}_{  \nu    \nu    \nu  }}{\omega_{  \nu  }^2} &=
		\frac{\mathbf{C}_{  \nu    \nu    \nu  }}{\omega_{  \nu  }^2}+
		\frac{\mathbf{D}_{  \nu    \nu    \nu  }}{\omega_{  \nu  }^2}+
		 \frac{128}{\pi ^{3/2}}\frac{(-1)^  \nu  }{(2   \nu  +3)^2}  \frac{2   \nu  +3}{\sqrt{(  \nu  +1) (  \nu  +2)}}\mathbf{F}_{  \nu  } ,
	\end{align*}
	for all integers $  \nu   \geq 1$, where
	\begin{align*}
		 \mathbf{F}_{  \nu  }= \sum_{\lambda=0}^{2  \nu  +1}\widetilde{\mathbf{F}}_{  \nu  }(\lambda),
	\end{align*}
	with
	\begin{align*}
		\widetilde{\mathbf{F}}_{  \nu  }(\lambda)=\frac{(\lambda +1)^3 (2   \nu  -\lambda +2) (2   \nu  +\lambda +4)}{(2   \nu  -2 \lambda -1) (2   \nu  -2 \lambda +1) (2   \nu  -2 \lambda +3) (2   \nu  +2 \lambda +3) (2   \nu  +2 \lambda +5) (2   \nu  +2 \lambda +7)}.
	\end{align*} 
	Now, Zeilberger's algorithm yields the 1-step recurrence relation
	\begin{align*}
		\mathbf{F}_{  \nu  +1}-\mathbf{F}_{  \nu  }=\mathbf{G}_{  \nu  },
	\end{align*}
	valid for all integers $  \nu  \geq 1$, where
	\begin{align*}
		\mathbf{G}_{  \nu  }=\frac{(  \nu  +2) \left(432   \nu  ^4+3456   \nu  ^3+10288   \nu  ^2+13504   \nu  +6601\right)}{48  (2   \nu  +3)^2 (2   \nu  +5)^2 (6   \nu  +11) (6   \nu  +13) }-\frac{(  \nu  +2) \left(4   \nu  ^2+16   \nu  +13\right)}{16 (2   \nu  +3)^2 (2   \nu  +5)^2 }.
	\end{align*} 
	Since the latter is given in a summable form, we obtain
	\begin{align*}
		\mathbf{F}_{  \nu  +1}=\mathbf{F}_{1}+ \sum_{j=1}^{  \nu  }\mathbf{G}_{j},
	\end{align*} 
	for all integers $  \nu  \geq 1$. On the one hand, one can easily compute
	\begin{align*}
		\mathbf{F}_{1}=\frac{4112}{225225}.
	\end{align*}
	One the other hand, Gosper's algorithm yields
	\begin{align*}
		\sum_{j=1}^{\infty}\mathbf{G}_{j} = \frac{\log (3)}{128}-\frac{56509}{7207200}.
	\end{align*}
	Consequently, we deduce that
	\begin{align*}
		\lim_{  \nu   \rightarrow \infty} \mathbf{F}_{  \nu  }=\frac{4112}{225225}+\frac{\log (3)}{128}-\frac{56509}{7207200}=\frac{1}{384} (4+\log (27)).
	\end{align*}
	In particular, the latter shows that
	\begin{align*}
		\lim_{  \nu   \rightarrow \infty}  \frac{(-1)^  \nu  }{(2   \nu  +3)^2}  \frac{2   \nu  +3}{\sqrt{(  \nu  +1) (  \nu  +2)}}\mathbf{F}_{  \nu  }  =0.
	\end{align*}
	Putting all together, yields
	\begin{align*}
	\lim_{\omega_{  \nu  } \rightarrow \infty}	\frac{\mathbf{E}_{  \nu    \nu    \nu  }}{\omega_{  \nu  }^2} &=
		\lim_{\omega_{  \nu  } \rightarrow \infty}	 \left(\frac{\mathbf{C}_{  \nu    \nu    \nu  }}{\omega_{  \nu  }^2}+
		\frac{\mathbf{D}_{  \nu    \nu    \nu  }}{\omega_{  \nu  }^2} \right)=-\frac{1}{2 \sqrt{\pi }},
	\end{align*}
	according to Lemmata \ref{LemmaLimitC3modes} and \ref{LemmaLimitD3modes}, that completes the proof.
\end{proof}
\noindent

\appendix

\section{List of formulas and identities}
\noindent
In this section, we list all formulas and identities used above.

\begin{itemize}[leftmargin=*]
%%%%%%%%%%%%
\item The symmetries of Jacobi polynomials \cite{MR0372517}:
\begin{align}
		 P_{n}^{ \left(\alpha,\beta \right)}(y) &=(-1)^n  P_{n}^{ \left(\beta,\alpha \right)}(-y),\label{AppendixSymmetriesJacobi} \\
		 C_{n}^{(\alpha)} (y) &=(-1)^n C_{n}^{(\alpha)} (-y),\label{AppendixSymmetriesGegenbauer}
	\end{align}
	valid for all integers $n \geq 0$, real $y\in[-1,1]$ and real $\alpha,\beta >-1$. 
%%%%%%%%%%%%
\item The relation between Jacobi and Gegenbauer polynomials (18.7.16 in \cite{MR2723248}),
	\begin{align}\label{AppendixRelationJacobitoGegenbauer}
	P_n^{\left(3/2,1/2\right)}(-\cos (2 x))= \frac{ \Gamma  \left(n+\frac{3}{2}\right) }{ \sqrt{\pi }  \Gamma  (n+3)   }\frac{1}{\sin(x)} C_{2 n+1}^{(2)}(\sin (x)),
\end{align}
valid for all integers $n \geq 0$ and real $x\in[0,\pi/2]$.
%%%%%%%%%%%%
\item The orthogonality of the Gegenbauer polynomials \cite{MR0372517}, 
\begin{align}
\int_{-1}^1C_m^{(1)}(y) C_n^{(1)}(y)  (1-y^2)^{1/2}  dy &= \frac{\pi }{2}  \mathds{1}(n=m), \label{AppendixOrthogonalityGegenbauerExponent1} \\
	  \int_{-1}^{1} C_{n}^{(2)} (y)C_{m}^{(2)} (y) (1-y^2)^{3/2} dy &= \xi _{n} \mathds{1}(n=m), \label{AppendixOrthogonalityGegenbauerExponent2} \\
	  \int_{-1}^1  C_m^{\left( 5/2  \right)}(y) C_n^{\left( 5/2  \right)}(y)  dy \left(1-y^2\right)^2&= \rho_{n} \mathds{1}(n=m),\label{AppendixOrthogonalityGegenbauerExponent5half} \\
	   \int_{-1}^1 C_n^{(3)}(y) C_m^{(3)}(y)   dy  (1-y^2)^{5/2} &=  \zeta_{n}\mathds{1}(n=l)\label{AppendixOrthogonalityGegenbauerExponent3} \\
	   \int_{-1}^{1}  C_{n}^{\left( 7/2  \right)}(y) C_{ l}^{\left( 7/2  \right)}(y) (1-y^2)^3 dy &=\theta_{ n} \mathds{1} (n=l),\label{AppendixOrthogonalityGegenbauerExponent7half} 
\end{align}
valid for all integers $n,m \geq 0$, where the coefficient is given by
\begin{align*}
	\xi_{n}   &=\frac{\pi }{8}  (n+1) (n+3), \\
	\rho_{n} &= \frac{ 2 (n+4) (n+3) (n+2) (n+1)  }{9 (2 n+5)}, \\
	\zeta_{n}  & = \frac{ \pi }{128} (n+1) (n+2) (n+4) (n+5), \\
	\theta _{n} & =  \frac{2 (n +6) (n +5) (n +4) (n +3) (n +2) (n +1)}{225 (2 n +7)}.
\end{align*}
\item The recurrence relations (Section 18.9 in \cite{MR2723248}),
 \begin{align}
	C_n^{(2)}(y) &=\frac{2}{n+2} \left( C_n^{(3)}(y)-C_{n-2}^{(3)}(y) \right),
	\label{AppendixRecrrenceRelation1} \\
	\left(1-y^2\right) C_{2 n+1}^{(2)}(y) &=  \frac{n+2}{2}  C_{2 n+1}^{(1)}(y)-
\frac{n+1}{2}  C_{2 n+3}^{(1)}(y) , \label{AppendixRecrrenceRelation2}\\
y C_{2 n+2}^{(1)}(y) &= \frac{1}{2} (C_{2 n+1}^{(1)}(y)+C_{2 n+3}^{(1)}(y)),\label{AppendixRecrrenceRelation3}  \\
\frac{d}{dy} \left(   (1-y^2 )^{3/2} C_{2 n+1}^{(2)}(y) \right) &= -2 (n+1) (n+2)   (1-y^2 )^{1/2}  C_{2( n+1)}^{(1)}(y),\label{AppendixRecrrenceRelation4} \\
C_{2 n+1}^{(2)}(y) &= \frac{4 n+4 }{2 n+1}y C_{2 n}^{(2)}(y)
	-\frac{2 n+3 }{2 n+1} C_{2 n-1}^{(2)}(y),\label{AppendixRecrrenceRelation5},
	% y C_{2 n}^{(2)}(y) &=  \frac{2 n+3 }{4 n+4}C_{2 n-1}^{(2)}(y)+\frac{2 n+1 }{4 n+4} C_{2 n+1}^{(2)}(y),\label{AppendixRecrrenceRelation6} 
\end{align}
valid for all integers $n \geq 0$ and real $y \in [-1,1]$. 
%%%%%%%%%%%%
\item The computation of integrals in terms of hypergeometric functions (relation 4, section 7.311 in \cite{MR3307944}) 
\begin{align}
	 \int_{-1}^{1}(1-y^2)^{3}	C_{2 m}^{(2)}(y) dy &=\sum _{l=0}^{2 m} \frac{96 \sqrt{\pi } (-1)^l (l+3) (l+2) (l+1)  \Gamma  (l+2 m+4)}{ \Gamma \left(l+\frac{5}{2}\right)  \Gamma  (l+8)  \Gamma  (-l+2 m+1)}, \label{AppendixComputationIntegralHyperGeometricFunction1} \\
	 \int_{-1}^{1}  (1-y^2)^{2} 
	  C_{2( l+1)}^{(1)}(y) dy &=\sum _{\mu=0}^{2(l+1)} \frac{(-1)^{\mu } 4^{\mu +3}  \Gamma  (2 l+\mu +4)}{(\mu +3) (\mu +4) (\mu +5)  \Gamma (2 \mu +2) \Gamma  (2 l-\mu +3)}  ,\label{AppendixComputationIntegralHyperGeometricFunction2}
\end{align}
valid for all integers $m,l \geq 0$.
%%%%%%%%%%%%
\item Dougall's formula (Proposition 3.2 in \cite{MR1789139}, relation 18.18.22 in \cite{MR2723248})
	\begin{align}\label{AppendixDougallFormula}
	\int_{-1}^{1} C_{l}^{(2)}(y)  C_{n}^{(2)}(y)C_{m}^{(2)}(y)(1-y^2)^{3/2} dy= D(l,n,m) 
	\mathds{1} \begin{pmatrix}
		|n-m| \leq l \leq n+m, \\
		n+m+l \in 2\mathbb{N}\cup \{0\} 
	\end{pmatrix} ,
	\end{align}
	valid for all integers $n,m,l \geq 0$, where the coefficients are given by
	 \begin{align*}
	 	D(l,n,m)  =\frac{\pi}{128}   (m+n-l+2) (m-n+l+2) (-m+n+l+2) (l+m+n+6).
	 \end{align*}
	 The relation \eqref{AppendixDougallFormula} is an equivalent formulation for the linearization formula (relation 18.18.22 in \cite{MR2723248})
	 \begin{align}\label{AppendixLinearizatonFormula}
	 	C_{n}^{(2)}(y)C_{m}^{(2)}(y) = \sum_{l=0}^{\min (n,m)} L_{k}(n,m) C_{m+n-2l}^{(2)}(y),
	 \end{align}
	 valid for all integers $n,m \geq 0$, where the coefficients are given by
	 \begin{align*}
	 	L_{k}(n,m)  =   -\frac{(l+1) (l-m-1) (l-n-1) (l-m-n-3)}{(-2 l+m+n+1) (-2 l+m+n+3)}.
	 \end{align*} 
%%%%%%%%%%%%
\item  The connection formulas for Gegenbauer polynomials (relation 18.18.16	 in \cite{MR2723248}),
\begin{align}
	C_{2 n}^{(2)}(y) &= \sum _{\mu_1=0}^n \eta_{n}^{(1)} (\mu_1 )C_{2 \mu_1}^{\left( 5/2  \right)}(y),\label{AppendixConnectionFormulaExponent5half} \\
	 C_{2 (l+1)}^{(1)}(y) &= \sum _{\mu_2 =0}^{l+1}\eta_{l}^{(2)} (\mu_2 ) C_{2 \mu_2 }^{\left( 5/2  \right)}(y) ,\label{AppendixConnectionFormulaExponent5halfAgain}\\
	 C_{2 n}^{(2)}(y)&=   \sum _{l=0}^n \eta^{(3)}_{n}(l) C_{2 l}^{\left( 7/2  \right)}(y),\label{AppendixConnectionFormulaExponent7half},
\end{align}
valid for all integers $n \geq 0$ and real $y \in [-1,1]$, where the coefficients are given by
\begin{align*} 
	\eta_{n}^{(1)} (\mu_1 ) &= -\frac{  3 (4 \mu_1+5)  \Gamma  \left(-\mu_1+n-\frac{1}{2}\right) \Gamma  (\mu_1+n+2) }{16\Gamma  (-\mu_1+n+1) \Gamma  \left(\mu_1+n+\frac{7}{2}\right)}, \\
	\eta_{l}^{(2)} (\mu_2 )  &= \frac{ 9 (4 \mu_2 +5) \Gamma  \left(l-\mu_2 -\frac{1}{2}\right) \Gamma (l+\mu_2 +2) }{32 \Gamma (l-\mu_2 +2) \Gamma \left(l+\mu_2 +\frac{9}{2}\right)}, \\
	\eta^{(3)}_{n}(l) &=\frac{45 (4 l+7)  \Gamma  \left(-l+n-\frac{3}{2}\right)  \Gamma  (l+n+2)}{64  \Gamma (-l+n+1)  \Gamma \left(l+n+\frac{9}{2}\right)}  .
\end{align*} 
\end{itemize}

 \section{First order linear recurrence relations}
\noindent
\begin{lemma}[Solution to first order linear recurrence relation]\label{AppendixLemma1stOrderLinearRecurrence}
	Let $\{F_{n}: n \geq 0\}$ and $\{G_{n}: n \geq 0\}$ be two sequences with $F_n \neq 0$, for all integers $n\geq 0$. Assume that $\{A_{n}: n \geq 0\}$ is a sequence that satisfies the first order linear recurrence relation
	\begin{align*}
	A_{n}=F_{n} A_{n-1}+G_{n},
\end{align*}
for all integers $n \geq 0$. Then, 
\begin{align*}
	A_{n} =\prod_{j=0}^{n} F_j \left(  \frac{A_{0}}{ F_0 }  +  \sum_{i=1}^{n} \frac{G_{i}}{\prod_{j=0}^{i}F_j} \right)  .
\end{align*}  
\end{lemma}
\begin{proof}
	Let $n \geq 0$ be any integer. Then, we set
	\begin{align*}
		B_n = \frac{A_{n}}{\prod_{j=0}^{n} F_j } 
	\end{align*}
	and obtain
	\begin{align*}
		A_{n}=F_{n} A_{n-1}+G_{n} &\Longrightarrow
		\frac{A_{n}}{\prod_{j=0}^{n} F_j } =\frac{F_{n}}{\prod_{j=0}^{n}F_j} A_{n-1}+ \frac{G_{n}}{\prod_{j=0}^{n}F_j}  
		\Longrightarrow
		\frac{A_{n}}{\prod_{j=0}^{n} F_j } =\frac{A_{n-1}}{\prod_{j=0}^{n-1}F_j}  + \frac{G_{n}}{\prod_{j=0}^{n}F_j} \\
		&\Longrightarrow
		B_{n} =B_{n-1}  + \frac{G_{n}}{\prod_{j=0}^{n}F_j} 
		\Longrightarrow
	B_{n}-B_0 = \sum_{i=1}^{n}	\left( B_{i} -B_{i-1} \right)=  \sum_{i=1}^{n} \frac{G_{i}}{\prod_{j=0}^{i}F_j} \\
	&\Longrightarrow
	A_n= \prod_{j=0}^{n} F_j B_{n} =\prod_{j=0}^{n} F_j \left( B_0  +  \sum_{i=1}^{n} \frac{G_{i}}{\prod_{j=0}^{i}F_j} \right),
	\end{align*}
	that completes the proof.
\end{proof}
\noindent

\begin{lemma}[Limit of the solution to a first order linear recurrence relation]\label{AppendixLemmaLimitofa1stOrderLinearRecurrence}
	Let $\{F_{n}: n \geq 0\}$ be a sequence such that the limit $ \lim_{n \rightarrow \infty} F_n \in \mathbb{R}$ exists. Assume that $\{A_{n}: n \geq 0\}$ is a sequence that satisfies the first order linear recurrence relation
	\begin{align}\label{RecurreneceLimitApendix}
	(1+  n  )A_{ n  }-(3+  n  )A_{  n  +1}= F_{  n  },
\end{align}
for all integers $n \geq 0$. Then, the limit $ \lim_{n \rightarrow \infty} A_n$ exists and
	\begin{align*}
		 \lim_{n \rightarrow \infty} A_n =  -\frac{1}{2} \lim_{n \rightarrow \infty} F_n.
	\end{align*} 
\end{lemma}
\begin{proof}
Let $\{F_{n}: n \geq 0\}$ be a sequence such that the limit $F_{\infty}:=\lim_{n \rightarrow \infty} F_n \in \mathbb{R}$ exists and let $\{A_{n}: n \geq 0\}$ be a sequence that satisfies the first order linear recurrence relation \eqref{RecurreneceLimitApendix}. Multiplying the recurrence relation with $n+2$ yields 
\begin{align*}
	(1+  n  )(2+  n  )A_{ n  }-(2+  n  )(3+  n  )A_{  n  +1}= F_{  n  }
\end{align*}
for all integers $n \geq 0$. Summing both sides of this identity, the left-hand side telescopes, giving
\begin{align*}
	2A_0-(n+1)(n+2)A_n=\sum_{k=0}^{n-1}(k+2)F_k,
\end{align*}
for all integers $n \geq 0$. This can be written as
\begin{align*}
	A_n=\frac{2A_0}{D_n}-\frac{S_n}{D_n},
\end{align*}
for all integers $n \geq 0$, where
\begin{align*}
	D_n:= (n+1)(n+2), \quad S_n:= \sum_{k=0}^{n-1}(k+2)F_k.
\end{align*}
Since $F_k \rightarrow F_{\infty}$, as $k \rightarrow +\infty$, we apply the Stolz-Ces\'aro theorem (the discrete analogue of L'Hopital's rule for limits of sequences), to obtain
\begin{align*}
	\lim_{n \rightarrow +\infty} A_n &=-\lim_{n \rightarrow +\infty}\frac{S_n}{D_n} = 
	-\lim_{n \rightarrow +\infty}\frac{S_{n+1}-S_n}{D_{n+1}-D_n}=
	-\lim_{n \rightarrow +\infty}\frac{(n+2)F_n}{2(n+2)}=-\frac{1}{2}F_{\infty}.
\end{align*}
which proves the claim. 
\end{proof}
\noindent

  \bibliographystyle{plain}
\bibliography{V6FourierEKG2}

@incollection {MR1448687,
    AUTHOR = {Paule, Peter and Riese, Axel},
     TITLE = {A {M}athematica {$q$}-analogue of {Z}eilberger's algorithm
              based on an algebraically motivated approach to
              {$q$}-hypergeometric telescoping},
 BOOKTITLE = {Special functions, {$q$}-series and related topics ({T}oronto,
              {ON}, 1995)},
    SERIES = {Fields Inst. Commun.},
    VOLUME = {14},
     PAGES = {179--210},
 PUBLISHER = {Amer. Math. Soc., Providence, RI},
      YEAR = {1997},
   MRCLASS = {33D15 (05A30 68Q40)},
  MRNUMBER = {1448687},
MRREVIEWER = {David M. Bressoud},
}

@article {MR2001619,
    AUTHOR = {Paule, Peter and Schneider, Carsten},
     TITLE = {Computer proofs of a new family of harmonic number identities},
   JOURNAL = {Adv. in Appl. Math.},
  FJOURNAL = {Advances in Applied Mathematics},
    VOLUME = {31},
      YEAR = {2003},
    NUMBER = {2},
     PAGES = {359--378},
      ISSN = {0196-8858},
   MRCLASS = {33F10 (11B65)},
  MRNUMBER = {2001619},
MRREVIEWER = {Scott Ahlgren},
       DOI = {10.1016/S0196-8858(03)00016-2},
       URL = {https://doi.org/10.1016/S0196-8858(03)00016-2},
}

@book {MR2768529,
    AUTHOR = {Kauers, Manuel and Paule, Peter},
     TITLE = {The concrete tetrahedron},
    SERIES = {Texts and Monographs in Symbolic Computation},
      NOTE = {Symbolic sums, recurrence equations, generating functions,
              asymptotic estimates},
 PUBLISHER = {SpringerWienNewYork, Vienna},
      YEAR = {2011},
     PAGES = {x+203},
      ISBN = {978-3-7091-0444-6},
   MRCLASS = {05-01 (05A15 05A16 33C05)},
  MRNUMBER = {2768529},
MRREVIEWER = {Takao Komatsu},
       DOI = {10.1007/978-3-7091-0445-3},
       URL = {https://doi.org/10.1007/978-3-7091-0445-3},
}

@article {MR1394550,
    AUTHOR = {Riese, Axel},
     TITLE = {A generalization of {G}osper's algorithm to bibasic
              hypergeometric summation},
   JOURNAL = {Electron. J. Combin.},
  FJOURNAL = {Electronic Journal of Combinatorics},
    VOLUME = {3},
      YEAR = {1996},
    NUMBER = {1},
     PAGES = {Research Paper 19, approx. 16},
   MRCLASS = {33D65 (33D20 68Q40)},
  MRNUMBER = {1394550},
       URL = {http://www.combinatorics.org/Volume_3/Abstracts/v3i1r19.html},
}

@article {MR1953199,
    AUTHOR = {Lyons, Russell and Paule, Peter and Riese, Axel},
     TITLE = {A computer proof of a series evaluation in terms of harmonic
              numbers},
   JOURNAL = {Appl. Algebra Engrg. Comm. Comput.},
  FJOURNAL = {Applicable Algebra in Engineering, Communication and
              Computing},
    VOLUME = {13},
      YEAR = {2002},
    NUMBER = {4},
     PAGES = {327--333},
      ISSN = {0938-1279},
   MRCLASS = {68W30},
  MRNUMBER = {1953199},
       DOI = {10.1007/s00200-002-0107-z},
       URL = {https://doi.org/10.1007/s00200-002-0107-z},
}

@article {MR1962799,
    AUTHOR = {Riese, Axel},
     TITLE = {q{M}ulti{S}um---a package for proving {$q$}-hypergeometric
              multiple summation identities},
   JOURNAL = {J. Symbolic Comput.},
  FJOURNAL = {Journal of Symbolic Computation},
    VOLUME = {35},
      YEAR = {2003},
    NUMBER = {3},
     PAGES = {349--376},
      ISSN = {0747-7171},
   MRCLASS = {33F10 (33D70 39A13)},
  MRNUMBER = {1962799},
MRREVIEWER = {Raimundas Vidunas},
       DOI = {10.1016/S0747-7171(02)00138-4},
       URL = {https://doi.org/10.1016/S0747-7171(02)00138-4},
}

@incollection {MR1395420,
    AUTHOR = {Paule, Peter and Schorn, Markus},
     TITLE = {A {M}athematica version of {Z}eilberger's algorithm for
              proving binomial coefficient identities},
      NOTE = {Symbolic computation in combinatorics $\Delta_1$ (Ithaca, NY,
              1993)},
   JOURNAL = {J. Symbolic Comput.},
  FJOURNAL = {Journal of Symbolic Computation},
    VOLUME = {20},
      YEAR = {1995},
    NUMBER = {5-6},
     PAGES = {673--698},
      ISSN = {0747-7171},
   MRCLASS = {39A10 (05A10 05A19 33C05 39-04)},
  MRNUMBER = {1395420},
MRREVIEWER = {Alain J. Phares},
       DOI = {10.1006/jsco.1995.1071},
       URL = {https://doi.org/10.1006/jsco.1995.1071},
}

@incollection {MR1880090,
    AUTHOR = {Riese, Axel},
     TITLE = {Fine-tuning {Z}eilberger's algorithm. {T}he methods of
              automatic filtering and creative substituting},
 BOOKTITLE = {Symbolic computation, number theory, special functions,
              physics and combinatorics ({G}ainesville, {FL}, 1999)},
    SERIES = {Dev. Math.},
    VOLUME = {4},
     PAGES = {243--254},
 PUBLISHER = {Kluwer Acad. Publ., Dordrecht},
      YEAR = {2001},
   MRCLASS = {33F10 (11B65)},
  MRNUMBER = {1880090},
       DOI = {10.1007/978-1-4613-0257-5\_15},
       URL = {https://doi.org/10.1007/978-1-4613-0257-5_15},
}

@article {MR485674,
    AUTHOR = {Gosper, Jr., R. William},
     TITLE = {Decision procedure for indefinite hypergeometric summation},
   JOURNAL = {Proc. Nat. Acad. Sci. U.S.A.},
  FJOURNAL = {Proceedings of the National Academy of Sciences of the United
              States of America},
    VOLUME = {75},
      YEAR = {1978},
    NUMBER = {1},
     PAGES = {40--42},
      ISSN = {0027-8424},
   MRCLASS = {10A40 (68A15)},
  MRNUMBER = {485674},
MRREVIEWER = {D. H. Lehmer},
       DOI = {10.1073/pnas.75.1.40},
       URL = {https://doi.org/10.1073/pnas.75.1.40},
}

@article {MR1048463,
    AUTHOR = {Zeilberger, Doron},
     TITLE = {A fast algorithm for proving terminating hypergeometric
              identities},
   JOURNAL = {Discrete Math.},
  FJOURNAL = {Discrete Mathematics},
    VOLUME = {80},
      YEAR = {1990},
    NUMBER = {2},
     PAGES = {207--211},
      ISSN = {0012-365X},
   MRCLASS = {33C20 (05A19 11B65 11Y16 33C05)},
  MRNUMBER = {1048463},
MRREVIEWER = {R. A. Askey},
       DOI = {10.1016/0012-365X(90)90120-7},
       URL = {https://doi.org/10.1016/0012-365X(90)90120-7},
}

@article {MR3270350,
    AUTHOR = {Craps, Ben and Evnin, Oleg and Vanhoof, Joris},
     TITLE = {Renormalization group, secular term resummation and {A}d{S}
              (in)stability},
   JOURNAL = {J. High Energy Phys.},
  FJOURNAL = {Journal of High Energy Physics},
      YEAR = {2014},
    NUMBER = {10},
     PAGES = {048, front matter+30},
      ISSN = {1126-6708},
   MRCLASS = {83C47},
  MRNUMBER = {3270350},
       DOI = {10.1007/JHEP10(2014)048},
       URL = {https://doi.org/10.1007/JHEP10(2014)048},
}

@book {MR1644447,
    AUTHOR = {Koepf, Wolfram},
     TITLE = {Hypergeometric summation},
    SERIES = {Advanced Lectures in Mathematics},
      NOTE = {An algorithmic approach to summation and special function
              identities},
 PUBLISHER = {Friedr. Vieweg \& Sohn, Braunschweig},
      YEAR = {1998},
     PAGES = {x+230},
      ISBN = {3-528-06950-3},
   MRCLASS = {33-02 (05A30 33C05 33C20 33F10)},
  MRNUMBER = {1644447},
MRREVIEWER = {Jiang Zeng},
       DOI = {10.1007/978-3-322-92918-1},
       URL = {https://doi.org/10.1007/978-3-322-92918-1},
}

@article{PhysRevD.92.084001,
  title = {Islands of stability and recurrence times in AdS},
  author = {Green, Stephen R. and Maillard, Antoine and Lehner, Luis and Liebling, Steven L.},
  journal = {Phys. Rev. D},
  volume = {92},
  issue = {8},
  pages = {084001},
  numpages = {23},
  year = {2015},
  month = {Oct},
  publisher = {American Physical Society},
  doi = {10.1103/PhysRevD.92.084001},
  url = {https://link.aps.org/doi/10.1103/PhysRevD.92.084001}
}

@article{PhysRevLett.113.071601,
  title = {Holographic Thermalization, Stability of Anti--de Sitter Space, and the Fermi-Pasta-Ulam Paradox},
  author = {Balasubramanian, Venkat and Buchel, Alex and Green, Stephen R. and Lehner, Luis and Liebling, Steven L.},
  journal = {Phys. Rev. Lett.},
  volume = {113},
  issue = {7},
  pages = {071601},
  numpages = {5},
  year = {2014},
  month = {Aug},
  publisher = {American Physical Society},
  doi = {10.1103/PhysRevLett.113.071601},
  url = {https://link.aps.org/doi/10.1103/PhysRevLett.113.071601}
}

@book {MR1379802,
    AUTHOR = {Petkov\v{s}ek, Marko and Wilf, Herbert S. and Zeilberger, Doron},
     TITLE = {{$A=B$}},
      NOTE = {With a foreword by Donald E. Knuth,
              With a separately available computer disk},
 PUBLISHER = {A K Peters, Ltd., Wellesley, MA},
      YEAR = {1996},
     PAGES = {xii+212},
      ISBN = {1-56881-063-6},
   MRCLASS = {05-01 (05A10 05A19 33C20 68R05)},
  MRNUMBER = {1379802},
MRREVIEWER = {Peter Paule},
}

@ARTICLE{SakisJacques,
       author = {{Chatzikaleas}, Athanasios and {Smulevici}, Jacques},
        title = "{Non-linear periodic waves on the Einstein cylinder}",
      journal = {arXiv e-prints},
         year = 2022,
        month = jan,
          eid = {arXiv:2201.05447},
        pages = {arXiv:2201.05447},
archivePrefix = {arXiv},
       eprint = {2201.05447},
 primaryClass = {math.AP},
       adsurl = {https://ui.adsabs.harvard.edu/abs/2022arXiv220105447C}
}

@book{MR3307944,
    AUTHOR = {Gradshteyn, I. S. and Ryzhik, I. M.},
     TITLE = {Table of integrals, series, and products},
   EDITION = {Eighth},
      NOTE = {Translated from the Russian,
              Translation edited and with a preface by Daniel Zwillinger and
              Victor Moll,
              Revised from the seventh edition [MR2360010]},
 PUBLISHER = {Elsevier/Academic Press, Amsterdam},
      YEAR = {2015},
     PAGES = {xlvi+1133},
      ISBN = {978-0-12-384933-5},
   MRCLASS = {00A22 (33-00)},
  MRNUMBER = {3307944},
}

@article {MR3527696,
    AUTHOR = {Cohl, Howard S.},
     TITLE = {Report from the open problems session at {OPSFA}13},
      NOTE = {Held in Gaithersburg, MD, June 4, 2015},
   JOURNAL = {SIGMA Symmetry Integrability Geom. Methods Appl.},
  FJOURNAL = {SIGMA. Symmetry, Integrability and Geometry. Methods and
              Applications},
    VOLUME = {12},
      YEAR = {2016},
     PAGES = {Paper No. 091, 12},
   MRCLASS = {33-02 (33C05 33C45 33C55 33F10 60B20)},
  MRNUMBER = {3527696},
       DOI = {10.3842/SIGMA.2016.071},
       URL = {https://doi-org.accesdistant.sorbonne-universite.fr/10.3842/SIGMA.2016.071},
}

@article {MR634149,
    AUTHOR = {Rahman, Mizan},
     TITLE = {A nonnegative representation of the linearization coefficients
              of the product of {J}acobi polynomials},
   JOURNAL = {Canadian J. Math.},
  FJOURNAL = {Canadian Journal of Mathematics. Journal Canadien de
              Math\'{e}matiques},
    VOLUME = {33},
      YEAR = {1981},
    NUMBER = {4},
     PAGES = {915--928},
      ISSN = {0008-414X},
   MRCLASS = {33A65},
  MRNUMBER = {634149},
MRREVIEWER = {R. A. Askey},
       DOI = {10.4153/CJM-1981-072-9},
       URL = {https://doi-org.accesdistant.sorbonne-universite.fr/10.4153/CJM-1981-072-9},
}

@article {MR3435559,
    AUTHOR = {Craps, Ben and Evnin, Oleg and Jai-akson, Puttarak and
              Vanhoof, Joris},
     TITLE = {Ultraviolet asymptotics for quasiperiodic {$\rm Ads_4$}
              perturbations},
   JOURNAL = {J. High Energy Phys.},
  FJOURNAL = {Journal of High Energy Physics},
      YEAR = {2015},
    NUMBER = {10},
     PAGES = {080, front matter+18},
      ISSN = {1126-6708},
   MRCLASS = {83C47},
  MRNUMBER = {3435559},
       DOI = {10.1007/JHEP10(2015)080},
       URL = {https://doi.org/10.1007/JHEP10(2015)080},
}

@article {MR3435560,
    AUTHOR = {Craps, Ben and Evnin, Oleg and Vanhoof, Joris},
     TITLE = {Ultraviolet asymptotics and singular dynamics of {A}d{S}
              perturbations},
   JOURNAL = {J. High Energy Phys.},
  FJOURNAL = {Journal of High Energy Physics},
      YEAR = {2015},
    NUMBER = {10},
     PAGES = {079, front matter+16},
      ISSN = {1126-6708},
   MRCLASS = {83C47},
  MRNUMBER = {3435560},
       DOI = {10.1007/JHEP10(2015)079},
       URL = {https://doi.org/10.1007/JHEP10(2015)079},
}

@book {MR2723248,
     TITLE = {N{IST} handbook of mathematical functions},
    EDITOR = {Olver, Frank W. J. and Lozier, Daniel W. and Boisvert, Ronald
              F. and Clark, Charles W.},
      NOTE = {With 1 CD-ROM (Windows, Macintosh and UNIX)},
 PUBLISHER = {U.S. Department of Commerce, National Institute of Standards
              and Technology, Washington, DC; Cambridge University Press,
              Cambridge},
      YEAR = {2010},
     PAGES = {xvi+951},
      ISBN = {978-0-521-14063-8},
   MRCLASS = {33-00 (00A20 65-00)},
  MRNUMBER = {2723248},
}

@Article{zbMATH07227949,
 Author = {Hunik-Kostyra, Dominika and Rostworowski, Andrzej},
 Title = {{AdS} instability: resonant system for gravitational perturbations of {{\(\mathrm{AdS}_5\)}} in the cohomogeneity-two biaxial {Bianchi} {IX} ansatz},
 FJournal = {Journal of High Energy Physics},
 Journal = {J. High Energy Phys.},
 ISSN = {1126-6708},
 Volume = {2020},
 Number = {6},
 Pages = {40},
 Note = {Id/No 2},
 Year = {2020},
 Language = {English},
 DOI = {10.1007/JHEP06(2020)002},
 zbMATH = {7227949},
 Zbl = {1437.83025}
}

@article {MR4150259,
    AUTHOR = {Moschidis, Georgios},
     TITLE = {A proof of the instability of {A}d{S} for the {E}instein-null
              dust system with an inner mirror},
   JOURNAL = {Anal. PDE},
  FJOURNAL = {Analysis \& PDE},
    VOLUME = {13},
      YEAR = {2020},
    NUMBER = {6},
     PAGES = {1671--1754},
      ISSN = {2157-5045},
   MRCLASS = {83C05 (35Q75 35Q83 83C57)},
  MRNUMBER = {4150259},
MRREVIEWER = {Kotik K. Lee},
       DOI = {10.2140/apde.2020.13.1671},
       URL = {https://doi-org.accesdistant.sorbonne-universite.fr/10.2140/apde.2020.13.1671},
}

@proceedings {MR1789139,
     TITLE = {Harmonic analysis and integral geometry},
    SERIES = {Chapman \& Hall/CRC Research Notes in Mathematics},
    VOLUME = {422},
 BOOKTITLE = {Proceedings of the session of the 1st {S}ummer {U}niversity of
              {S}afi held in {S}afi, {J}uly 15--20, 1998},
    EDITOR = {Picardello, Massimo A.},
 PUBLISHER = {Chapman \& Hall/CRC, Boca Raton, FL},
      YEAR = {2001},
     PAGES = {vi+171},
      ISBN = {1-58488-183-6},
   MRCLASS = {43-06 (58-06)},
  MRNUMBER = {1789139},
}

@ARTICLE{Moschidis2020NullDust,
    author  = {Moschidis, Georgios},
    title   = {A proof of the instability of {AdS}
               for the {Einstein}--null dust system
               with an inner mirror},
    journal = {Analysis \& PDE},
    volume  = {13},
    number  = {6},
    pages   = {1671--1754},
    year    = {2020},
    doi     = {10.2140/apde.2020.13.1671},
    eprint  = {1704.08681},
    archivePrefix = {arXiv},
    primaryClass  = {gr-qc}
}

@ARTICLE{2017arXiv170408685M,
       author = {{Moschidis}, Georgios},
        title = "{The Einstein--null dust system in spherical symmetry with an inner mirror: structure of the maximal development and Cauchy stability}",
      journal = {arXiv e-prints},
         year = 2017,
        month = apr,
          eid = {arXiv:1704.08685},
        pages = {arXiv:1704.08685},
archivePrefix = {arXiv},
       eprint = {1704.08685},
 primaryClass = {gr-qc},
       adsurl = {https://ui.adsabs.harvard.edu/abs/2017arXiv170408685M}
}

@article{PhysRevLett.111.051102,
  title = {Time-Periodic Solutions in an Einstein AdS--Massless-Scalar-Field System},
  author = {Maliborski, Maciej and Rostworowski, Andrzej},
  journal = {Phys. Rev. Lett.},
  volume = {111},
  issue = {5},
  pages = {051102},
  numpages = {5},
  year = {2013},
  month = {Aug},
  publisher = {American Physical Society},
  doi = {10.1103/PhysRevLett.111.051102},
  url = {https://link.aps.org/doi/10.1103/PhysRevLett.111.051102}
}

@article {MR3512262,
    AUTHOR = {Evnin, Oleg and Jai-akson, Puttarak},
     TITLE = {Detailed ultraviolet asymptotics for {A}d{S} scalar field
              perturbations},
   JOURNAL = {J. High Energy Phys.},
  FJOURNAL = {Journal of High Energy Physics},
      YEAR = {2016},
    NUMBER = {4},
     PAGES = {054, front matter+19},
      ISSN = {1126-6708},
   MRCLASS = {83C47},
  MRNUMBER = {3512262},
       DOI = {10.1007/jhep04(2016)054},
       URL = {https://doi.org/10.1007/jhep04(2016)054},
}

@ARTICLE{MoschidisNew,
       author = {{Moschidis}, Georgios},
        title = "{A proof of the instability of AdS for the Einstein-massless Vlasov system}",
      journal = {Invent. math.},
         year = 2022,
          eid = {},
        pages = {},
archivePrefix = {},
       eprint = {},
 primaryClass = {},
       adsurl = {https://doi.org/10.1007/s00222-022-01152-7}
}

@book {MR3828498,
    AUTHOR = {Moschidis, Georgios},
     TITLE = {Two {I}nstability {R}esults in {G}eneral {R}elativity},
      NOTE = {Thesis (Ph.D.)--Princeton University},
 PUBLISHER = {ProQuest LLC, Ann Arbor, MI},
      YEAR = {2018},
     PAGES = {286},
      ISBN = {978-0438-04763-1},
   MRCLASS = {Thesis},
  MRNUMBER = {3828498},
       URL =
              {http://gateway.proquest.com/openurl?url_ver=Z39.88-2004&rft_val_fmt=info:ofi/fmt:kev:mtx:dissertation&res_dat=xri:pqm&rft_dat=xri:pqdiss:10813054},
}

@article{PhysRevLett.111.041102,
  title = {Globally Regular Instability of 3-Dimensional Anti--De Sitter Spacetime},
  author = {Bizo\ifmmode \acute{n}\else \'{n}\fi{}, Piotr and Ja\l{}mu\ifmmode \dot{z}\else \.{z}\fi{}na, Joanna},
  journal = {Phys. Rev. Lett.},
  volume = {111},
  issue = {4},
  pages = {041102},
  numpages = {5},
  year = {2013},
  month = {Jul},
  publisher = {American Physical Society},
  doi = {10.1103/PhysRevLett.111.041102},
  url = {https://link.aps.org/doi/10.1103/PhysRevLett.111.041102}
}

@incollection {MR3329398,
    AUTHOR = {Bizo\'{n}, Piotr and Rostworowski, Andrzej},
     TITLE = {Instability of anti-de {S}itter spacetime},
 BOOKTITLE = {General relativity, cosmology and astrophysics},
    SERIES = {Fundam. Theor. Phys.},
    VOLUME = {177},
     PAGES = {239--244},
 PUBLISHER = {Springer, Cham},
      YEAR = {2014},
   MRCLASS = {83C25},
  MRNUMBER = {3329398},
}

@article {MR3491620,
    AUTHOR = {Craps, Ben and Evnin, Oleg},
     TITLE = {Ad{S} (in)stability: an analytic approach},
   JOURNAL = {Fortschr. Phys.},
  FJOURNAL = {Fortschritte der Physik. Progress of Physics},
    VOLUME = {64},
      YEAR = {2016},
    NUMBER = {4-5},
     PAGES = {336--344},
      ISSN = {0015-8208},
   MRCLASS = {83C47},
  MRNUMBER = {3491620},
       DOI = {10.1002/prop.201500067},
       URL = {https://doi.org/10.1002/prop.201500067},
}

@article{PhysRevLett.114.071102,
  title = {Stability of Anti--de Sitter Space in Einstein-Gauss-Bonnet Gravity},
  author = {Deppe, Nils and Kolly, Allison and Frey, Andrew and Kunstatter, Gabor},
  journal = {Phys. Rev. Lett.},
  volume = {114},
  issue = {7},
  pages = {071102},
  numpages = {5},
  year = {2015},
  month = {Feb},
  publisher = {American Physical Society},
  doi = {10.1103/PhysRevLett.114.071102},
  url = {https://link.aps.org/doi/10.1103/PhysRevLett.114.071102}
}

@article{PhysRevLett.115.081103,
  title = {Resonant Dynamics and the Instability of Anti--de Sitter Spacetime},
  author = {Bizo\ifmmode \acute{n}\else \'{n}\fi{}, Piotr and Maliborski, Maciej and Rostworowski, Andrzej},
  journal = {Phys. Rev. Lett.},
  volume = {115},
  issue = {8},
  pages = {081103},
  numpages = {5},
  year = {2015},
  month = {Aug},
  publisher = {American Physical Society},
  doi = {10.1103/PhysRevLett.115.081103},
  url = {https://link.aps.org/doi/10.1103/PhysRevLett.115.081103}
}

@article {MR3205859,
    AUTHOR = {Bizo\'{n}, Piotr},
     TITLE = {Is {A}d{S} stable?},
   JOURNAL = {Gen. Relativity Gravitation},
  FJOURNAL = {General Relativity and Gravitation},
    VOLUME = {46},
      YEAR = {2014},
    NUMBER = {5},
     PAGES = {Art. 1724, 14},
      ISSN = {0001-7701},
   MRCLASS = {83C75},
  MRNUMBER = {3205859},
       DOI = {10.1007/s10714-014-1724-0},
       URL = {https://doi.org/10.1007/s10714-014-1724-0},
}

@article {MR4016456,
    AUTHOR = {Biasi, Anxo and Craps, Ben and Evnin, Oleg},
     TITLE = {Energy returns in global {$\rm AdS_4$}},
   JOURNAL = {Phys. Rev. D},
  FJOURNAL = {Physical Review D},
    VOLUME = {100},
      YEAR = {2019},
    NUMBER = {2},
     PAGES = {024008, 17},
      ISSN = {2470-0010},
   MRCLASS = {83C05},
  MRNUMBER = {4016456},
       DOI = {10.1103/physrevd.100.024008},
       URL = {https://doi.org/10.1103/physrevd.100.024008},
}

@article {MR3708611,
    AUTHOR = {Bizo\'{n}, Piotr and Rostworowski, Andrzej},
     TITLE = {Gravitational turbulent instability of {$\rm AdS_5$}},
   JOURNAL = {Acta Phys. Polon. B},
  FJOURNAL = {Jagellonian University. Institute of Physics. Acta Physica
              Polonica B},
    VOLUME = {48},
      YEAR = {2017},
    NUMBER = {8},
     PAGES = {1375--1381},
      ISSN = {0587-4254},
   MRCLASS = {83C25},
  MRNUMBER = {3708611},
       DOI = {10.5506/APhysPolB.48.1375},
       URL = {https://doi.org/10.5506/APhysPolB.48.1375},
}

@article {MR2978943,
    AUTHOR = {Dias, \'{O}scar J. C. and Horowitz, Gary T. and Santos, Jorge E.},
     TITLE = {Gravitational turbulent instability of anti-de {S}itter space},
   JOURNAL = {Classical Quantum Gravity},
  FJOURNAL = {Classical and Quantum Gravity},
    VOLUME = {29},
      YEAR = {2012},
    NUMBER = {19},
     PAGES = {194002, 7},
      ISSN = {0264-9381},
   MRCLASS = {83C30},
  MRNUMBER = {2978943},
       DOI = {10.1088/0264-9381/29/19/194002},
       URL = {https://doi.org/10.1088/0264-9381/29/19/194002},
}

@article{PhysRevLett.119.191103,
  title = {Nonspherically Symmetric Collapse in Asymptotically AdS Spacetimes},
  author = {Bantilan, Hans and Figueras, Pau and Kunesch, Markus and Romatschke, Paul},
  journal = {Phys. Rev. Lett.},
  volume = {119},
  issue = {19},
  pages = {191103},
  numpages = {6},
  year = {2017},
  month = {Nov},
  publisher = {American Physical Society},
  doi = {10.1103/PhysRevLett.119.191103},
  url = {https://link.aps.org/doi/10.1103/PhysRevLett.119.191103}
}

@article{PhysRevLett.107.031102,
  title = {Weakly Turbulent Instability of Anti--de Sitter Spacetime},
  author = {Bizo\ifmmode \acute{n}\else \'{n}\fi{}, Piotr and Rostworowski, Andrzej},
  journal = {Phys. Rev. Lett.},
  volume = {107},
  issue = {3},
  pages = {031102},
  numpages = {4},
  year = {2011},
  month = {Jul},
  publisher = {American Physical Society},
  doi = {10.1103/PhysRevLett.107.031102},
  url = {https://link.aps.org/doi/10.1103/PhysRevLett.107.031102}
}

@ARTICLE{DafermosHolzegel,
       author = {{Dafermos}, M and {Holzegel}, G},
        title = { 
 Dynamic instability of solitons in $4+1$ dimesnional gravity with negative cosmological constant.},
      journal = {https://www.dpmms.cam.ac.uk/~md384/ADSinstability.pdf},
         year = 2006,
          eid = {},
        pages = {},
archivePrefix = {},
       eprint = {},
 primaryClass = {},
       adsurl = {}
}

@ARTICLE{DafermosTalk,
       author = {{Dafermos}, M},
        title = {The black hole stability problem.},
      journal = {Talk at the Newton Institute, Cambridge, http://www- old.newton.ac.uk/webseminars/pg+ws/2006/gmx/1010/dafermos/},
         year = 2006,
          eid = {},
        pages = {},
archivePrefix = {},
       eprint = {Talk given at the Newton Institute, Cambridge},
 primaryClass = {},
       adsurl = {}
}

@misc{chatzikaleas2023timeperiodicsolutionsnekhoroshev,
      title={Time periodic solutions and Nekhoroshev stability to non-linear massive Klein-Gordon equations in Anti-de Sitter}, 
      author={Athanasios Chatzikaleas and Jacques Smulevici},
      year={2023},
      eprint={2304.12784},
      archivePrefix={arXiv},
      primaryClass={math.AP},
      url={https://arxiv.org/abs/2304.12784}, 
}

@ARTICLE{2020arXiv200411049C,
    author  = {Chatzikaleas, Athanasios},
    title   = {On the Fourier Analysis of the
               {Einstein}--{Klein}--{Gordon} System:
               Growth and Decay of the Fourier Constants},
    journal = {Annales Henri Poincar{\'e}},
    year    = {2024},
    volume  = {25},
    number  = {6},
    pages   = {3009--3079},
    doi     = {10.1007/s00023-023-01393-z},
    eprint  = {2004.11049},
    archivePrefix = {arXiv},
    primaryClass  = {math.AP}
}

@misc{Mathematica,
  author = {Wolfram Research{,} Inc.},
  title = {Mathematica, {V}ersion 12.2},
  url = {https://www.wolfram.com/mathematica},
  note = {Champaign, IL, 2020}
}

@article {MathematicaFastZeilbergerPackage,
    AUTHOR = {Paule, Peter and Schorn, Markus and Riese, Axel},
     TITLE = {Fast Zeilberger Package, {V}ersion 3.61},
   JOURNAL = {Copyright Research Institute for Symbolic Computation (RISC), Johannes Kepler University, Linz, Austria},
}

@article {MathematicaMultiSumPackage,
    AUTHOR = {Wegschaider, Kurt and Riese, Axel and Zimmermann, Burkhard},
     TITLE = {MultiSum Package, {V}ersion 2.3},
   JOURNAL = {Copyright Research Institute for Symbolic Computation (RISC), Johannes Kepler University, Linz, Austria}, 
}

@book {MR0372517,
    AUTHOR = {Szeg\H{o}, G\'{a}bor},
     TITLE = {Orthogonal polynomials},
   EDITION = {Fourth},
      NOTE = {American Mathematical Society, Colloquium Publications, Vol.
              XXIII},
 PUBLISHER = {American Mathematical Society, Providence, R.I.},
      YEAR = {1975},
     PAGES = {xiii+432},
   MRCLASS = {42A52 (33A65)},
  MRNUMBER = {0372517},
}

\end{document}